\documentclass[a4paper,12pt]{report}
\usepackage{amsmath,amsthm,amssymb,color,epsfig,bm}
\usepackage[french]{babel}
\usepackage[T1]{fontenc}
\usepackage[latin1]{inputenc}
\usepackage[all]{xy}

\numberwithin{equation}{section}

{\theoremstyle{definition}
\newtheorem{definitiona}{D{\'e}finition}[section]

\newtheorem{remark}[definitiona]{Remarque}
\newtheorem{example}[definitiona]{Exemple}}
\newtheorem{proposition}[definitiona]{Proposition}
\newtheorem{lemma}[definitiona]{Lemme}
\newtheorem{theorem}[definitiona]{Th{\'e}or{\`e}me}
\newtheorem{corollary}[definitiona]{Corollaire}

\newenvironment{definition}{\begin{definitiona}}{\mbox{} \end{definitiona}}

\renewcommand{\theenumi}{\alph{enumi}}
\renewcommand{\labelenumi}{\rm\theenumi)}

\newcommand{\N}{\mathbb{N}}
\newcommand{\Z}{\mathbb{Z}}
\newcommand{\Q}{\mathbb{Q}}
\newcommand{\R}{\mathbb{R}}
\newcommand{\C}{\mathbb{C}}
\newcommand{\F}{\mathbb{F}}
\newcommand{\T}{\mathbb{T}}

\newcommand{\fc}{\mathbf{1}}
\newcommand{\cS}{\mathcal{S}}
\newcommand{\cP}{\mathcal{P}}
\newcommand{\cL}{\mathcal{L}}
\newcommand{\cG}{\mathcal{G}}
\newcommand{\cA}{\mathcal{A}}
\newcommand{\cM}{\mathcal{M}}
\newcommand{\cF}{\mathcal{F}}
\newcommand{\cH}{\mathcal{H}}
\newcommand{\cK}{\mathcal{K}}
\newcommand{\cB}{\mathcal{B}}
\newcommand{\gA}{\mathfrak{A}}
\newcommand{\gB}{\mathfrak{B}}

\def\k{C_c(G,A)}
\def\N{{\Bbb N}}
\def\B{{\Bbb B}}
\def\Z{{\Bbb Z}}
\def\Q{{\Bbb Q}}
\def\R{{\Bbb R}}
\def\C{{\Bbb C}}

\def\F{{\Bbb F}}
\def\K{{\Bbb K}}

\newcommand{\Aut}{{\rm Aut}}
\newcommand{\ie}{{\it i.e.}\/ }
\newcommand{\cf}{{\it cf.}\/ }

\begin{document}

\thispagestyle{empty}

\quad
\vspace{-1.5cm}


\centerline{\large{\textbf{UNIVERSIT{\'E} PARIS DIDEROT - PARIS 7}}}

\vspace{0.5cm}

\centerline{\large{UFR  Math{\'e}matiques}}

\vspace{0.25cm}
\vspace{0.75cm}

\centerline{\Large{\textbf{\textsc{Th{\`e}se de Doctorat}}}}

\vspace{0.5cm}

\centerline{\large{Discipline : Math{\'e}matiques}}

\vspace{2cm}
\centerline{\large{Pr{\'e}sent{\'e}e par}}

\vspace{0.5cm}

\centerline{\large{\textbf{ATHINA MAGEIRA}}}

\vspace{0.5cm}

\centerline{\large{pour obtenir le grade de Docteur de l'Universit{\'e}
    Paris 7}}
\vspace{2cm}
\vspace{1cm}

\centerline{\Large{\textbf{$\boldsymbol{C^*}$-\textsc{Alg{\`e}bres Gradu{\'e}es Par Un
        Semi-Treillis}}}}

\vspace{1cm}
\vspace{2.5cm}
\leftline{\normalsize{
\begin{tabular}{lll}
\normalsize{Th{\`e}se soutenue le 25 avril 2007 devant le Jury compos{\'e}
  de}\vspace{0.40cm}\\
\textbf{Mme. Anne BOUTET DE MONVEL}\vspace{0.20cm}\\
\textbf{M. Vladimir GEORGESCU} $\;$  Rapporteur\vspace{0.20cm}\\
\textbf{M. Andrei IFTIMOVICI}\vspace{0.20cm}\\
\textbf{M. Jean RENAULT}\vspace{0.20cm}\\
\textbf{M. Georges  SKANDALIS} $\;$  Directeur \vspace{0.20cm}\\
\textbf{M. Alain VALETTE} $\;$  Rapporteur\\
\end{tabular}
}}

\newpage
\pagenumbering{roman}\setcounter{page}{0}
\thispagestyle{empty}
${}$

\newpage
\pagenumbering{roman}\setcounter{page}{1}
\tableofcontents

\renewcommand{\labelitemi}{$\bullet $}
\everymath={\displaystyle}

\chapter*{Introduction}
\pagenumbering{arabic}
\addcontentsline{toc}{chapter}{Introduction}
Cette th{\`e}se est consacr{\'e}e {\`a} l'{\'e}tude des $C^*$-alg{\`e}bres gradu{\'e}es par un semi-treillis.

\bigskip Dans le cadre de leur travail sur le probl{\`e}me {\`a} $N$ corps, A. Boutet de Monvel et V. Georgescu ont {\'e}t{\'e} amen{\'e}s {\`a} introduire dans \cite{bg1}, \cite{bg4} et \cite{bg2} l'{\'e}tude des $C^*$-alg{\`e}bres gradu{\'e}es par un
semi-treillis fini. 

\bigskip L'utilisation des $C^*$-alg{\`e}bres dans le probl{\`e}me {\`a} $N$ corps quantique est assez r{\'e}cente (\cf \cite{abg} - voir aussi par exemple \cite{amp} pour d'autres utilisations r{\'e}centes).

D'autre part, R.G. Froese et I. Herbst ont introduit dans \cite{fh} la notion d'un semi-treillis toujours en relation avec le probl{\`e}me {\`a} $N$ corps. Cette notion {\'e}tait par la suite utilis{\'e}e et d{\'e}velopp{\'e}e par W.O. Amrein, A. Boutet de Monvel et V. Georgescu dans \cite{abg1}, \cite{abg2} et \cite{abg3} pour donner une description d{\'e}taill{\'e}e des propri{\'e}t{\'e}s spectrales des hamiltoniens
des syst{\`e}mes quantiques {\`a} plusieurs corps (en particulier ils ont
{\'e}tudi{\'e} une classe d'hamiltoniens \og de type $A$\fg\ qui appara{\^\i}t pour la
premi{\`e}re fois dans le livre de S. Agmon \cite{agm}).

\bigskip C'est donc tout naturellement que A. Boutet de Monvel et V. Georgescu
 ont {\'e}t{\'e} amen{\'e}s {\`a} introduire dans \cite{bg1}, \cite{bg4} et \cite{bg2}
 l'{\'e}tude des $C^*$-alg{\`e}bres gradu{\'e}es par un
semi-treillis fini. Les composantes (parties homog{\`e}nes) de ces alg{\`e}bres
correspondent aux \og niveaux d'int{\'e}raction\fg. L'utilisation de ces
$C^*$-alg{\`e}bres gradu{\'e}es leur a permis de retrouver des r{\'e}sultats classiques  sur le spectre essentiel des
hamiltoniens qui \og engendrent\fg\ ces alg{\`e}bres, comme par exemple le th{\'e}or{\`e}me de Hunziker-Van Winter-Zhislin (HVZ) et de
donner une g{\'e}n{\'e}ralisation de l'{\'e}quation de Weinberg-van Winter (WVW)
introduite dans les ann{\'e}es soixante (voir aussi dans \cite{rs2} pour
 une premi{\`e}re approche). De plus, dans \cite{bg5} et dans un cadre plus
 g{\'e}n{\'e}ral dans \cite{bg2} ils se sont servis des
 $C^*$-alg{\`e}bres gradu{\'e}es pour faire le calcul de
 l'estimation de Mourre pour des syst{\`e}mes {\`a} $N$ corps (\cf aussi \cite{bg3}). On retrouve les $C^*$-alg{\`e}bres gradu{\'e}es dans l'article de A. Boutet
de Monvel, V. Georgescu et A. Soffer \cite{bg6} pour l'{\'e}tude des
hamiltoniens d'un syst{\`e}me {\`a} $N$-corps avec des interactions tr{\`e}s
singuli{\`e}res. 

 On pourra consulter  \cite{abg} pour une pr{\'e}sentation plus globale et syst{\'e}matique des r{\'e}sultats rapidement cit{\'e}s ci-dessus.

\bigskip Les travaux de M. Damak et V. Georgescu dans \cite{dg1}, \cite{dg2} ainsi que V. Georgescu et  A. Iftimovici dans \cite{gi1}, \cite{gi2} et \cite{gi3} proposent un {\'e}largissement du cadre et une syst{\'e}matisation de l'{\'e}tude des $C^*$-alg{\`e}bres gradu{\'e}es. Dans ces travaux, les auteurs {\'e}tudient les $C^*$-alg{\`e}bres gradu{\'e}es par des semi-treillis pouvant {\^e}tre infinis. 

Un important r{\'e}sultat de cette s{\'e}rie d'articles consiste {\`a} reconna{\^\i}tre la $C^*$-alg{\`e}bre gradu{\'e}e consid{\`e}r{\'e}e dans \cite{bg1} comme un produit crois{\'e}. Ceci est d{\'e}montr{\'e} par M. Damak et V. Georgescu dans  \cite{dg1} en s'appuyant sur un r{\'e}sultat donn{\'e} par V. Georgescu et  A. Iftimovici (voir th{\'e}or{\`e}me 3.12 dans \cite{gi1}). Il s'agit du
 produit crois{\'e} de la $C^*$-alg{\`e}bre form{\'e}e des potentiels d'int{\'e}raction (qui est une $C^*$-alg{\`e}bre
 gradu{\'e}e et commutative) par le groupe des translations. Cette
 nouvelle forme a {\'e}t{\'e} tr{\`e}s utile pour d{\'e}terminer le quotient de cette $C^*$-alg{\`e}bre par
 une alg{\`e}bre d'op{\'e}rateurs compacts et par cons{\'e}quent
 pour obtenir des nouveaux r{\'e}sultats
 dans la th{\'e}orie spectrale des hamiltoniens d'un syst{\`e}me physique. Ces r{\'e}sultats sont d{\'e}velopp{\'e}s dans \cite{dg1}, \cite{gi1} et \cite{gi3}.

En dehors de ce produit crois{\'e}, un autre exemple de $C^*$-alg{\`e}bres
gradu{\'e}es par un semi-treillis, \emph{la $C^*$-alg{\`e}bre
  symplectique} (qui {\'e}tait d{\'e}j{\`a} apparue dans \cite{bg2}), a {\'e}t{\'e} {\'e}tudi{\'e} par V. Georgescu et A. Iftimovici dans \cite{gi2}. 

\bigskip \bigskip Dans cette th{\`e}se, nous proposons une {\'e}tude syst{\'e}matique des $C^*$-alg{\`e}bres gradu{\'e}es
par un semi-treillis quelconque.\begin{itemize}
\item  On simplifie quelques
axiomes des $C^*$-alg{\`e}bres gradu{\'e}es. 
\item On reconstruit ces alg{\`e}bres et on donne une pr{\'e}sentation alg{\'e}brique en fonction des composantes et de leur produit.
\item On {\'e}tablit la stabilit{\'e} des $C^*$-alg{\`e}bres gradu{\'e}es pour des op{\'e}rations comme le produit crois{\'e} et le produit tensoriel.
\item  On {\'e}tudie des propri{\'e}t{\'e}s classiques des
$C^*$-alg{\`e}bres pour les $C^*$-alg{\`e}bres gradu{\'e}es (commutativit{\'e},
nucl{\'e}arit{\'e}, exactitude) ainsi que leur $K$-th{\'e}orie.
\item On propose enfin l'{\'e}tude
de quelques exemples.
\end{itemize}

\bigskip \bigskip Nous pr{\'e}sentons maintenant nos r{\'e}sultats un peu plus en d{\'e}tail.

Soit $\gA$ une $C^*$-alg{\`e}bre. On dit qu'elle est \emph{gradu{\'e}e} par
  un semi-treillis $\cL$ ou bien qu'elle est \emph{$\cL$-gradu{\'e}e} si l'on s'est donn{\'e} une famille
  lin{\'e}airement ind{\'e}pendante et totale $(A_i)_{i\in \cL} $ de
  sous-$C^*$-alg{\`e}bres de $\gA$, que l'on appellera les \emph{composantes} de $A$  telles que $A_i A_j\subset
  A_{i\wedge j}$ pour tout $i,j\in \cal L$. 
  
  Soit $\gA$ une $C^*$-alg{\`e}bre gradu{\'e}e par un semi-treillis.  On dit que $(\gA,(A_i)_{i\in \cL})$ est une $C^*$-alg{\`e}bre $\cL$-gradu{\'e}e o{\`u} pour $i\in \cL$ on a not{\'e} $A_i$ la composante de $\gA$ correspondante.

\medskip On peut remarquer tout de suite qu'une $C^*$-alg{\`e}bre gradu{\'e}e par un semi-treillis $\cL$
est limite inductive de ses sous-alg{\`e}bres gradu{\'e}es par les
sous-semi-treillis finis de $\cL$. De plus, une $C^*$-alg{\`e}bre gradu{\'e}e pr{\'e}sente plusieurs suites exactes scind{\'e}es qui permettent - dans le cas d'un semi-treillis fini - de la comprendre inductivement. Dans un sens, on peut consid{\'e}rer qu'une graduation d'une $C^*$-alg{\`e}bre par un semi-treillis est une fa{\c c}on  d'organiser une famille de suites exactes scind{\'e}es \og compatibles entre elles\fg. Il s'ensuit que toute construction
 \og fonctorielle\fg\ de $C^*$-alg{\`e}bres qui est compatible avec les suites
 exactes scind{\'e}es et les limites inductives, va aussi {\^e}tre compatible
 avec les $C^*$-alg{\`e}bres gradu{\'e}es. Par exemple, si $(\gA,(A_i)_{i\in
   \cL})$ est une $C^*$-alg{\`e}bre gradu{\'e}e et si un groupe localement
 compact $G$ op{\`e}re sur $\gA$ en pr{\'e}servant les $A_i$, alors le produit
 crois{\'e} (plein ou r{\'e}duit) $\gA\rtimes G$ est gradu{\'e} par les
 $A_i\rtimes G$. Il en va de m{\^e}me pour les produits tensoriels de
 $C^*$-alg{\`e}bres. De plus, le m{\^e}me principe montre que la $K$-th{\'e}orie
 d'une $C^*$-alg{\`e}bre gradu{\'e}e est somme directe des $K$-th{\'e}ories de ses
 composantes homog{\`e}nes.

\medskip Soit $(\gA,(A_i)_{i\in \cL})$ une $C^*$-alg{\`e}bre gradu{\'e}e. La restriction du produit au niveau de ses composantes, fournit des applications $q_{i,j}:A_i\times A_j\to A_{i\wedge j}$ (pour $i,j\in \cL$) appel{\'e}es \emph{applications de structure}  et des morphismes de $C^*$-alg{\`e}bres $\varphi_{i,j}:A_j\to M(A_i)$
(avec $i,j\in \cL$ tels que $i\leq j$) appel{\'e}s \emph{morphismes de
  structure} - ici $M(A_i)$ d{\'e}signe l'alg{\`e}bre des multiplicateurs de
$A_i$. On {\'e}nonce imm{\'e}diatement les propri{\'e}t{\'e}s alg{\'e}briques de la
famille $(q_{i,j})_{i,j\in \cL}$ et de la famille
$(\varphi_{i,j})_{i\le j}$ qui traduisent l'associativit{\'e} du produit
et les propri{\'e}t{\'e}s de l'involution de la $C^*$-alg{\`e}bre $\gA$, et il est
alors facile de d{\'e}crire le passage de l'une {\`a} l'autre de ces
familles. Nous d{\'e}montrons en fait que ces applications et morphismes
de structure nous permettent de reconstruire les $C^*$-alg{\`e}bres gradu{\'e}es en ce sens qu'{\`a} une famille $(A_i)_{i\in \cL}$ de $C^*$-alg{\`e}bres et une famille d'applications  $(q_{i,j})_{i,j\in \cL}$ (ou $(\varphi_{i,j})_{i\le j}$) v{\'e}rifiant les propri{\'e}t{\'e}s mentionn{\'e}es ci-dessus correspond une et une seule $C^*$-alg{\`e}bre gradu{\'e}e ({\`a} isomorphisme pr{\`e}s) qui admet les $A_i$ comme  composantes et les $q_{i,j}$ comme applications de structure (ou les $\varphi_{i,j}$ comme morphismes de structure).

\medskip 
Puisqu'une $C^*$-alg{\`e}bre gradu{\'e}e est enti{\`e}rement d{\'e}crite par ses
composantes homog{\`e}nes, il est naturel d'essayer de lire des propri{\'e}t{\'e}s de cette alg{\`e}bre
uniquement en termes des dites composantes. On d{\'e}montre ainsi qu'une
$C^*$-alg{\`e}bre gradu{\'e}e est commutative, nucl{\'e}aire ou exacte si et
seulement si ses composantes v{\'e}rifient cette m{\^e}me
propri{\'e}t{\'e}. Dans le cas commutatif, on peut aussi faire un lien entre
le spectre d'une alg{\`e}bre gradu{\'e}e $\gA$ et ceux de ses
composantes. Dans de bons cas, le spectre des composantes d{\'e}crit
compl{\`e}tement celui de $\gA$.

\medskip 
Les repr{\'e}sentations d'une $C^*$-alg{\`e}bre gradu{\'e}e (ainsi que les
morphismes {\`a} valeurs dans une autre $C^*$-alg{\`e}bre) sont enti{\`e}rement d{\'e}crites en termes des repr{\'e}sentations de ses composantes. On peut alors donner les conditions nec{\'e}ssaires et suffisantes pour qu'un homomorphisme d'une $C^*$-alg{\`e}bre gradu{\'e}e {\`a} valeurs dans une autre $C^*$-alg{\`e}bre soit injectif ou surjectif en termes de ses restrictions aux composantes. Si le semi-treillis poss{\`e}de un plus petit {\'e}l{\'e}ment $i_0$ on {\'e}tudie plus particuli{\`e}rement  l'injectivit{\'e} de l'homomorphisme $\gA\to M(A_{i_0})$. Ces r{\'e}sultats
nous permettent d'une part de simplifier l'{\'e}tude du produit crois{\'e}
trait{\'e} par M. Damak
et V. Georgescu dans \cite{dg1} et d'autre part d'{\'e}tudier de nombreux 
exemples de $C^*$-alg{\`e}bres gradu{\'e}es.

Nous {\'e}tudions en particulier deux exemples commutatifs dont nous d{\'e}terminons le spectre.\begin{itemize}
\item  Pour le premier, on consid{\`e}re un sous-semi treillis $\cG$ du treillis (pour l'inclusion) de sous espaces d'un espace vectoriel de dimension finie $E$ et dont les composantes sont $A_L=C_0(E/L)$ (pour $L\in\cG$). Nous {\'e}tudions particuli{\`e}rement le cas o{\`u} $E$ est un plan $\cP$ et $\cG=\{\{0\},\delta_1,...,\delta_n,\cP\}$ o{\`u} $\delta_1,\ldots, \delta_n$ sont 
$n$ droites  de $\cP$. 
 En utilisant la th{\'e}orie des
 formes normales d'une surface de Riemann (voir \cite {ma} et
 \cite{spr}) on montre que le spectre de la $C^*$-alg{\`e}bre gradu{\'e}e associ{\'e}e est hom{\'e}omorphe au tore
${\mathcal T}_g$ {\`a} $g$ trous o{\`u} $g=\frac{n} {2}$ si $n$ est pair et dans le cas o{\`u}
$n$ est impair ce spectre est un tore {\`a} $g=E(\frac{n} {2})$ trous pinc{\'e} (\ie
 deux de ses points sont id{\'e}ntifi{\'e}s). 
 \item Dans le deuxi{\`e}me exemple, on
 consid{\`e}re un semi-treillis $\cL$ et on prend toutes les composantes $A_i$ {\'e}gales {\`a} $\C$. On identifie le spectre de $\cL$ avec l'ensemble des sous-semi treillis finissants (non vides) de $\cL$. En particulier, lorsque $\cL=\Q$, on montre alors que le spectre de la  $C^*$-alg{\`e}bre gradu{\'e}e dont les composantes sont les $A_i$ est en
 bijection avec l'ensemble $\R\textstyle\coprod\Q\textstyle\coprod\{-\infty\}$.
 \end{itemize}

\bigskip \bigskip Le texte de la th{\`e}se se d{\'e}compose de la mani{\`e}re suivante:

\begin{itemize}
\item Dans le premier chapitre, on rappelle les principales d{\'e}finitions et propri{\'e}t{\'e}s des $C^*$-alg{\`e}bres. On y rappelle en particulier les notions de multiplicateurs de $C^*$-alg{\`e}bres, de produits crois{\'e}s d'une $C^*$-alg{\`e}bre par une action d'un groupe localement compact et de produit tensoriel de $C^*$-alg{\`e}bres.

\item Les $C^*$-alg{\`e}bres gradu{\'e}es par un semi-treillis sont introduites dans le deuxi{\`e}me chapitre. On explore les morphismes de $C^*$-alg{\`e}bres gradu{\'e}es et on discute l'injectivit{\'e} et la surjectivit{\'e} de ces morphismes; on explicite la relation entre les $C^*$-alg{\`e}bres gradu{\'e}es et les suites exactes scind{\'e}es de  $C^*$-alg{\`e}bres. Ensuite, on introduit les applications de structure et les morphismes de structure d'une $C^*$-alg{\`e}bre gradu{\'e}e et on {\'e}tudie les propri{\'e}t{\'e}s alg{\'e}briques de ces familles d'applications. On conclut avec un cas particulier de morphisme de $C^*$-alg{\`e}bres gradu{\'e}es dont on traite l'injectivit{\'e} qui nous sera utile dans l'{\'e}tude d'exemples.

\item Dans le troisi{\`e}me chapitre on montre comment une $C^*$-alg{\`e}bre gradu{\'e}e peut {\^e}tre enti{\`e}rement  reconstruite {\`a} partir de ses composantes et des morphismes de structure.

\item Dans le chapitre 4 on montre que le produit crois{\'e} de $C^*$-alg{\`e}bres et le produit tensoriel de $C^*$-alg{\`e}bres ont un bon comportement pour les $C^*$-alg{\`e}bres gradu{\'e}es.

\item Au chapitre 5 on {\'e}tudie la commutativit{\'e}, nucl{\'e}arit{\'e} et
exactitude d'une $C^*$-alg{\`e}bre gradu{\'e}e et on exprime ses groupes de K-th{\'e}orie en termes de ceux des composantes. Enfin, on {\'e}tudie le spectre des $C^*$-alg{\`e}bres gradu{\'e}es commutatives.

\item Enfin le chapitre 6 est consacr{\'e} {\`a} l'{\'e}tude de quelques exemples de
$C^*$-alg{\`e}bres gradu{\'e}es.  
\end{itemize} 
\newpage
\quad \thispagestyle{empty}

\newpage

\chapter{Pr{\'e}liminaires, rappels}

\section{$C^*$-alg{\`e}bres}

On rappelle ici bri{\`e}vement un certain nombre de d{\'e}finitions et
propri{\'e}t{\'e}s des $C^*$-alg{\`e}bres qui seront utiles dans la suite. Nous ne
donnons pas de d{\'e}monstrations. Celles-ci peuvent {\^e}tre trouv{\'e}es dans
\cite{mur}, \cite{pe} et \cite{wo}.

\begin{definition} Soit $A$ une alg{\`e}bre complexe. Une
  \emph{involution} de $A$ est une application antilin{\'e}aire  $a\mapsto
  a^*$ de $A$ dans $A$ telle que pour tout $a,b\in A$ on ait
  ${(a^*)}^*=a$ et ${(ab)}^*=b^*a^*$.

 On appelle le couple $(A,*)$ une alg{\`e}bre
  \emph{involutive} ou une \emph{$*$-alg{\`e}bre}.

 Une \emph{alg{\`e}bre de Banach involutive} est une alg{\`e}bre involutive
  $A$ munie d'une norme
  sous-multiplicative qui est compl{\`e}te, dont l'involution est
  isom{\'e}trique.

 Une \emph{$C^*$-alg{\`e}bre} $A$ est
  une alg{\`e}bre de Banach involutive telle que pour tout $a\in A$ on ait
  $||a^*a||={||a||}^2$.

Plus g{\'e}n{\'e}ralement, une $C^*$-(semi)-norme sur une alg{\`e}bre involutive
$A$ est une (semi)-norme sous-multiplicative $N:A\to \R_+$
satisfaisant l'{\'e}galit{\'e} $N(a^*a)=N(a)^2$ pour tout $a\in A$. On v{\'e}rifie
alors que l'on a $N(a^*)=N(a)$. En ce sens, une $C^*$-alg{\`e}bre est une
alg{\`e}bre de Banach involutive dont la norme est une $C^*$-norme. 
\end{definition}

\begin{definition}Soit $A$ une alg{\`e}bre complexe unif{\`e}re et $a\in A$. On note
  $A^{-1}$ l'ensemble des {\'e}l{\'e}ments inversibles de $A$. Le
  \emph{spectre} de $a$ dans $A$ est le sous-ensemble
  $\mbox{Sp}_A(a)=\{\lambda\in \C\,;\,a-\lambda\notin A^{-1}\}$ de $\C$.
\end{definition}

Si $A$ est une alg{\`e}bre non unif{\`e}re, on peut la plonger dans une
alg{\`e}bre unif{\`e}re $\widetilde A$ qui comme espace vectoriel est
isomorphe {\`a} $A\times \C$ et dont la loi du produit est d{\'e}finie par:
$(a,\lambda)(b,\mu)=(ab+\lambda a+\mu b,\lambda\mu)$ (pour tout
$a,b\in A,\lambda,\mu\in \C$). C'est une alg{\`e}bre unif{\`e}re (d'unit{\'e}
$(0,1)$) qui contient une copie de $A$, $A\times \{0\}$. On identifie
alors $A$ avec son image dans $\widetilde A$. Pour $a\in A$ et
$\lambda\in \C$ on note
$a+\lambda$ l'{\'e}l{\'e}ment $(a,\lambda)$ de $\widetilde A$. 

Si $A$ est une alg{\`e}bre de Banach qui n'est pas unif{\`e}re, l'alg{\`e}bre
unif{\`e}re $\widetilde A$ admet une structure d'alg{\`e}bre de Banach pour la
norme d{\'e}finie par: $||(a,\lambda)||=||a||+|\lambda|$ (pour $a\in
A,\lambda\in \C$).

\begin{remark}\label{unif}\begin{enumerate}
\item Si $A$ est une alg{\`e}bre non unif{\`e}re, pour tout $a\in A$, $0\in
  \mbox{Sp}_{\widetilde A}(a)$.

\item Si l'alg{\`e}bre $A$ poss{\`e}de d{\'e}j{\`a} un {\'e}l{\'e}ment unit{\'e} not{\'e} $e$,
  l'application $(a,\lambda)\mapsto (a+\lambda e,\lambda)$ est un
  isomorphisme entre $\widetilde A$ et l'alg{\`e}bre produit $A\times
  \C$. Ceci montre que $\mbox{Sp}_{\widetilde A}(a)=\mbox{Sp}_A(a)\cup
  \{0\}$ pour tout $a\in A$.

\item Si $\pi:A\to B$ est un homomorphisme d'alg{\`e}bres, il existe un
  unique homomorphisme unital $\widetilde \pi:\widetilde A\to
  \widetilde B$ d{\'e}fini par $\widetilde \pi(a+\lambda)=\pi(a)+\lambda$
  dont $\pi$ est la restriction. Il est injectif (resp. surjectif) si
  $\pi$ l'est. 

\end{enumerate}
\end{remark}

\begin{proposition}Soit $A$ une $C^*$-alg{\`e}bre (non n{\'e}cessairement
  unif{\`e}re). Munie de l'involution $(a+\lambda)^*=a^*+\lambda$,
  l'alg{\`e}bre $\widetilde A$ est une $C^*$-alg{\`e}bre. En d'autres termes
  $\widetilde A$ admet une (n{\'e}cessairement unique) norme qui en fait
  une $C^*$-alg{\`e}bre.\newline\indent \hfill$\square$\end{proposition}

 Donc l'alg{\`e}bre $A$ est une
  sous-alg{\`e}bre involutive ferm{\'e}e de $\widetilde A$.

\begin{definition} Soit $A$ une $C^*$-alg{\`e}bre et $B$ un sous-ensemble de $A$, on
  pose $B^*=\{b^*\,|\,b\in B\}$, on dit que $B$ est \emph{autoadjoint} si
  $B^*=B$.

 Un {\'e}l{\'e}ment $a\in A$ est \emph{autoadjoint} si $a=a^*$ et il
  est \emph{normal} si $a^*a=aa^*$.

  Si $A$ poss{\`e}de un {\'e}l{\'e}ment unit{\'e}, on dit que $a\in A$  est \emph{unitaire} si
  $a^*a=aa^*=1$.

 \end{definition}

\begin{definition} Soient $A$, $B$ des alg{\`e}bres involutives. Un
  homomorphisme $\varphi:A\to B$ est appel{\'e} un
  \emph{$*$-homomorphisme} ou un homomorphisme \emph{involutif} s'il est stable par l'involution,
  c'est-{\`a}-dire si $\varphi(a^*)={(\varphi(a))}^*$ pour tout $a\in A$. Si $\varphi:A\to B$ est un
  $*$-homomorphisme, alors $\varphi(A)$ est une sous-$*$-alg{\`e}bre de
  $B$. Si $A$ et $B$ sont des $C^*$-alg{\`e}bres on dit aussi que
  $\varphi$ est un \emph{morphisme} de $C^*$-alg{\`e}bres. 
\end{definition}

\begin{example} Soit $\cH$ un espace de Hilbert. On note $\B(\cH)$ l'ensemble des applications lin{\'e}aires born{\'e}es de $\cH$
dans lui-m{\^e}me (op{\'e}rateurs sur $\cH$). L'alg{\`e}bre $\B(\cH)$ est une
$C^*$-alg{\`e}bre.

 Soit $A$ une $C^*$-alg{\`e}bre. Un morphisme de
$C^*$-alg{\`e}bres $A\to \B(\cH)$ s'appelle une \emph{repr{\'e}sentation} de
$A$ dans $\cH$.
\end{example}

\paragraph{Th{\'e}or{\`e}me de Gel'fand.} Soit $A$ une alg{\`e}bre de Banach. On appelle
  \emph{caract{\`e}re} de $A$ tout homomorphisme d'alg{\`e}bres continu et non
  nul de $A$ dans $\C$.

On appelle \emph{spectre} d'une alg{\`e}bre de Banach commutative $A$ et
on note $\mbox{Sp}(A)$ l'ensemble de ses caract{\`e}res. On munit $\mbox{Sp}(A)$ de la
topologie de la convergence simple. Autrement dit la topologie de $\mbox{Sp}(A)$ est la topologie la moins fine pour laquelle les applications
$\chi\mapsto \chi(x)$ de $\mbox{Sp}(A)$ dans $\C$ sont continues (pour tout
$x\in A$). 

\begin{proposition}(Transformation de Gel'fand). Soit $A$ une
  alg{\`e}bre de Banach commutative.
\begin{enumerate}
\item Pour tout $a\in A$ on a
  $\mbox{Sp}_A(a)=\{\chi(a)\,\mbox{o{\`u}}\;\chi\in \mbox{Sp}(A)\}$.
\item Si $A$ est unif{\`e}re, $\mbox{Sp}(A)$ est compact.
\item L'application $G:A\to C(\mbox{Sp}(A))$ donn{\'e}e par $G(a):\chi\to
  \chi(a)$ est un homomorphisme d'alg{\`e}bres de Banach qui est continu.\newline\indent \hfill$\square$
\end{enumerate}

\end{proposition}

\begin{theorem}(Gel'fand). Soit $A$ une $C^*$-alg{\`e}bre commutative et
  unif{\`e}re.  La transformation de Gel'fand $G:A\to C(\mbox{Sp}(A))$ de $A$ est
  un isomorphisme isom{\'e}trique de $C^*$-alg{\`e}bres.\newline\indent \hfill$\square$ \end{theorem}

\begin{remark}\begin{enumerate}
\item Soit $A$ une alg{\`e}bre de Banach non unif{\`e}re et $\chi\in
 \mbox{Sp}(A)$. On peut {\'e}tendre $\chi$ de
 mani{\`e}re unique {\`a} un caract{\`e}re
 $\widetilde \chi$ de $\widetilde A$ en posant $\widetilde \chi
  (a+\lambda)=\chi(a)+\lambda$ (pour tout $a\in A,\lambda\in\C$). Ceci
 nous permet de voir que le spectre de $A$ est la partie du spectre de
 $\widetilde A$ form{\'e}e des caract{\`e}res qui ne sont pas nuls sur $A$. 

\item Soit $A$ une alg{\`e}bre de Banach commutative non n{\'e}cessairement
  unif{\`e}re. Le spectre
  de $A$ est un espace localement compact dont le compactifi{\'e}
  d'Alexandroff est le spectre de $\widetilde A$. La transformation de
  Gel'fand de $A$ est l'homomorphisme $a\mapsto G(a)$ de $A$ sur
  $C_0(\mbox{Sp}(A))$ qui est donn{\'e} par $G(a)(\chi)=\chi(a)$ pour tout
  $a\in A$ et $\chi\in \mbox{Sp}(A)$. Si $A$ est une $C^*$-alg{\`e}bre
  commutative, le th{\'e}or{\`e}me de Gel'fand correspondant nous dit que
  l'application $G$ est un isomorphisme
  isom{\'e}trique de $C^*$-alg{\`e}bres.
\end{enumerate}
\end{remark} 

\paragraph{Calcul fonctionnel continu.} Un cas particulier du th{\'e}or{\`e}me de
Gel'fand est le suivant: soit $x$ un {\'e}l{\'e}ment normal d'une
$C^*$-alg{\`e}bre unif{\`e}re $A$. Notons $B$ l'adh{\'e}rence dans $A$ de
l'ensemble $\{P(x,x^*),\,P\in \C[X,Y]\}$. C'est une sous-$C^*$-alg{\`e}bre
commutative de $A$ contenant l'unit{\'e}, donc isomorphe {\`a}
$C(\mbox{Sp}(B))$. Son spectre $\mbox{Sp}(B)$ s'id{\'e}ntifie {\`a}
$\mbox{Sp}_A(x)$ via l'hom{\'e}omorphisme $\chi\mapsto \chi(x)$. La
transformation de Gel'fand de $B$ est un isomorphisme isom{\'e}trique
$G:B\to C(\mbox{Sp}_A(x))$ qui associe {\`a} $x$ la fonction $z$ qui d{\'e}signe
l'inclusion de $\mbox{Sp}_A(x)$ dans $\C$. Pour $f\in
C(\mbox{Sp}_A(x))$ on pose $f(x)=G^{-1}(f)$.

\begin{proposition}\begin{enumerate}
\item Soient $x$ un {\'e}l{\'e}ment normal d'une $C^*$-alg{\`e}bre unif{\`e}re $A$ et
  $f\in C(\mbox{Sp}_A(x))$. On a $\mbox{Sp}_A
  f(x)=f(\mbox{Sp}_A(x))$, l'{\'e}l{\'e}ment $f(x)$ de $A$ est normal et pour
  toute fonction $g\in C(\mbox{Sp}_A f(x))$ on a $(g\circ
  f)(x)=g(f(x))$.
\item Soient $\pi:A\to B$ un morphime unital de $C^*$-alg{\`e}bres
  (unif{\`e}res), $x\in A$ un {\'e}l{\'e}ment normal de $A$ et $f$ une fonction
  continue sur $\mbox{Sp}_A(x)$. Alors $f(\pi(x))=\pi(f(x))$.
\newline\indent \hfill$\square$
\end{enumerate} 
\end{proposition}  

Soient $A$ une $C^*$-alg{\`e}bre non unif{\`e}re, $x$ un {\'e}l{\'e}ment normal de $A$ et $f\in
C(\mbox{Sp}_A(x))$, alors $f(x)\in \widetilde A$. A l'aide de la
proposition prec{\'e}dente (partie $b$) appliqu{\'e}e au morphisme
$\epsilon:(a,\lambda)\mapsto \lambda$ de $\widetilde A$ dans $\C$ on a $f(x)\in A$ si
et seulement si $f(0)=0$.

\paragraph{El{\'e}ments positifs; la relation d'ordre d'une
  $C^*$-alg{\`e}bre.}Un {\'e}l{\'e}ment $a\in A$ est \emph{positif} s'il est
  autoadjoint et son spectre $\mbox{Sp}_A(a)$ est un sous-ensemble de 
  ${\R}^{+}$. 
\begin{proposition}Soit $A$ une $C^*$-alg{\`e}bre.
\begin{enumerate}
\item Pour un {\'e}l{\'e}ment autoadjoint $h$ de $A$ les conditions suivantes
  sont {\'e}quivalentes:
\begin{enumerate}
\item $\mbox{Sp}_A(h)\subset \R^{+}.$
\item Il existe $k\in A$ tel que $k=k^*$ et $k^2=h$.
\item Il existe $a\in A$ tel que $a^*a=h$.
\end{enumerate}

\item Les {\'e}l{\'e}ments autoadjoints de $A$ v{\'e}rifiant ces conditions
  forment un c{\^o}ne convexe saillant de $A$.\newline\indent \hfill$\square$
\end{enumerate}
\end{proposition}

 Le c{\^o}ne des {\'e}l{\'e}ments positifs de $A$ est not{\'e}
  $A_{+}$. On d{\'e}finit une relation d'ordre sur l'ensemble des {\'e}l{\'e}ments
  autoadjoints de $A$ par $b-a\in A_{+}$ que l'on note $a\leq b$.\\
 
\paragraph{Unit{\'e}s approch{\'e}es.} On appelle \emph{unit{\'e} approch{\'e}e} d'une alg{\`e}bre de
  Banach $A$ une famille $(u_i)_{i\in I}$ d'{\'e}l{\'e}ments de $A$ ind{\'e}x{\'e}e
  par un ensemble $I$ muni d'un ordre filtrant croissant telle que,
  pour tout $a\in A$ on ait $a=\lim au_i$ et $a=\lim u_i a$.

Si $A$ est une $C^*$-alg{\`e}bre, on dit qu'une unit{\'e} approch{\'e}e
$(u_i)_{i\in I}$ de $A$ est croissante si pour
$i\leq i'$ on a $u_i\leq u_{i'}$.

\begin{proposition}Toute $C^*$-alg{\`e}bre admet une unit{\'e} approch{\'e}e croissante.
\newline\indent \hfill$\square$
\end{proposition}

\paragraph{Morphismes de $C^*$-alg{\`e}bres.} Un morphisme de
$C^*$-alg{\`e}bres est continu, de norme inf{\'e}rieure ou {\'e}gale {\`a}
$1$. Autrement dit, on a la proposition suivante:
\begin{proposition}Tout morphisme de $C^*$-alg{\`e}bres est contractant.\newline\indent \hfill$\square$
\end{proposition}

\begin{proposition} \label{isom} Soit $\varphi:A\to B$ un morphisme de
$C^*$-alg{\`e}bres, alors pour tout $y\in \varphi(A)$ il existe $x\in A$
tel que $\varphi(x)=y$ et $||x||=||y||$.\end{proposition}

\begin{proof}On peut supposer que $A$ et $B$ sont unif{\`e}res. Soit $y\in \varphi(A)$, il existe $z\in A$ tel que
  $\varphi(z)=y$. Posons $f(s)=\inf(1,{\frac{||y||}{\sqrt s}})$ pour 
  $s\geq 0$ ($f(0)=1$) et $x=zf(z^*z)$. Pour $s\in \mbox{Sp}(y^*y)$ on a $s\leq ||y^*y||=||y||^2$ donc
  $f(s)=1$ de sorte que $f(y^*y)=1$. Puisque $\varphi$ est un
  morphisme de $C^*$-alg{\`e}bres on a $\varphi(x)=yf(y^*y)=y$. Par la
  proposition prec{\'e}dente $\varphi$ est contractant donc
  $||y||=||\varphi(x)||\leq ||x||$. Par
  ailleurs $x^*x=g(z^*z)$ o{\`u} $g(s)=sf(s)^2=\inf(s,||y||^2)$, donc
  $||x||^2=||x^*x||\leq ||y||^2$ d'o{\`u} l'{\'e}galit{\'e} $||x||=||y||$.
 \end{proof}

On en deduit un r{\'e}sultat classique de $C^*$-alg{\`e}bres qui sera tr{\`e}s utile
dans la suite.

\begin{proposition}\label{coniso} Tout morphisme injectif de
  $C^*$-alg{\`e}bres est isom{\'e}trique. Tout morphisme de $C^*$-alg{\`e}bres est d'image ferm{\'e}e. 
\newline\indent \hfill$\square$\end{proposition}

\paragraph{Id{\'e}aux et quotients d'une $C^*$-alg{\`e}bre.}
Un id{\'e}al {\`a} gauche (resp. {\`a} droite) $I$, d'une alg{\`e}bre $A$ est un sous-espace vectoriel
de $A$ tel que $a\in A$ et $b\in I\Rightarrow ab\in I$ (resp. $ba\in
I$). Un id{\'e}al bilat{\`e}re (on l'appellera simplement \emph{id{\'e}al}) est un
id{\'e}al {\`a} gauche et {\`a} droite.

\medskip

Soit $I$ un id{\'e}al d'une alg{\`e}bre $A$, alors $A/I$ est une
alg{\`e}bre munie de la multiplication $(a+I)(b+I)=ab+I$. Soit $I$ un id{\'e}al ferm{\'e} d'une alg{\`e}bre norm{\'e}e $A$, alors $A/I$ est une alg{\`e}bre
munie de la norme \emph{quotient}: $||a+I||=\inf\limits_{b\in
  I}||a+b||$.

\medskip

Soient $A$, $B$ des alg{\`e}bres involutives. Remarquons que le noyau d'un homomorphisme involutif $\varphi:A\to B$
est un id{\'e}al autoadjoint de $A$.
 
\begin{proposition} Soit $I$ un id{\'e}al ferm{\'e} d'une $C^*$-alg{\`e}bre $A$, alors
$I$ est autoadjoint et le quotient $A/I$ est une $C^*$-alg{\`e}bre pour la
norme quotient.\newline\indent \hfill$\square$
\end{proposition}

\begin{proposition}\label{idssa}Soit $I$ un id{\'e}al ferm{\'e} d'une
  $C^*$-alg{\`e}bre $A$ et $B$ une sous-$C^*$-alg{\`e}bre de $A$. Alors $I+B$
  est une $C^*$-alg{\`e}bre.\newline\indent \hfill$\square$
\end{proposition}

 Un id{\'e}al $I$ d'une $C^*$-alg{\`e}bre $A$ est \emph{essentiel} si pour tout
autre id{\'e}al non-nul $J$ de $A$ on a $I\cap J\neq \{0\}$.

\paragraph{Doubles centralisateurs; multiplicateurs.} Soit $A$ une
$C^*$-alg{\`e}bre.\\ On appelle \emph{double
  centralisateur} de $A$ un couple $(L,R)$
d'applications $L,R:A\to A$ qui verifie $$R(a)b=aL(b)$$ pour tout
$a,b\in A$. On note ${\cal {DC}}(A)$ l'ensemble des doubles
centralisateurs de $A$.

\begin{proposition} Si $(L,R)$ est un double centralisateur de $A$
  alors $L(ab)=L(a)b$ et $R(ab)=aR(b)$ pour tout $a,b\in A$. Les applications $L$ et
  $R$ sont lin{\'e}aires born{\'e}s avec $||L||=||R||$.\newline\indent
  \hfill$\square$\end{proposition}

  L'ensemble ${\cal {DC}}(A)$ est une alg{\`e}bre de Banach munie de la norme
  $||(L,R)||=||L||=||R||$ et les op{\'e}rations suivantes :
  $$(L_1,R_1)+(L_2,R_2)=(L_1+L_2,R_1+R_2),$$
  $$z(L,R)=(zL,zR)\,(\mbox{pour}\,z\in \C),$$
  $$(L_1,R_1)(L_2,R_2)=(L_1L_2,R_2R_1).$$ 
Si $L:A\to A$ est une application lin{\'e}aire, on d{\'e}finit $L^*:A\to A$ par $L^*(a)={(L(a^*))}^*$. Alors
$L^*$ est lin{\'e}aire et l'application $L\to L^*$ est une application
antilin{\'e}aire isometrique de $\B(A)$ dans lui-m{\^e}me telle que
$L^{**}=L$ et ${(L_1L_2)}^*=L_1^*L_2^*.$ Si $(L,R)\in {\cal {DC}}(A)$
alors ${(L,R)}^*=(R^*,L^*)\in {\cal {DC}}(A)$. L'application $(L,R)\mapsto
{(L,R)}^*$ est une involution sur ${\cal {DC}}(A)$.\\ On a la proposition suivante:

\begin{proposition} Si $A$ est une $C^*$-alg{\`e}bre, alors ${\cal
    {DC}}(A)$ est une $C^*$-alg{\`e}bre.\newline\indent \hfill$\square$\end{proposition}

\begin{example} Soit $A$ un id{\'e}al ferm{\'e} d'une $C^*$-alg{\`e}bre $B$ et soit
 $b\in B$, alors $L_b:a\mapsto ba$ et $R_b:a\mapsto ab$ sont des applications
 de $A$ dans $A$ et le couple $(L_b,R_b)\in {\cal {DC}}(A)$. On
 remarque que pour tout $b\in B$, $||(L_b,R_b)||\leq ||b||$.
\end{example}

Si $A$ est une $C^*$-alg{\`e}bre, l'application $A\to {\cal {DC}}(A)$
donn{\'e}e par $a\mapsto (L_a,R_a)$ est un morphisme de $C^*$-alg{\`e}bres
isom{\'e}trique appel{\'e} le plongement canonique de $A$ dans ${\cal {DC}}(A)$. Donc, on
id{\'e}ntifie $A$ avec son image dans ${\cal {DC}}(A)$ qui est un id{\'e}al
ferm{\'e}. 
\begin{proposition} \label{injmul}Soit $A$ un id{\'e}al ferm{\'e} d'une $C^*$-alg{\`e}bre
 $B$. L'application $\mu:b\mapsto (L_b,R_b)$ est un homomorphisme
  $\mu:B\to {\cal {DC}}(A)$ dont la restriction {\`a} $A$ est le
  plongement canonique de $A$ dans ${\cal {DC}}(A)$. De plus, $\mu$ est injectif si et seulement si
  $A$ est essentiel dans $B$.\newline\indent \hfill$\square$ 
\end{proposition}

Soit $A$ et $B$ des $C^*$-alg{\`e}bres. Un morphisme $\pi:A\to
{\cal {DC}}(B)$ est \emph{non-d{\'e}g{\'e}n{\'e}r{\'e}} si $\pi(A)B=\{\pi(a)b$ o{\`u} $a\in A$,
$b\in B\}$ est dense dans $B$.

\begin{proposition}\label{cohen}Soit $A$ un id{\'e}al ferm{\'e} d'une $C^*$-alg{\`e}bre $D$ et
  $\pi:A\to {\cal {DC}}(B)$ un morphisme non-d{\'e}g{\'e}n{\'e}r{\'e}. Il existe une
    unique extension $\widetilde 
    {\pi}:D\to {\cal {DC}}(B)$ de $\pi$. En particulier, tout
    morphisme de $C^*$-alg{\`e}bres $\pi:A\to {\cal 
    {DC}}(B)$ non-d{\'e}g{\'e}n{\'e}r{\'e} s'{\'e}tend de mani{\`e}re unique en $\widetilde {\pi}:{\cal
    {DC}}(A)\to {\cal {DC}}(B)$.\newline\indent \hfill$\square$  
\end{proposition}

\begin{remark}Soit $A$ une $C^*$-alg{\`e}bre (non n{\'e}cessairement
    unif{\`e}re) et $\chi$ un caract{\`e}re de $A$, alors $\chi$ est un morphisme non-d{\'e}g{\'e}n{\'e}r{\'e} $\chi:A\to
    {\cal {DC}}(\C)=\C$ et s'{\'e}tend donc 
    de mani{\`e}re unique en un homomorphisme $\widetilde {\chi}:{\cal
    {DC}}(A)\to \C$. On a $\widetilde {\chi}(T)\chi(a)=\chi(Ta)$ 
    pour tout $T\in {\cal {DC}}(A)$ et $a\in A$.\\
\end{remark}

\begin{remark}Soit $\pi:A\to {\cal {DC}}(B)$ un morphisme non-d{\'e}g{\'e}n{\'e}r{\'e}
  de $C^*$-alg{\`e}bres. Alors il existe une application continue
  $\pi^*:\mbox{Sp}(B)\to \mbox{Sp}(A)$ satisfaisant
  $\pi^*(\chi)(a)=\chi(\pi(a)b)$ pour tout $a\in A$ et $b\in B$ tel
  que $\chi(b)=1$.\\
\end{remark}

Soient $A$ une alg{\`e}bre involutive et $\cH$ un espace hilbertien. On dit qu'une
repr{\'e}sentation $L:A\to \B(\cH)$ de $A$ est \emph{non-d{\'e}g{\'e}ner{\'e}e} si 
l'ensemble $\{L(a)h:a\in A,h\in \cH\}$ engendre un sous-espace dense
de $\cH$. On dit que $L$ est \emph{fid{\`e}le} si $\ker L=\{0\}$. 

\begin{proposition}\label{extrep} Soit $A$ un id{\'e}al ferm{\'e} d'une
  $C^*$-alg{\`e}bre $B$ et $\pi_A$ une repr{\'e}sentation non-d{\'e}g{\'e}ner{\'e}e de
  $A$. Alors $\pi_A$ se prolonge de fa{\c c}on unique en une repr{\'e}sentation
  (non-d{\'e}g{\'e}ner{\'e}e) $\pi_B$ de $B$.\newline\indent \hfill$\square$
\end{proposition}

Supposons que $A$ est une sous-$C^*$-alg{\`e}bre de
$\B(\cH)$ et que $A\cH=\cH$. 

\begin{definition}On appelle \emph{alg{\`e}bre de
  multiplicateurs} de $A$ l'ensemble $$M(A):=\{x\in
\B(\cH)\,|\,xA\subset A\,\mbox{et}\,Ax\subset A\}.$$
\end{definition}

 L'alg{\`e}bre $M(A)$ est une
$C^*$-alg{\`e}bre contenant $A$ comme id{\'e}al essentiel.

\begin{example}Soit $\cH$ un espace de Hilbert. On note
  $\K(\cH)\subset \B(\cH)$ l'ensemble des op{\'e}rateurs compacts d{\'e}finis sur
  $\cH$. En fait, $\K(\cH)$ est un id{\'e}al ferm{\'e} de $\B(\cH)$,
  repr{\'e}sent{\'e} de fa{\c c}on  non-degener{\'e}e et l'on a $M(\K(\cH))=\B(\cH)$.\\ 
\end{example}

Suite {\`a} la proposition \ref{injmul} on a:

\begin{proposition}Si $A$ est repr{\'e}sent{\'e} de fa{\c c}on fid{\`e}le et
  non-d{\'e}g{\'e}n{\'e}r{\'e}e, alors l'application $x\mapsto (L_x,R_x)$ de $M(A)\to
  {\cal {DC}}(A)$ est un isomorphisme de $C^*$-alg{\`e}bres.\newline\indent
    \hfill$\square$\end{proposition}

\begin{remark}L'application r{\'e}ciproque est donn{\'e}e par la proposition
    \ref{extrep} appliqu{\'e} {\`a} $B={\cal {DC}}(A)$. \end{remark}

On id{\'e}ntifiera dans la suite $M(A)$ avec ${\cal {DC}}(A)$.

\section{Produits tensoriels}
On va donner rapidement quelques r{\'e}sultats classiques sur les produits
tensoriels de $C^*$-alg{\`e}bres qui seront utiles dans la suite
(\cf \cite{mur}, \cite{tak1}, \cite{was} et \cite{wo}).

 Soient $A$ et $B$ des espaces vectoriels. On note $A\odot B$ le \emph{produit tensoriel alg{\'e}brique} de $A$ et
 $B$ et $a\otimes b$ le \emph{tenseur {\'e}l{\'e}mentaire} qui est un {\'e}l{\'e}ment
 de $A\odot B$ avec $a\in A$ et $b\in B$. 

 Si $A$ et $B$ sont deux alg{\`e}bres, la formule $(a_1\otimes
  b_1,a_2\otimes b_2)\mapsto (a_1a_2)\otimes (b_1b_2)$ pour les
  tenseurs {\'e}l{\'e}mentaires permet de munir $A\odot B$ d'une structure
  d'alg{\`e}bre. Si $A$ et $B$ sont des alg{\`e}bres involutives, la formule
  $(a\otimes b)^*=a^*\otimes b^*$ permet de d{\'e}finir une involution sur
  $A\odot B$.

  Soient $A$, $B$ et $C$ des alg{\`e}bres involutives. Soient
  $\varphi:A\to C$, $\psi:B\to C$ des homomorphismes involutifs tels que $\varphi(a)\psi(b)=\psi(b)\varphi(a)$ pour
  tout $a\in A$, $b\in B$, alors l'application lin{\'e}aire $\varphi\times
  \psi:A\odot B\to C$ satisfaisant $(\varphi\times \psi)(a\otimes
  b)=\varphi(a)\psi(b)$ est un homomorphisme involutif.

 Soient $A,B,C,D$ des alg{\`e}bres involutives et $\varphi:A\to C$,
  $\psi:B\to D$ des homomorphismes involutifs. Alors l'application
  lin{\'e}aire 
   $\varphi\odot \psi:A\odot B\to C\odot D$
  v{\'e}rifiant $(\varphi\odot \psi)(a\otimes b)=\varphi(a)\otimes
  \psi(b)$ est un homomorphisme involutif.\\

Soient $A$ et $B$ des $C^*$-alg{\`e}bres et $||\cdot||_{\alpha}$ une
$C^*$-norme sur l'ag{\`e}bre involutive $A\odot B$. Le compl{\'e}t{\'e} de $A\odot
B$ par rapport {\`a} la norme $||\cdot||_{\alpha}$ est le \emph{produit
  tensoriel} de $A$ par $B$, not{\'e} $A\otimes_{\alpha}B$. C'est une
$C^*$-alg{\`e}bre.\\

Soient $\cH$ et $\cK$ des espaces hilbertiens, le \emph{produit tensoriel
  hilbertien} $\cH\otimes \cK$ est l'espace hilbertien obtenu en
  compl{\'e}tant $\cH\odot \cK$ pour le produit scalaire d{\'e}fini par
  $<h_1\otimes k_1,h_2\otimes k_2>=<h_1,h_2><k_1,k_2>$.

Soient $\cH_1$ et $\cH_2$ des espaces de Hilbert. Pour $T_1\in
  \B(\cH_1),\,T_2\in \B(\cH_2)$ l'application lin{\'e}aire $T_1\odot T_2$
  de $\cH_1\odot \cH_2$ dans lui-m{\^e}me d{\'e}finie par $(T_1\odot
  T_2)(h_1\otimes h_2)=T_1h_1\otimes T_2h_2$ s'{\'e}tend en un op{\'e}rateur $(T_1\otimes
  T_2)\in \B(\cH_1\otimes \cH_2)$ et l'on a $||T_1\otimes T_2||\leq
  ||T_1||\,||T_2||$. On d{\'e}crit ainsi l'injection naturelle de
  $\B(\cH_1)\odot \B(\cH_2)$ comme une sous-alg{\`e}bre involutive de
  $\B(\cH_1\otimes \cH_2)$. De plus on a $||T_1\otimes
  T_2||=||T_1||\,||T_2||$.

\medskip
 
Si $\pi_1:A_1\to \B(\cH_1)$ et $\pi_2:A_2\to
  \B(\cH_2)$ sont des repr{\'e}sentations de $C^*$-alg{\`e}bres, il existe un
  unique homomorphisme involutif $\pi_1\odot \pi_2:A_1\odot A_2\to
  \B(\cH_1\otimes \cH_2)$ v{\'e}rifiant $(\pi_1\odot \pi_2)(a_1\otimes
  a_2)=\pi_1(a_1)\otimes \pi_2(a_2)$ pour tout $a_1\in A_1$ et $a_2\in
  A_2$. Si $\pi_1$ et $\pi_2$ sont fid{\`e}les, $\pi_1\odot\pi_2$ est
  injectif.

\paragraph{Produit tensoriel maximal; produit tensoriel minimal.} En
  g{\'e}n{\'e}ral, on
  n'a pas unicit{\'e} de $C^*$-norme sur le produit tensoriel
  alg{\'e}brique de $C^*$-alg{\`e}bres. Ici on va discuter les deux cas
  extr{\^e}mes de $C^*$-normes.

\begin{proposition}Soient $A_1$ et $A_2$ des $C^*$-alg{\`e}bres. Toute
  $C^*$-semi-norme $||\cdot||_\nu$ sur $A_1\odot A_2$ v{\'e}rifie
  $||a_1\otimes a_2||_\nu\leq ||a_1||\,||a_2||$ pour tout $a_1\in
  A_1,\,a_2\in A_2$.\newline\indent \hfill$\square$ 
\end{proposition}

 Pour $t\in A_1\odot A_2$ on d{\'e}finit 
$$||t||_{\max}=\sup\{||t||_\nu\,\mbox{o{\`u}}\,||\cdot||_\nu\,\mbox{est
  une}\,C^*\mbox{-semi-norme sur}\,A_1\odot A_2\}.$$
Suite {\`a} la proposition pr{\'e}c{\'e}dente $||\cdot||_{\max}$ est finie et on peut voir
  que c'est une $C^*$-norme sur $A_1\odot A_2$ qui majore toutes les
  $C^*$-semi-normes de $A_1\odot A_2$ et qui v{\'e}rifie: $||a_1\otimes
  a_2||_{\max}=||a_1||\,||a_2||$, pour tout $a_1\in A_1$, $a_2\in A_2$.

 Cette
  $C^*$-norme est appel{\'e}e la $C^*$-norme \emph{maximale} de $A_1\odot
  A_2$. Le compl{\'e}t{\'e} de $A_1\odot A_2$ pour la norme maximale, not{\'e}e $A_1\otimes_{\max}A_2$ est une
  $C^*$-alg{\`e}bre que l'on appelle le \emph{produit tensoriel
  maximal} de $A_1$ et $A_2$. 

\begin{proposition}\label{mormax}Si $\psi_1:A_1\to B$, $\psi_2:A_2\to B$ sont des morphismes de
$C^*$-alg{\`e}bres tels que
$\psi_1(a_1)\psi_2(a_2)=\psi_2(a_2)\psi_1(a_1)$ pour tout $a_1\in
A_1$ et $a_2\in A_2$, alors $\psi_1\times \psi_2$ s'{\'e}tend en un morphisme
$\psi_1\times \psi_2:A_1\otimes_{\max}A_2\to B$.

Soient $\varphi_1:A_1\to B_1$ et $\varphi_2:A_2\to
  B_2$ des morphismes de $C^*$-alg{\`e}bres, alors $\varphi_1\odot
  \varphi_2$ s'{\'e}tend en un morphisme $\varphi_1\otimes_{\max}\varphi_2:A_1\otimes_{\max}A_2\to
  B_1\otimes_{\max}B_2$. Si $\varphi_1$ et $\varphi_2$ sont surjectifs
  alors $\varphi_1\otimes_{\max}\varphi_2$ est surjectif.

Si $A_1$ et $A_2$ sont des id{\'e}aux ferm{\'e}s de $B_1$ et $B_2$
respectivement alors
$\varphi_1\otimes_{\max}\varphi_2$ est un morphisme injectif. 
\newline\indent
\hfill$\square$  
\end{proposition}

\begin{proposition}Soient $A_1$ et $A_2$ des $C^*$-alg{\`e}bres et
  $\pi_1:A_1\to \B(\cH_1)$, $\pi_2:A_2\to \B(\cH_2)$, $\rho_1:A_1\to
  \B(\cK_1)$, $\rho_2:A_2\to \B(\cK_2)$ des repr{\'e}sentations
  fid{\`e}les. Alors, $||(\pi_1\odot \pi_2)(t)||=||(\rho_1\odot
  \rho_2)(t)||$ pour tout $t\in A_1\odot A_2$.\newline\indent
  \hfill$\square$ 
\end{proposition}

On d{\'e}finit ainsi une $C^*$-norme sur $A_1\odot A_2$ par $t\mapsto ||(\pi_1\odot \pi_2)(t)||=||t||_{\min}$, o{\`u}
  $\pi_1$ et $\pi_2$ sont des repr{\'e}sentations fid{\`e}les de $A_1$ et
  $A_2$ qu'on appelle norme \emph{minimale} ou \emph{spatiale} et elle v{\'e}rifie $||a_1\otimes
  a_2||_{\min}=||a_1||\,||a_2||$, pout tout $a_1\in A_1$, $a_2\in
  A_2$.\\
Le compl{\'e}t{\'e} de $A_1\odot A_2$ pour cette norme
  est une $C^*$-alg{\`e}bre qu'on la note $A_1\otimes_{\min}
  A_2$ et on l'appelle le \emph{produit tensoriel minimal ou 
  spatial} de $A_1$ et $A_2$.

\begin{proposition}\label{repfid} Soient $\pi_1$ et $\pi_2$ des repr{\'e}sentations (non
  n{\'e}cessairement fid{\`e}les) de $A_1$ et $A_2$, alors on a $||(\pi_1\odot
  \pi_2)(t)||\leq ||t||_{\min}$ pour tout $t\in A_1\odot
  A_2$. Autrement dit, il existe une
  unique repr{\'e}sentation $\pi_1\otimes \pi_2:A_1\otimes_{\min} A_2\to
  \B(\cH_1\otimes \cH_2)$ v{\'e}rifiant $(\pi_1\otimes \pi_2)(a_1\otimes
  a_2)=\pi_1(a_1)\otimes \pi_2(a_2)$ pour tout $a_1\in A_1$ et $a_2\in
  A_2$.

 Si $\pi_1$ et $\pi_2$ sont fid{\`e}les, $\pi_1\otimes\pi_2$ est
  une repr{\'e}sentation fid{\`e}le de $A_1\otimes_{\min}A_2$.\newline\indent
  \hfill$\square$ 
\end{proposition}

\begin{proposition}\label{morspatial} Soient $\varphi_1:A_1\to B_1$ et $\varphi_2:A_2\to
  B_2$ deux morphismes de $C^*$-alg{\`e}bres, alors il existe un unique
  morphisme
  $\varphi_1\otimes_{\min}\varphi_2:A_1\otimes_{\min}A_2\to
  B_1\otimes_{\min}B_2$ tel que
  $(\varphi_1\otimes_{\min}\varphi_2)(a_1\otimes a_2)=\varphi_1(a_1)\otimes_{\min}\varphi_2(a_2)$
  pour tout $a_1\in A_1$ et $a_2\in A_2$. De plus, si $\varphi_1$ et
  $\varphi_2$ sont injectifs, alors
  $\varphi_1\otimes_{\min}\varphi_2$ est injectif, si $\varphi_1$ et
  $\varphi_2$ sont surjectifs, alors
  $\varphi_1\otimes_{\min}\varphi_2$ est surjectif.\newline\indent
  \hfill$\square$ 

\end{proposition}

\begin{proposition} \textnormal{(\cf \cite{tak1})} La norme spatiale est minimum. Autrement dit, si $A$
  et $B$ sont des $C^*$-alg{\`e}bres, toute $C^*$-norme
  $||\cdot||_{\alpha}$ sur $A\odot B$ domine la norme spatiale, \ie
  $||t||_{\min}\leq ||t||_\alpha$ pour tout $t\in A\odot
  B$.\newline\indent \hfill$\square$  
\end{proposition} 

\begin{proposition}\label{multens} Soient $A_1$ et $A_2$ des
  $C^*$-alg{\`e}bres. Soit $||\cdot||_\alpha$ une $C^*$-norme sur
  $A_1\odot A_2$. Il existe des morphismes injectifs
  $i_{A_1}:M(A_1)\to M(A_1\otimes_\alpha A_2)$ et $i_{A_2}:M(A_2)\to
  M(A_1\otimes_\alpha A_2)$ tels que
  $i_{A_1}(a_1)i_{A_2}(a_2)=i_{A_2}(a_2)i_{A_1}(a_1)=a_1\otimes a_2$ pour tout $a_1\in A_1$,
  $a_2\in A_2$.\newline\indent \hfill$\square$ 
\end{proposition}  

\begin{remark}Soient $A$ et $B$ des $C^*$-alg{\`e}bres. Par les propositions
  \ref{multens} et \ref{mormax} on obtient un morphisme
  $M(A)\otimes_{\max} M(B)\to M(A\otimes_{\max} B)$.

Soient $\pi_1:A\to \B(\cH_1)$ et $\pi_2:B\to \B(\cH_2)$ des
repr{\'e}sentations fid{\`e}les. Par la proposition \ref{repfid} la
repr{\'e}sentation $\pi_1\otimes \pi_2$ de $A\otimes_{\min} B$ est aussi
fid{\`e}le. Notons $\widetilde{\pi_1}$ l'extension de $\pi_1$ {\`a} $M(A)$ et
$\widetilde{\pi_2}$ celle de $\pi_2$ {\`a} $M(B)$. Alors
$\widetilde{\pi_1}$ et $\widetilde{\pi_2}$ sont fid{\`e}les et on
obtient une repr{\'e}sentation fid{\`e}le $\widetilde{\pi_1}\otimes
\widetilde{\pi_2}$ de $M(A)\otimes_{\min} M(B)$. En id{\'e}ntifiant
$A\otimes_{\min} B$ avec son image par $\pi_1\otimes \pi_2$ on
remarque que l'image de $\widetilde{\pi_1}\otimes \widetilde{\pi_2}$
est contenue dans $M(A\otimes_{\min} B)$. Autrement dit on a
$M(A)\otimes_{\min} M(B)\subset M(A\otimes_{\min} B)$.   

\end{remark}

\begin{corollary}\label{mormult}Soient $\varphi_1:A_1\to M(B_1)$ et
  $\varphi_2:A_2\to M(B_2)$ deux morphismes de $C^*$-alg{\`e}bres. En composant $\varphi_1\otimes_{\min}\varphi_2$
  avec l'inclusion $M(B_1)\otimes_{\min}M(B_2)\to
  M(B_1\otimes_{\min}B_2)$ on a un morphisme $A_1\otimes_{\min}A_2\to
  M(B_1\otimes_{\min}B_2)$ que l'on note aussi
  $\varphi_1\otimes_{\min}\varphi_2$. Remarquons que le morphisme
  $\varphi_1\otimes_{\min}\varphi_2$ est injectif ou non-d{\'e}g{\'e}n{\'e}r{\'e} si
  $\varphi_1$ et $\varphi_2$ v{\'e}rifient la m{\^e}me propri{\'e}t{\'e}.\newline\indent
  \hfill$\square$ 
\end{corollary}

\begin{remark}Soient $B_1,B_2$ des $C^*$-alg{\`e}bres et $A_1$ (resp. $A_2$)
  un id{\'e}al ferm{\'e} de $B_1$ (resp. $B_2$), alors $A_1\odot A_2$ est
  un id{\'e}al de $B_1\odot B_2$.

 Posons $\alpha=\min$ ou $\max$. Notons $i_1:A_1\to B_1$ et
  $i_2:A_2\to B_2$ les inclusions naturelles, puisque
  $i_1\otimes_\alpha i_2:A_1\otimes_\alpha A_2\to B_1\otimes_\alpha B_2$ est
  continue alors $(i_1\otimes_\alpha i_2)(A_1\otimes_\alpha A_2)$ est un
  id{\'e}al ferm{\'e} de $B_1\otimes_\alpha B_2$. De plus, $i_1\otimes_\alpha
  i_2$ est injectif d'apr{\`e}s les propositions \ref{mormax} et \ref{morspatial}, donc on peut id{\'e}ntifier $A_1\otimes_\alpha A_2$ {\`a} un
  id{\'e}al ferm{\'e} de $B_1\otimes_\alpha B_2$. 
\end{remark}

\begin{remark}Soient $A$ et $B$ des $C^*$-alg{\`e}bres et $||\cdot||_\alpha$,
  $||\cdot||_\beta$ des $C^*$-normes sur $A\odot B$ telles que
  $||\cdot||_\alpha\leq ||\cdot||_\beta$, alors $A\otimes_{\alpha} B$ est un quotient de
  $A\otimes_{\beta} B$.

On en deduit que pour tout $C^*$-norme $||\cdot||_\alpha$ sur $A\odot B$,
$A\otimes_{\min} B$ est un quotient de $A\otimes_{\alpha} B$ qui
lui-m{\^e}me est un quotient de $A\otimes_{\max} B$. 

\end{remark}

\begin{proposition}\label{exmax} Soit $D$ une $C^*$-alg{\`e}bre et $0\to
  A\xrightarrow{i}B\xrightarrow{p}C\to 0$ une suite exacte de
  $C^*$-alg{\`e}bres, alors 
$$0\to A\otimes_{\max}D\xrightarrow{i\otimes_{\max} Id}B\otimes_{\max}
D\xrightarrow{p\otimes_{\max} Id}C\otimes_{\max}D\to 0$$
est une suite exacte de $C^*$-alg{\`e}bres. De plus, elle est scind{\'e}e si la suite \\ $0\to
  A\xrightarrow{i}B\ \mathop{\rightleftarrows}\limits^{p}_{\sigma}C\to
  0$ est exacte scind{\'e}e.\newline\indent
\hfill$\square$  
\end{proposition}

\begin{proposition}\label{exspa}  Soit $D$ une $C^*$-alg{\`e}bre et $0\to
  A\xrightarrow{i}B\ \mathop{\rightleftarrows}\limits^{p}_{\sigma}C\to
  0$ une suite exacte scind{\'e}e de $C^*$-alg{\`e}bres, alors on a une
  suite exacte scind{\'e}e de produits tensoriels minimaux 
 $$\xymatrix{
0 \ar[r]^{} & A\otimes_{\min}D \ar[r]^{i\otimes_{\min} Id} &
B\otimes_{\min}D  \ar@<2pt>[r]^{p\otimes_{min}
  Id} &  C\otimes_{\min}D \ar@<2pt>[l]^{\sigma\otimes_{\min} 
  Id} \ar[r] & 0.
}$$\newline\indent
\hfill$\square$

\end{proposition}

\paragraph{$C^*$-alg{\`e}bres nucl{\'e}aires.} Une $C^*$-alg{\`e}bre $A$ est
\emph{nucl{\'e}aire} si pour toute $C^*$-alg{\`e}bre $B$, il existe une unique
$C^*$-norme sur $A\odot B$, autrement dit si le morphisme naturel
$A\otimes_{\max}B\to A\otimes_{\min}B$ est injectif (dans ce cas $A\otimes_{\max}B\simeq A\otimes_{\min}B$).

Rappelons quelques propri{\'e}t{\'e}s des $C^*$-alg{\`e}bres nucl{\'e}aires.

\begin{itemize}

\item Un id{\'e}al ferm{\'e} $I$ d'une $C^*$-alg{\`e}bre nucl{\'e}aire $A$ est
  nucl{\'e}aire et le quotient $A/I$ est aussi nucl{\'e}aire.

\item L'extension d'une $C^*$-alg{\`e}bre nucl{\'e}aire par une $C^*$-alg{\`e}bre
  nucl{\'e}aire est nucl{\'e}aire. Plus pr{\'e}cisement: si la suite $0\to
  A\to B\to C\to 0$ est exacte et $A$ et $C$ sont
  nucl{\'e}aires, alors $B$ est nucl{\'e}aire.

\item Toute $C^*$-alg{\`e}bre commutative est nucl{\'e}aire.
 
\item Une limite inductive de $C^*$-alg{\`e}bres nucl{\'e}aires est nucl{\'e}aire.

\end{itemize}
On pourrait consulter  \cite{tak2}, \cite{was} ou \cite{wo} pour les
d{\'e}monstrations de ces r{\'e}sultats.

\paragraph{$C^*$-alg{\`e}bres exactes.} Une $C^*$-alg{\`e}bre $D$ est
  \emph{exacte} si pour toute suite exacte de $C^*$-alg{\`e}bres $0\to
  A\xrightarrow{i}B\xrightarrow{p}C\to 0$, la suite $$0\to A\otimes_{\min}D\xrightarrow{i\otimes_{\min} Id}B\otimes_{\min}
D\xrightarrow{p\otimes_{\min} Id}C\otimes_{\min}D\to 0$$ est exacte.

Rappelons quelques propri{\'e}tes des $C^*$-alg{\`e}bres exactes.

\begin{itemize}

\item Une sous-alg{\`e}bre d'une $C^*$-alg{\`e}bre exacte est exacte. 

\item Une $C^*$-alg{\`e}bre nucl{\'e}aire est exacte.

\item Un quotient d'une $C^*$-alg{\`e}bre exacte est exact (Kirchberg). 

\item Une limite inductive de $C^*$-alg{\`e}bres exactes est exacte.

\end{itemize}  

On peut trouver ces propri{\'e}t{\'e}s dans \cite{was}. Signalons l'important
r{\'e}sultat qu'ont r{\'e}cemment d{\'e}montr{\'e} Kirchberg et Philips
dans \cite{kp}. Une $C^*$-alg{\`e}bre s{\'e}parable est exacte si et seulement
si elle peut {\^e}tre plong{\'e}e dans une $C^*$-alg{\`e}bre nucl{\'e}aire.

\section{Produits crois{\'e}s}

\paragraph{$C^*$-syst{\`e}mes dynamiques.} On appele 
\emph{$C^*$-syst{\`e}me dynamique} un triplet $(A,G,\alpha)$ o{\`u} $A$ est
  une $C^*$-alg{\`e}bre, $G$ un groupe localement compact et $\alpha$
  une application $g\mapsto \alpha_g$ qui est un
  homomorphisme continu de $G$ dans le groupe $\Aut(A)$ des automorphismes
  involutifs de $A$ muni de la topologie de la convergence simple (\ie pour tout $a\in A$ l'application
$g\mapsto {\alpha}_g(a)$ est continue).

Soit $(A,G,\alpha)$ un $C^*$-syst{\`e}me dynamique. Notons $\lambda$ une
  mesure de Haar {\`a} gauche sur $G$. Notons $\Delta:G\to {\mathbb R}^*_{+}$ la fonction
  modulaire de $G$ d{\'e}finie par:
$$\Delta(r)\int_G h(sr)d\lambda(s)=\int_G h(s)d\lambda(s)\;\mbox{pour tout}\,
 h\in C_c(G),r\in G.$$ Rappelons que pour tout $f\in
 C_c(G)$,$$\int_Gf(s^{-1})\Delta(s^{-1})d\lambda(s)=\int_Gf(s)d\lambda(s).$$

\paragraph{L'alg{\`e}bre involutive.} L'espace $C_c(G,A)$ des fonctions
  continues $\phi:G\to A$ {\`a} support compact est muni d'une structure
  d'alg{\`e}bre involutive donn{\'e}e par:
\begin{enumerate}
\item un produit appel{\'e} convolution $$(f*g)(s)=\int_Gf(r){\alpha}_r(g(r^{-1}s))d\lambda(r),\;f,g\in
C_c(G,A),\;s\in G,$$
\item une involution
$$f^*(s)=\Delta(s^{-1}){\alpha}_s({f(s^{-1})}^*),\; f\in
C_c(G,A)$$ 
\end{enumerate}
 On d{\'e}finit une norme sur $C_c(G,A)$
par:$${||f||}_1=\int_G||f(s)||d\lambda(s).$$
Elle verifie ${||f||}_1={||f^*||}_1$ et ${||f*g||}_1\leq
{||f||}_1{||g||}_1$ gr{\^a}ce aux propri{\'e}t{\'e}s de la mesure de Haar. On
l'appellera norme $L^1$. Le compl{\`e}t{\'e} de
$C_c(G,A)$ pour la norme $L^1$ est une alg{\`e}bre de Banach involutive
not{\'e}e $L^1(G,A)$.\\

\paragraph{Repr{\'e}sentations covariantes.} Soit $(A,G,\alpha)$ un
$C^*$-syst{\`e}me dynamique. On appelle
\emph{repr{\'e}sentation covariante} de $(A,G,\alpha)$ une paire $({\pi}_A,U)$ o{\`u} ${\pi}_A:A\to \mathbb B(\cal H)$ est une
repr{\'e}sentation de $A$ et $U:G\to \mathbb U(\cal H)$ est un
homomorphisme $*$-fortement continu de $G$ dans le groupe des
unitaires de $\B(\cH)$ (on l'appelera aussi une \emph{repr{\'e}sentation
  unitaire} de $G$) not{\'e}e $U(g)=U_g$ et telle que,
${\pi}_A({\alpha}_s(a))=U_s{\pi}_A(a)U_s^*$ pour tout $a\in
A,s\in G.$

On dira que $(\pi_A,U)$ est non-d{\'e}g{\'e}n{\'e}r{\'e}e si $\pi_A$ est non-d{\'e}g{\'e}n{\'e}r{\'e}e.

\begin{example}\emph{Repr{\'e}sentations r{\'e}guli{\`e}res {\`a} gauche}\\
 Soit ${\rho}_A:A\to \mathbb B(\cal H)$ une repr{\'e}sentation de
$A$ sur $\cal H$. Le couple $(\widetilde{{\rho}_A},U)$ de repr{\'e}sentations
sur l'espace de Hilbert $L^2(G,\cal H)$ d{\'e}finies par :\\
$\widetilde{{\rho}_A}:A\to \mathbb B(L^2(G,\cal H))$
telle que
$\widetilde{{\rho}_A}(a)h(r)={\rho}_A({\alpha}_r^{-1}(a))(h(r))$,
$a\in A,r\in G,h\in L^2(G,\cH)$ et
$U:G\to \mathbb U(L^2(G,\cal H))$ telle que $U_{s}h(r)=h(s^{-1}r)$
pour $r,s\in G,h\in L^2(G,\cH)$ est une repr{\'e}sentation covariante de $(A,G,\alpha)$ qu'on appelle repr{\'e}sentation
r{\'e}guli{\`e}re {\`a} gauche associ{\'e}e {\`a} $\rho_A$. De plus elle est non-d{\'e}g{\'e}ner{\'e}e si ${\rho}_A$
l'est.
\end{example}

\begin{proposition}\label{rc}Soit $({\pi}_A,U)$ une repr{\'e}sentation covariante de $(A,G,\alpha)$ sur
$\cal H$. Alors $$({\pi}_A\times U)(f)=\int_G{\pi}_A(f(s))U_s
d\lambda(s)\;,\,f\in C_c(G,A)$$ d{\'e}finit une repr{\'e}sentation involutive de
$C_c(G,A)$ sur $\cal H$ qui est born{\'e}e pour la norme $L^1$. Par
cons{\'e}quent elle d{\'e}finit une repr{\'e}sentation involutive de l'alg{\`e}bre de
Banach involutive $L^1(G,A)$. 

La repr{\'e}sentation $({\pi}_A\times U)$ est non-d{\'e}g{\'e}ner{\'e}e
si ${\pi}_A$ est non-d{\'e}g{\'e}ner{\'e}e. 

De plus, l'application $({\pi}_A,U)\mapsto {\pi}_A\times U$ est une bijection entre les repr{\'e}sentations
covariantes non-d{\'e}g{\'e}ner{\'e}es de $(A,G,\alpha)$ et les repr{\'e}sentations
non-d{\'e}gener{\'e}es de $L^1(G,A)$.\newline\indent \hfill$\square$ 
\end{proposition}

\paragraph{Produit crois{\'e}.}Soit $(A,G,\alpha)$ un
$C^*$-syst{\`e}me dynamique. On d{\'e}finit une norme sur $C_c(G,A)$ par
$${||f||}_u=\sup\{||({\pi}_A\times
U)(f)||\;\mbox{o{\`u}}\;({\pi}_A,U)\;\mbox{repr{\'e}sentation\;covariante\;de}\;(A,G,\alpha)\}.$$
On l'appelle norme \emph{universelle} et elle major{\'e}e par la norme
$L^1$. Le complet{\'e} de $C_c(G,A)$ ou de $L^1(G,A)$ pour la norme
universelle, not{\'e} $A{\rtimes}_{\alpha}G$ est le \emph{produit crois{\'e}
  maximal} de $A$ par $G$ et c'est une $C^*$-alg{\`e}bre.

Remarquons que l'on a (par la proposition \ref{rc})
$${||f||}_u=\sup\{||L(f)||\;\mbox{o{\`u}}\;L\;\mbox{repr{\'e}sentation\;de}\;L^1(G,A)\}.$$
En d'autres termes, $A{\rtimes}_{\alpha}G$ est la $C^*$-alg{\`e}bre \emph{enveloppante} de $L^1(G,A)$.\\

 Le compl{\'e}t{\'e} de $C_c(G,A)$ ou de $L^1(G,A)$ pour la norme $${||f||}_r=\sup\{||(\widetilde{{\rho}_A}\times
U)(f)||\;\mbox{o{\`u}}\;{\rho}_A\;\mbox{repr{\'e}sentation\;de}\;A\}$$ est le
\emph{produit crois{\'e} reduit} de $A$ par $G$ qu'on le note
$A{\rtimes}_{r,{\alpha}}G$ et qui est une $C^*$-alg{\`e}bre. 

\begin{proposition}\label{fidele}
Soit $(A,G,\alpha)$ un $C^*$-systeme dynamique et $\pi_A$ une
repr{\'e}sentation fid{\`e}le de $A$ alors $\widetilde{\pi_A}\times U$ est une
repr{\'e}sentation fid{\`e}le de $A{\rtimes}_{r,{\alpha}}G$.\newline\indent \hfill$\square$
\end{proposition}

Donnons trois autres fa{\c c}ons d'{\'e}noncer cette proposition. On suppose qu'on a un $C^*$-systeme dynamique $(A,G,\alpha)$ et une
repr{\'e}sentation fid{\'e}le $\pi_A$ de $A$, alors
 \begin{itemize}
\item ${||f||}_r=||\widetilde{{\pi}_A}\times U(f)||\;$ pour tout $f\in
C_c(G,A)$.
\item $A{\rtimes}_{r,{\alpha}}G=A{\rtimes}_{\alpha}G/\ker(\widetilde{{\pi}_A}\times
U)$.
\item $A{\rtimes}_{r,{\alpha}}G\simeq (\widetilde{{\pi}_A}\times
U)(A{\rtimes}_{\alpha}G)$.
\end{itemize}

\bigskip

\begin{example}\label{comt}(\cf \cite{ta}) Soit $G$ un groupe localement compact
  et $C_0(G)$ la $C^*$-alg{\`e}bre des fonctions continues d{\'e}finies sur
  $G$ qui tendent vers z{\`e}ro {\`a} l'infini. Alors le produit crois{\'e}
  $C_0(G)\rtimes G$ pour l'action continue de $G$ donn{\'e}e par les
  translations {\`a} gauche est isomorphe {\`a} la $C^*$-alg{\`e}bre
  $\K(L^2(G))$ des op{\'e}rateurs compacts d{\'e}finis sur $L^2(G)$.
\end{example}

\paragraph{Morphismes {\'e}quivariants.} Soient $A,B$ des $C^*$-alg{\`e}bres
munies des actions $\alpha$, $\beta$ d'un groupe $G$ et soit $\varphi:A\to B$ un
morphisme de $C^*$-alg{\`e}bres. On dit que $\varphi$ est \emph{{\'e}quivariant} si
$\varphi[{\alpha}_g(a)]={\beta}_g[\varphi(a)]$ pour tout $g\in G$ et
tout $a\in A$.

Pour $g\in G$ notons $\widetilde{\beta_g}$ l'extension de $\beta_g$ 
{\`a} l'alg{\`e}bre des multiplicateurs $M(B)$ qui est donn{\'e}e par:
$\widetilde{\beta_g}(b)\beta_g(b')=\beta_g(bb')$ pour tout $b\in
M(B),\,b'\in B$ et $g\in G$.
On dit qu'un morphisme $\pi:A\to M(B)$ est {\'e}quivariant si
$\pi(\alpha_g(a))=\widetilde{\beta_g}(\pi(a))$ pour tout $g\in
G, a\in A$, autrement dit si
$\pi(\alpha_g(a))\beta_g(b)=\beta_g(\pi(a)b)$ pour tout $g\in
G, a\in A$ et $b\in B$.  

\begin{proposition}\label{crmax}
Soient $(A,G,\alpha)$ et $(B,G,\beta)$ deux $C^*$-syst{\'e}mes dynamiques et
$\varphi:A\to B$ un morphisme {\'e}quivariant. Alors il existe un unique
morphisme ${\varphi}_*: A{\rtimes}_{\alpha}G\to B{\rtimes}_{\beta}G$
donn{\'e} par ${\varphi}_*(f)(s)=\varphi(f(s))$ pour tout $f\in C_c(G,A)$. Plus
pr{\'e}cis{\'e}ment, on a un diagramme commutatif

$$\xymatrix{
           C_c(G,A) \ar[r]^{\varphi_*} \ar[d]_{} & C_c(G,B) \ar[d]^{} \\
           A{\rtimes}_{\alpha}G \ar[r]_{{\varphi}_{*}} & B{\rtimes}_{\beta}G
         }$$
o{\`u} les fl{\`e}ches verticales sont les applications canoniques dans le
           compl{\`e}t{\'e} et o{\`u} l'on a not{\'e} $\varphi_*:C_c(G,A)\to C_c(G,B)$
           l'application $f\mapsto \varphi\circ f$.
\end{proposition}

\begin{proof}
Puisque $\varphi$ est {\'e}quivariant on v{\'e}rifie imm{\'e}diatement que
 l'application $\varphi_*:C_c(G,A)\to C_c(G,B)$ est un homomorphisme
 involutif.

 De plus, puisque $\varphi$ est un morphisme de $C^*$-alg{\`e}bres, on a
${||\varphi_*(f)||}_{1}=\int_G||\varphi(f(s))||d\lambda(s)\leq
  \int_G||f(s)||d\lambda(s)={||f||}_{1}$. Donc $\varphi_*$ admet un unique
  prolongement {\`a} $L^1(G,A)$.

 Soit $L$ une
  repr{\'e}sentation de $L^1(G,B)$, alors $L\circ \varphi_*$ est une
  repr{\'e}sentation de $L^1(G,A)$ et $||L\circ \varphi_*(f)||\leq
  {||f||}_{u}$ pour tout $f\in L^1(G,A)$. On en deduit que
  ${||\varphi_*(f)||}_{u}\leq {||f||}_{u}$ pour tout $f\in
  L^1(G,A)$. Donc il existe une
unique extension de $\varphi_*$ {\`a} un morphisme ${\varphi}_*: A{\rtimes}_{\alpha}G\to
B{\rtimes}_{\beta}G$.  \end{proof}

On a {\'e}galement,

\begin{proposition}\label{crred}
Soient $(A,G,\alpha)$ et $(B,G,\beta)$ deux $C^*$-syst{\'e}mes dynamiques et
$\varphi:A\to B$ un morphisme {\'e}quivariant. Alors il existe un unique
morphisme ${\varphi}_{*,r}: A{\rtimes}_{r,{\alpha}}G\to
B{\rtimes}_{r,{\beta}}G$ donn{\'e} par
${\varphi}_{*,r}(f)(s)=\varphi(f(s))$ pour tout $f\in C_c(G,A)$.

\end{proposition}

\begin{proof} 
Soit ${\pi}_B$ une repr{\'e}sentation fid{\'e}le de $B$ et soit
$\widetilde{{\pi}_B}\times U$ la repr{\'e}sentation de
$B{\rtimes}_{\beta}G$ correspondante {\`a} la repr{\'e}sentation 
$(\widetilde{{\pi}_B},U)$. En posant ${\pi}_B\circ \varphi\equiv
{\pi}_A$ on a par un simple calcul que $(\widetilde{{\pi}_B}\times U)\circ
  {\varphi}_*=(\widetilde{{\pi}_A}\times U)$ donc
  ${||{\varphi}_*(f)||}_{r}\leq {||f||}_{r}$ pour tout $f\in
  C_c(G,A)$ et $\varphi_*$ admet un unique prolongement {\`a} $A{\rtimes}_{r,\alpha}G$. \end{proof}

\begin{corollary}\label{comm}
Soient $(A,G,\alpha)$ et $(B,G,\beta)$ deux $C^*$-syst{\'e}mes dynamiques et
$\varphi:A\to B$ un morphisme {\'e}quivariant. Soient
${\lambda}_A\;($resp.${\lambda}_B)$ la surjection canonique
$A{\rtimes}_{\alpha}G\to
A{\rtimes}_{r,\alpha}G\;($resp.$B{\rtimes}_{\beta}G\to
B{\rtimes}_{r,\beta}G)$, alors le diagramme suivant est commutatif:
$$\xymatrix{
           A{\rtimes}_{\alpha}G \ar[r]^{{\varphi}_*} \ar[d]_{{\lambda}_A} & B{\rtimes}_{\beta}G \ar[d]^{{\lambda}_B} \\
           A{\rtimes}_{r,\alpha}G \ar[r]_{{\varphi}_{*,r}} & B{\rtimes}_{r,\beta}G
         .}$$\newline\indent \hfill$\square$
\end{corollary}

Il est clair que la construction est naturelle au sens suivant

\begin{proposition}\textbf{(Naturalit{\'e} des produits crois{\'e}s)}\label{comp}\\
Soient $(A,G,\alpha)$, $(B,G,\beta)$, $(C,G,\gamma)$ des
$C^*$-syst{\`e}mes dynamiques et soient $\varphi:A\to B$, $\psi:B\to
C$ des morphismes equivariants. On a
${({\mbox{Id}}_A)}_*={\mbox{Id}}_{A{\rtimes}_{\alpha}G}$,
${({\mbox{Id}}_A)}_{*,r}={\mbox{Id}}_{A{\rtimes}_{r,\alpha}G}$ et
${(\psi\circ \varphi)}_*=\psi_*\circ \varphi_*$, ${(\psi\circ
  \varphi)}_{*,r}=\psi_{*,r}\circ \varphi_{*,r}$.
\newline\indent \hfill$\square$\end{proposition}

 On dit qu'une sous-$C^*$-alg{\`e}bre $B$ de $A$ est
 \emph{$G$-invariante} si elle est stable par les automorphismes
 ${\alpha}_g\;,g\in G$.

\begin{proposition}\label{inj} Soit $(B,G,\beta)$ un $C^*$-syst{\`e}me
 dynamique et soit $A$ un id{\'e}al ferm{\'e} $G$-invariant de $B$.
 Notons $i:A\hookrightarrow B$ l'inclusion et encore
  ${\beta}_g$ la restriction de $\beta_g$ {\`a} $A$. Alors le morphisme
  $i_*:A{\rtimes}_{\beta}G\to B{\rtimes}_{\beta}G$ est injectif.
\end{proposition}

\begin{proof} Soit ${||\cdot||}_{u,B}$ la norme universelle sur $C_c(G,B)$ et
${||\cdot||}_{u,A}$ la norme universelle sur $C_c(G,A)$.

Pour tout $f\in A{\rtimes}_\beta G$ on a $||i_*(f)||\leq ||f||$ parce
que $i_*$ est un morphisme de $C^*$-alg{\`e}bres.

D'autre part soit $\pi\times U$ une repr{\'e}sentation fid{\`e}le
non-d{\'e}g{\'e}ner{\'e}e de $A{\rtimes}_{\beta}G$. La repr{\'e}sentation $\pi:A\to
\Bbb B(\cal H)$ admet une unique extension ${\pi}_B:B\to \Bbb B(\cal H)$
et par unicit{\'e} de l'extension $({\pi}_B,U)$ est une repr{\'e}sentation
covariante de $(B,G,\beta)$. Donc on a
${||f||}_{u,A}=||\pi\times U(f)||=||{\pi}_B\times U(f)||\leq
{||f||}_{u,B}$ pour tout $f\in \k$. On en deduit que
${||f||}_{u,A}={||f||}_{u,B}$ pour tout $f\in\k$ donc $i_*$ est
injectif.\end{proof}

\begin{proposition}\label{surj}
Soit $p:A\to B$ un morphisme {\'e}quivariant surjectif alors
$p_{*}:A{\rtimes}_{\alpha}G\to B{\rtimes}_{\beta}G$ et
$p_{*,r}:A{\rtimes}_{r,\alpha}G\to B{\rtimes}_{r,\beta}G$ sont des morphismes
surjectifs.\end{proposition}
\begin{proof} On note $\cal B$ l'espace de fonctions de la forme
$\sum\limits_{i=1}^n{\varphi}_{i} b_i\;,{\varphi}_i\in C_c(G),b_i\in
  B$. L'espace $\cal B$ est
  dense dans $C_c(G,B)$ pour la topologie de la convergence uniforme {\`a}
  support contenu dans un compact fix{\'e}, donc $\cal
  B$ est dense dans $C_c(G,B)$ pour la norme $L^1$ et par suite il est dense
  dans $B{\rtimes}_{\beta}G$.

On a de plus que
  $p_*(A{\rtimes}_{\alpha}G)$ est une sous-alg{\`e}bre ferm{\'e}e de
  $B{\rtimes}_{\beta}G$ car $p_*$ est un morphisme de
  $C^*$-alg{\`e}bres. Puisque l'image Im$p_*$ contient $\cal B$ qui lui
  m{\^e}me est dense dans
  $B{\rtimes}_{\beta}G$, $p_*$ est surjectif.

 A l'aide de \ref{comm}
  on a que $p_{*,r}$ est aussi surjectif.\end{proof}

\begin{proposition}\label{rinj} Si $i:A\to B$ est un morphisme
  {\'e}quivariant injectif de $C^*$-alg{\`e}bres alors, $i_{*,r}:A{\rtimes}_{r,\alpha}G\to B{\rtimes}_{r,\beta}G$ est
injectif.\end{proposition}

\begin{proof} Le resultat est immediat suite {\`a} la proposition \ref{fidele}
\end{proof}

\bigskip

Soit $(B,G,\beta)$ un $C^*$-systeme dynamique et $A$ un
id{\'e}al (ferm{\'e}) $G$-invariant de $B$. Notons $i:A\to B$ l'inclusion. Alors $C_c(G,A)$ est un
ideal de $C_c(G,B)$. Puisque $i_*:A{\rtimes}_{\beta}G\to
B{\rtimes}_{\beta}G$ et $i_{*,r}:A{\rtimes}_{r,\beta}G\to
B{\rtimes}_{r,\beta}G$ sont continus alors $i_*(A{\rtimes}_{\beta}G)$
(resp. $i_{*,r}(A{\rtimes}_{r,\beta}G)$) est un id{\'e}al ferm{\'e} de
$B{\rtimes}_{\beta}G$ (resp. $B{\rtimes}_{r,\beta}G$). Puisque $i_*$, $i_{*,r}$ sont injectifs on id{\'e}ntifiera 
$A{\rtimes}_{\beta}G$ (resp. $A{\rtimes}_{r,\beta}G$) {\`a} un id{\'e}al
ferm{\'e} de $B{\rtimes}_{\beta}G$
(resp.$B{\rtimes}_{r,\beta}G$).

\bigskip 

\begin{proposition}\label{mu} Soient $(A,G,\alpha)$ et $(B,G,\beta)$
  deux $C^*$-syst{\'e}mes dynamiques. Soit $\pi:A\to M(B)$ un morphisme
  {\'e}quivariant. Alors il existe des morphismes
  uniques  $\widetilde{\pi}:A{\rtimes}_{\alpha}G\to
  M(B{\rtimes}_{\beta}G)$ et $\widetilde{\pi_r}:A{\rtimes}_{r,\alpha}G\to
  M(B{\rtimes}_{r,\beta}G)$ donn{\'e}s par

  $$(\widetilde{\pi}(f)h)(s)=(\widetilde{\pi_r}(f)h)(s)=\int_G\pi(f(r))\beta_r(h(r^{-1}s))d\lambda(r),$$
  pour tout $f\in C_c(G,A),h\in C_c(G,B)$ et $s\in G$. De plus, si $\pi$ est injective alors
  $\widetilde{\pi_r}$ l'est aussi. 
\end{proposition}

\begin{proof}

 Soit $L$ une repr{\'e}sentation fid{\`e}le non-d{\'e}g{\'e}n{\'e}r{\'e}e de
  $B{\rtimes}_{\beta}G$. Par la proposition \ref{rc}, la
  repr{\'e}sentation $L$ est de la forme $L=\varphi_B\times U$
  o{\`u} $(\varphi_B,U)$ est une repr{\'e}sentation covariante non-d{\'e}g{\'e}n{\'e}r{\'e}e
  de $(B,G,\beta)$. Notons $\overline{L}=\overline{\varphi_B\times
  U}:M(B\rtimes_{\beta}G)\to \B(\cH)$ son extension {\`a} l'alg{\`e}bre des
  multiplicateurs qui est aussi fid{\`e}le.  

\medskip

On d{\'e}finit
  $\varphi_A=\overline{\varphi_B}\circ \pi$ o{\`u}  $\overline{\varphi_B}$
  est l'unique extension de la repr{\'e}sentation $\varphi_B:B\to \B(\cH)$
  {\`a} l'alg{\`e}bre des multiplicateurs $M(B)$. Alors, $\varphi_A$ est
  une repr{\'e}sentation de $A$ et avec un simple calcul, en utilisant que
  $\pi$ est {\'e}quivariante, on montre que la
  repr{\'e}sentation $(\varphi_A,U)$ est covariante. Donc par la
  proposition \ref{rc}, on peut d{\'e}finir une repr{\'e}sentation
  $\varphi_A\times U$ de $A{\rtimes}_{\alpha}G$. L'image de
  $A{\rtimes}_{\alpha}G$ est dans les multiplicateurs et l'on obtient un homomorphisme $\widetilde{\pi}:A{\rtimes}_{\alpha}G\to
  M(B{\rtimes}_{\beta}G)$ tel que $\varphi_A\times U=\overline{\varphi_B\times
  U}\circ \widetilde{\pi}$.
\medskip

Pour le produit crois{\'e} reduit, soit $\varphi_B$ une repr{\'e}sentation
  fid{\`e}le et non-degener{\'e}e de $B$ et $\overline{\varphi_B}$ son
  extension {\`a} l'alg{\`e}bre des multiplicateurs qui est aussi fid{\`e}le. Par la proposition
  \ref{fidele}, $L=\widetilde{\varphi_B}\times U$ est une
  repr{\'e}sentation fid{\`e}le non-deg{\'e}n{\'e}r{\'e}e de $B{\rtimes}_{r,\beta}G$, avec $(\widetilde{\varphi_B},U)$ la repr{\'e}sentation
  r{\'e}guli{\`e}re {\`a} gauche de $(B,G,\beta)$ associ{\'e}e {\`a} $\varphi_B$. Il est clair que
  $\widetilde{\varphi_A}=(\overline{\widetilde{\varphi_B}}\circ
  \pi,U)$ est la repr{\'e}sentation r{\'e}guli{\`e}re {\`a} gauche de $(A,G,\alpha)$
  associ{\'e}e {\`a} $\varphi_A=\overline{\varphi_B}\circ \pi$, qui est
  covariante. Donc, $\widetilde{\varphi_A}\times
  U$ est une repr{\'e}sentation de $A\rtimes_{r,\alpha}G$ et on obtient un homomorphisme $\widetilde{\pi_r}:A{\rtimes}_{r,\alpha}G\to
  M(B{\rtimes}_{r,\beta}G)$ tel que $\widetilde{\varphi_A}\times
  U=\overline{\widetilde{\varphi_B}\times U}\circ \widetilde{\pi_r}$.

Si $\pi$ est injective $\varphi_A=\overline{\varphi_B}\circ \pi$ est
  fid{\`e}le. Donc la repr{\'e}sentation $\widetilde{\varphi_A}\times
  U$ est fid{\`e}le et on en deduit que $\widetilde{\pi_r}$ est injective. 

\end{proof}

\begin{proposition} Soient
           $(A,G,\alpha),\,(B,G,\beta),\,(C,G,\gamma)$ des $C^*$-syst{\`e}mes
dynamiques. On suppose qu'on a
une suite exacte de $C^*$-alg{\`e}bres,$$0\to
A\xrightarrow{i}B\xrightarrow{p}C\to 0$$ o{\`u} $i,p$ sont des morphismes
{\'e}quivariants. Alors les produits crois{\'e}s maximaux forment aussi une
suite exacte:
$$0\to
A{\rtimes}_{\alpha}G\xrightarrow{i_{*}}B{\rtimes}_{\beta}G\xrightarrow{p_*}C{\rtimes}_{\gamma}G\to
0.$$\end{proposition}
 
\begin{proof}  Par la proposition \ref{inj}, puisque
  $A{\rtimes}_{\alpha}G$ est un ideal de $B{\rtimes}_{\beta}G$ on a
  que $i_*$ est injectif. Par la proposition \ref{surj}, $p_*$ est
  surjectif.

 Montrons que $\ker p_*=$Im$i_*=A{\rtimes}_{\alpha}G$. Nous devons
  d{\'e}montrer que le morphisme
  $j:B{\rtimes}_{\beta}G/A{\rtimes}_{\alpha}G\to C{\rtimes}_{\gamma}G$
  d{\'e}duit de $p_*$ est injectif. Soit $L$ une repr{\'e}sentation fid{\`e}le 
  non-d{\'e}g{\'e}ner{\'e}e de
$B{\rtimes}_{\beta}G/A{\rtimes}_{\alpha}G$ alors $L$ est de la forme
suivante: $L={\pi}_B\times U$ repr{\'e}sentation non-d{\'e}g{\'e}ner{\'e}e de $B{\rtimes}_{\beta}G$
(o{\`u} $({\pi}_B,U)$ est une repr{\'e}sentation covariante de $(B,G,\beta)$
et ${\pi}_B$ repr{\'e}sentation non-d{\'e}g{\'e}ner{\'e}e de $B$) telle que $({\pi}_B\times
U)(A{\rtimes}_{\alpha}G)=0$.

On a $\pi_B(A)({\pi}_B\times
U)(B{\rtimes}_{\beta}G)\subset ({\pi}_B\times
U)(A{\rtimes}_{\alpha}G)$. Puisque $({\pi}_B\times
U)(A{\rtimes}_{\alpha}G)\\=0$ alors $\pi_B(A)({\pi}_B\times
U)(B{\rtimes}_{\beta}G)=0$, or la repr{\'e}sentation ${\pi}_B\times U$
de $B{\rtimes}_{\beta}G$ est non-d{\'e}g{\'e}ner{\'e}e donc
${{\pi}_{B}|}_{A}=0$. Il existe donc une unique repr{\'e}sentation ${\pi}_C$
de $C$ telle que $\pi_B={\pi}_C\circ p$. Or, $({\pi}_C,U)$ est une
repr{\'e}sentation covariante de $(C,G,\gamma)$ car: d'une part on a $\pi_B({\beta}_g(b))={\pi}_C(p({\beta}_g(b)))={\pi}_C({\gamma}_g(p(b)))$
pour tout $g\in G,b\in B$ et d'autre part on a
$\pi_B({\beta}_g(b))=U_g\pi_B(b)U_g^*=U_g({\pi}_C(p(b))U_g^*$ pour
tout $g\in G$. Alors, $L=({\pi}_C\times U)\circ j$. Comme $L$ est
  fid{\`e}le, $j$ est injectif.

Il est clair que $p_*\circ i_*=0$.

\end{proof}

\begin{corollary}\label{exactmax}
Si on suppose que la suite exacte $$0\to A\xrightarrow{i}
B\ \mathop{\rightleftarrows}\limits^p_\sigma\ C \to 0$$ 
est scind{\'e}e (avec $i$, $p$, et $\sigma$ des morphismes {\'e}quivariants), alors on obtient aussi une suite
exacte scind{\'e}e des produits crois{\'e}s maximaux,
$$0\to A{\rtimes}_{\alpha}G\xrightarrow{i_{*}} B{\rtimes}_{\beta}G\
\mathop{\rightleftarrows}\limits^{p_{*}}_{{\sigma}_{*}}\
C{\rtimes}_{\gamma}G \to 0.$$
\end{corollary}

\begin{proof} C'est {\'e}vident puisque par la proposition \ref{comp} on a ${(p\circ
  \sigma)}_*=p_*\circ{\sigma}_*={\mbox{Id}}_{C{\rtimes}_{\gamma}G}$.\end{proof}

 En d'autre termes, on a un isomorphisme lin{\'e}aire $$B{\rtimes}_{\beta}G=i_{*}(A{\rtimes}_{\alpha}G)\bm\oplus
{\sigma}_{*}(C{\rtimes}_{\gamma}G)\simeq
(A{\rtimes}_{\alpha}G)\bm\oplus (C{\rtimes}_{\gamma}G).$$

\begin{proposition}
 Soient $(A,G,\alpha)\;,\;(B,G,\beta)\;,\;(C,G,\gamma)$ des $C^*$-systemes
dynamiques. On suppose qu'on a
une suite exacte scind{\'e}e de $C^*$-alg{\`e}bres,$$0\to
A\xrightarrow{i}B \ \mathop{\rightleftarrows}\limits^{p}_\sigma\ C \to
0$$ o{\`u} $i,p$ et $\sigma$ sont des morphismes
{\'e}quivariants. Alors les produits crois{\'e}s reduits forment aussi une
suite exacte scind{\'e}e:
$$0\to
A{\rtimes}_{r,\alpha}G\xrightarrow{i_{*,r}}B{\rtimes}_{r,\beta}G \
\mathop{\rightleftarrows}\limits^{p_{*,r}}_{{\sigma}_{*,r}}\ C{\rtimes}_{r,\gamma}G\to
0.$$
\end{proposition}
 
\begin{proof} On a vu que $i_{*,r}$ est injective (proposition
  \ref{rinj}) et on a  $p_{*,r}\circ
 {\sigma}_{*,r}={\mbox{Id}}_{C{\rtimes}_{r,\gamma}G}.$ 

Notons $\lambda_K$
 avec $K=A,B,C$ les surjections canoniques respectivement de
 $K{\rtimes}_k G\to K{\rtimes}_{r,k}G$,
 $k=\alpha,\beta,\gamma$. Puisqu'on a une suite exacte scind{\'e}e 
           $$B{\rtimes}_{\beta}G=i_{*}(A{\rtimes}_{\alpha}G)\bm\oplus{\sigma}_{*}(C{\rtimes}_{\gamma}G),$$
           le diagramme   
$$\xymatrix{
           0 \ar[r]^{} & A{\rtimes}_{\alpha}G \ar[r]^{i_*}
           \ar[d]_{\lambda_A} & B{\rtimes}_{\beta}G \ar@<2pt>[r]^{p_*}
           \ar[d]_{{\lambda}_B} & C{\rtimes}_{\gamma}G
           \ar@<2pt>[l]^{\sigma_*} \ar[r]^{} 
           \ar[d]^{{\lambda}_C} &0 \\
           0 \ar[r]^{} & A{\rtimes}_{r,\alpha}G \ar[r]^{i_{*,r}} &
           B{\rtimes}_{r,\beta}G \ar@<2pt>[r]^{p_{*,r}} &
           C{\rtimes}_{r,\gamma}G \ar@<2pt>[l]^{\sigma_{*,r}} \ar[r]^{} & 0 
         }$$
           est commutatif et comme $\lambda_B$ est surjective, on a
           $$B{\rtimes}_{r,\beta}G=i_{*,r}(A{\rtimes}_{r,\alpha}G)\bm\oplus {\sigma}_{*,r}(C{\rtimes}_{r,\gamma}G).$$

 \end{proof}

\section{K-th{\'e}orie de $C^*$-alg{\`e}bres} 
La $K$-th{\'e}orie de $C^*$-alg{\`e}bres est un domaine qui s'est beaucoup
developp{\'e} depuis une trentaine d'ann{\'e}es. Il y a plusieurs livres
consacr{\'e}s {\`a} la $K$-th{\'e}orie, par exemple on pourrait consulter les
livres \cite{bl}, \cite{ka} et \cite{wo} pour se familiariser avec
cette th{\'e}orie. 

 La $K$-th{\'e}orie $K_0$ et $K_1$ est un foncteur covariant de la
  cat{\'e}gorie des $C^*$-alg{\`e}bres dans la cat{\'e}gorie des groupes
  commutatifs.

Ce n'est pas le but de ce
travail de la pr{\'e}senter. Signalons ici les deux propri{\'e}t{\'e}s {\'e}l{\'e}mentaires, qui nous seront utiles, sans du tout rentrer
dans sa d{\'e}finition. 

\begin{itemize}  

\item Si $A$ est une $C^*$-alg{\`e}bre qui est une limite inductive des
  $C^*$-alg{\`e}bres \ie $A=\varinjlim A_i$, alors
  $K_k(A)$, $k=0,1$ est une limite inductive des groupes $K_0,\,K_1$
  correspodants, \ie $K_k(A)=\varinjlim K_k(A_i)$
  pour $k=0,1$.

\item Si $0\to
  A\to B\ \mathop{\rightleftarrows}\limits^{}_{}C\to
  0$ est une suite exacte scind{\'e}e de $C^*$-alg{\`e}bres, alors $0\to
  K_k(A)\to K_k(B)\ \mathop{\rightleftarrows}\limits^{}_{} K_k(C)\to
  0$ pour $k=0,1$, est une suite exacte scind{\'e}e de groupes.  

\end{itemize}

\newpage
\quad \thispagestyle{empty}

\newpage

\chapter{$C^*$-alg{\`e}bres gradu{\'e}es}

 \section{Semi-treillis}

Soit $(\cL,\le)$ un ensemble ordonn{\'e}. On dit que $(\cL,\le)$ est  un \emph{semi-treillis}  si pour tout $k,\ell\in \cL$, l'ensemble $\{m\in \cL;\ m\le k$ et $m\le \ell\}$ poss{\`e}de un plus grand {\'e}l{\'e}ment not{\'e} $k\wedge \ell$. On dira parfois (improprement) que $\cL$ est  un semi-treillis.

Soit $(\cL,\le)$ un semi-treillis. 

\begin{itemize}
\item Si $F$ est une partie finie non vide de $\cL$ l'ensemble $\{k\in \cL;\forall \ell\in F, \ k\le \ell\}$ poss{\`e}de un plus grand {\'e}l{\'e}ment not{\'e} $\wedge_{\ell \in F} \ell$. Cela se voit imm{\'e}diatement {\`a} l'aide d'une r{\'e}currence sur le nombre d'{\'e}l{\'e}ments de $F$.

\item Une partie $\cM$ de $\cL$ est appel{\'e}e un \emph{sous-semi-treillis} de $\cL$ si pour tout $k,\ell\in \cM$, on a $k\wedge \ell\in \cM$.

\item Une \emph{partie commen{\c c}ante} de $\cL$ est une partie $\cM$ de $\cL$ telle que pour tout $a\in \cM$ on ait  $\{b\in \cL;\ b\le a\}\subset \cM$. Toute partie commen{\c c}ante de $\cL$ est un sous-semi-treillis de $\cL$.

\item Un \emph{sous-semi-treillis finissant} de $\cL$ est un sous-semi-treillis $\cM$ qui soit une partie finissante \ie telle que, pour tout $a\in \cM$, on ait  $\{b\in \cL;\ a\le b\}\subset \cM$.

\item Soit $k\in{\cal L}$. L'ensemble $\cL_k=\{\ell\in \cL;\ k\le \ell\}$ est un sous-semi-treillis finissant de ${\cal L}$. Remarquons que le compl{\'e}mentaire d'une partie finissante est commen{\c c}ante. En particulier l'ensemble $\cL'_k=\cL\setminus \cL_k$ est un sous-semi-treillis de $\cL$.

\item L'ensemble des sous-semi-treillis de $\cL$ est stable par intersection. En particulier, si $\cM$ est une partie de $\cL$, il existe un plus petit sous-semi-treillis de $\cL$ contenant $\cM$; on l'appelle le sous-semi-treillis \emph{engendr{\'e} par $\cM$}. Il est en fait facile de caract{\'e}riser ce sous-semi-treillis: c'est l'ensemble des $\wedge _{a\in F} a$ pour $F$ parcourrant l'ensemble des parties finies non-vides de $\cM$.

\item Toute partie finie $\cM$ de $\cL$ est contenue dans un
  sous-semi-treillis fini de $\cL$. En effet, il est clair que le
  sous-semi-treillis engendr{\'e} par $\cM$ est fini.

\item On dit que $\cL$ est un \emph{bon} semi-treillis si toute partie
  totalement ordonn{\'e}e non vide de $\cL$ est bien ordonn{\'e}e (\ie tout
  sous-ensemble non vide de cette partie admet un plus petit {\'e}l{\'e}ment). 
\end{itemize}

\begin{remark}\label{prtr}
Soit $(\cL^k)_{k=1}^{n}$, $n\in \N$ une famille finie de semi-treillis. Le
  produit $\cL=\prod\limits_{k=1}^n\cL^k$ est muni de l'ordre
  :
 $$\mbox{pour tout}\,a=(a^k),\,b=(b^k)\in \cL\,\mbox{alors},\,a\leq b\Leftrightarrow a^k\leq
  b^k\,\mbox{pour tout}\, k.$$
De plus, pour tout $a,b\in \cL$ on a $a\wedge b=(a^k\wedge b^k)_{k=1}^n$, donc $\cL$ est un semi-treillis. 
\end{remark}

\begin{remark}\label{bon}Soit $\cL$ un semi-treillis. On a les equivalences suivantes:
\begin{enumerate}
\item L'ensemble $\cL$ est un bon semi-treillis.
\item Toute partie totalement ordonn{\'e}e non vide de $\cL$ admet un
  plus petit {\'e}l{\'e}ment.
\item Tout sous-semi-treillis non vide de $\cL$ admet un plus petit {\'e}l{\'e}ment.
\item Tout sous-semi-treillis finissant non vide de $\cL$ admet un plus petit {\'e}l{\'e}ment.
\end{enumerate}

C'est evident que $a)\Leftrightarrow b)$.

Montrons $b)\Rightarrow c)$. Soit $\cF$ un sous-semi-treillis non vide
de $\cL$. Par le lemme de Zorn il existe une partie $A$ totalement ordonn{\'e}e maximale de
$\cF$. Notons $a$ son plus petit {\'e}l{\'e}ment et prenons $x\in
\cF$. Puisque $\cF$ est un sous-semi-treillis $a\wedge x\in
\cF$. Comme l'ensemble $\{a\wedge x\}\cup A$ est totalement ordonn{\'e}e et
$A$ est maximal alors $a\wedge x\in A$ et donc $a\leq a\wedge
x$. Cela implique $a\wedge x=a$. D'autre part, on a $a\wedge x\leq x$ et par cons{\'e}quent $a\leq
x$. En d'autres termes $a$ est le plus petit {\'e}l{\'e}ment de $\cF$. 

Puisqu'une partie totalement ordonn{\'e}e de $\cL$ est un
sous-semi-treillis de $\cL$ on a $c)\Rightarrow b)$.

C'est evident que $c)\Rightarrow d)$. 

Pour montrer que $d)\Rightarrow c)$, on consid{\`e}re $I$ un sous-semi-treillis de $\cL$ et soit $J=\{j\in \cL;\exists
i\in I; i\leq j\}$. Alors $J$ est un sous-semi-treillis finissant de
$\cL$. Par hypoth{\`e}se $J$ admet un plus petit {\'e}l{\'e}ment, notons-le
$j_0$. Alors puisque $j_0\in J$ il existe $i_0\in I$ tel que $i_0\leq j_0$. Mais puisque
par construction $I\subset J$, $i_0\in J$. Donc $j_0\leq i_0$ et on a
l'{\'e}galit{\'e}, autrement dit $i_0$ est le plus petit {\'e}l{\'e}ment de $I$.    
\end{remark}

\section{$C^*$-alg{\`e}bres gradu{\'e}es}

\paragraph{Rappels}
\begin{itemize}
\item Soit $E$ un espace vectoriel. On dit qu'une famille $(E_i)_{i\in
    I}$ de sous-espaces vectoriels de $E$ est \emph{en somme directe}
  (ou lin{\'e}airement ind{\'e}pendante) si pour toute famille
    ${(S_i)}_{i\in I}$ telle que $S_i\in
E_i,\,i\in I$ et $S_i\neq 0$ pour seulement un nombre fini au
maximum de $i$, alors si $\sum\limits_{i\in I} S_i=0$ on a $S_i=0$
    pour tout $i\in I.$\\

\item Soit $E$ un espace vectoriel norm{\'e}. Une partie de $E$ est dite \emph{totale} si le sous-espace vectoriel de $E$ qu'elle engendre est dense dans $E$. Une famille $(T_i)_{i\in I}$ de parties de $E$ est dite totale si $\bigcup\limits_{i\in I} T_i$ est une partie totale de $E$.
\end{itemize}
\bigskip Soit $(\cL,\le)$ un semi-treillis.

\begin{definition} Une $C^*$-alg{\`e}bre $\gA$ est dite \emph{gradu{\'e}e} par
  $\cL$ ou \emph{$\cL$-gradu{\'e}e}  si l'on s'est donn{\'e} une famille
  lin{\'e}airement ind{\'e}pendante et totale $(A_i)_{i\in \cL} $ de
  sous-$C^*$-alg{\`e}bres de $\gA$ telles que $A_i A_j\subset
  A_{i\wedge j}$ pour tout $i,j\in \cal L$.
\end {definition}

On dira alors que $(\gA,(A_i)_{i\in \cL})$ est une $C^*$-alg{\`e}bre
gradu{\'e}e et on appellera les $A_i$, $i\in \cL$ les \emph{composantes} de $\gA$.

\begin{definition}\label{sg} Soit $\gB$ une sous-$C^*$-alg{\`e}bre de $\gA$. Posons
  $B_i=A_i\cap \gB$. On dit que $\gB$ est une sous-alg{\`e}bre $\cL$-gradu{\'e}e si
  $\overline {\bigoplus\limits_{i\in \cL} B_i}=\gB$.

\end{definition}

Soit $(\gA,(A_i)_{i\in \cL})$ une $C^*$-alg{\`e}bre gradu{\'e}e. Pour toute
partie $\cM$ de $\cL$ notons $\gA_\cM$ le sous-espace de $\cL$ somme
de la famille $A_i$ pour $i\in \cM$. 

On a imm{\'e}diatement:

\begin{proposition} \label{ideal}\begin{enumerate}
\renewcommand\theenumi{\alph{enumi}}
\item Si $\cM$ est un sous-semi-treillis de $\cL$, l'ensemble $\overline{\gA_\cM}$ est une sous-$C^*$-alg{\`e}bre de $\gA$. Elle est gradu{\'e}e par $\cM$.

\item Si $\cM$ est une partie commen{\c c}ante de $\cL$, l'ensemble
  $\overline{\gA_\cM}$ est un id{\'e}al ferm{\'e} de $\gA$.
\end{enumerate}
\indent\hfill$\square$
\end {proposition}

Le r{\'e}sultat suivant sera tr{\`e}s utile dans toute la suite:

\begin{proposition}\label{ferm{\'e}} Si $\cF$ est un sous-semi-treillis \emph{fini} de $\cL$, l'ensemble $\gA_\cF$ est ferm{\'e} dans $\gA$. C'est une sous-$C^*$-alg{\`e}bre de $\gA$.
\end {proposition}

\begin{proof}On va le d{\'e}montrer par r{\'e}currence sur le nombre
  d'{\'e}l{\'e}ments de $\cF$.

 Si $\cal F$ poss{\`e}de un seul {\'e}l{\'e}ment ${\cal F}=\{i\}$
alors $A_i$ est ferm{\'e} par hypoth{\`e}se.

 Supposons que c'est vrai pour tout sous-semi-treillis
 poss{\`e}dent un nombre d'{\'e}l{\'e}ments inf{\'e}rieur ou {\'e}gal {\`a} $n-1$.
 On va le montrer pour $\cal F$ {\`a} $n$ {\'e}l{\'e}ments.

Soit $i$ un {\'e}l{\'e}ment maximal de $\cal F$ et soit $S={\cal
  F}\backslash \{i\}$. L'ensemble $S$ est un sous-semi-treillis de $\cal F$ {\`a} $n-1$
{\'e}l{\'e}ments.

On a $\overline{\gA_\cF}=\overline{\gA_S\bm\oplus A_i}$ et on remarque que $S$ est une partie
commen{\c c}ante de $\cF$, donc par la proposition pr{\'e}cedente
$\overline {\gA_S}$ est un id{\'e}al ferm{\'e} de $\overline {\gA_\cF}$. Par
l'hypoth{\`e}se de r{\'e}currence $\gA_S$ est un id{\'e}al ferm{\'e} de
$\overline{\gA_S\bm\oplus A_i}$. 

L'application $A_i\to
\overline{\gA_S\bm\oplus A_i}/\gA_S$ est un isomorphisme de
$C^*$-alg{\`e}bres puisque son noyau est nul et son image est dense. 

Soit $x\in
\overline{\gA_S\bm\oplus A_i}$, il existe $y\in A_i$ tel que la classe de
$x$ modulo $\gA_S$ soit {\'e}gale {\`a} la classe de $y$ modulo $\gA_S$. Alors
$x-y\in \gA_S\Rightarrow x\in y+\gA_S\Rightarrow x\in \gA_S\bm\oplus A_i$
et $\gA_\cF$ est ferm{\'e} dans $\gA$.

\end{proof}

Notons $\F_\cL$ l'ensemble des sous-semi-treillis finis de $\cL$. \\
Si $\cF_1,\,\cF_2\in \F_\cL$ alors il existe $\cF_3\in \F_\cL$ tel que
$\cF_1\subset \cF_3$ et $\cF_2\subset \cF_3$. Par d{\'e}finition pour
$\cF,\,\cM\in \F_\cL$ tels que $\cM\subset \cF$ on a $\gA_{\cM}\subset \gA_{\cF}$. En d'autres termes les ${(\gA_\cF)}_{\cF\in \F_\cL}$
forment un syst{\`e}me inductif de sous-$C^*$-alg{\`e}bres de $\gA$.\\

Puisque $\bigcup\limits_{\cF\in \F_\cL}\gA_\cF=\sum\limits_{i\in
 \cL}A_i=\gA_\cL$ contient les $A_i,\,i\in \cL$ et ${(A_i)}_{i\in
 \cL}$ est une famille totale alors $\gA_\cL$ est dense
 dans $\gA$. On en deduit que la $C^*$-alg{\`e}bre $\gA$ est la limite
 inductive des ${(\gA_\cF)}_{\cF\in \F_\cL}$.\\

\begin{proposition}\label{inter}
Soient $\cL_1$ et $\cL_2$ des semi-treillis et $(\gA,(A_i)_{i\in
  \cL_1})$ (resp. $(\gA,(B_j)_{j\in \cL_2})$) une $C^*$-alg{\`e}bre
  $\cL_1$-gradu{\'e}e (resp. $\cL_2$-gradu{\'e}e). Supposons que
la somme $\sum\limits_{(i,j)\in \cL_1\times \cL_2}(A_i\cap B_j)$ est
dense dans $\gA$. Alors, $\gA$ est aussi une $C^*$-alg{\`e}bre $\cL_1\times
  \cL_2$-gradu{\'e}e  dont les composantes sont les $A_i\cap B_j$, $i\in
  \cL_1,j\in \cL_2$.
\end{proposition}
 
\begin{proof}Montrons que la famille $(A_i\cap B_j)_{(i,j)\in
  \cL_1\times \cL_2}$ est lin{\'e}airement ind{\'e}pendente. Soit
  $(S_{ij})_{(i,j)\in \cL_1\times \cL_2}$ une famille telle que
  $S_{ij}\in (A_i\cap B_j)$ pour tout $i\in \cL_1,j\in \cL_2$ et
  $\sum\limits_{(i,j)\in \cL_1\times \cL_2}S_{ij}=0$. Alors,
  $\sum\limits_{i\in \cL_1}\sum\limits_{j\in \cL_2}S_{ij}=0$ et
  puisque $S_{ij}\in A_i$ et la famille $(A_i)_{i\in \cL_1}$ est lin{\'e}airement ind{\'e}pendent
  on a $\sum\limits_{j\in \cL_2}S_{ij}=0$ pour tout $i\in \cL_1$. Mais puisque $S_{ij}\in B_j$ et la famille $(B_j)_{j\in \cL_2}$ est
  lin{\'e}airement ind{\'e}pendente, donc $S_{ij}=0$ pour tout $j\in \cL_2$ et
  $i\in \cL_1$. 

Par hypoth{\`e}se elle est aussi totale, donc il suffit de montrer que
$(A_i\cap B_j)\cdot (A_\ell\cap B_k)\subset (A_{i\wedge \ell}\cap
B_{j\wedge k})$ pour tout $i,\ell\in \cL_1$ et $j,k\in \cL_2$. Pour
cela, soient $S_{ij}\in (A_i\cap B_j)$ et $D_{\ell k}\in (A_\ell\cap
B_k)$. Puisque $S_{ij}\in A_i$, $D_{\ell k}\in A_\ell$ et $(\gA,(A_i)_{i\in
  \cL_1})$ est gradu{\'e}e, on a $S_{ij} D_{\ell k}\in A_{i\wedge
  \ell}$. On a {\'e}galement $S_{ij} D_{\ell k}\in B_{j\wedge k}$,
d'o{\`u} le r{\'e}sultat.

\end{proof}

\section{Morphismes de $C^*$-alg{\`e}bres gradu{\'e}es}

\begin{definition}Soit $\cL$ un semi-treillis et soient
 $(\gA,(A_i)_{i\in \cL})$ et $(\gB,(B_i)_{i\in \cL})$ des
 $C^*$-alg{\`e}bres $\cL$-gradu{\'e}es. On dit qu'un homomorphisme
 $\psi:\gA\to \gB$ est un \emph{homomorphisme de $C^*$-alg{\`e}bres gradu{\'e}es}
 si $\psi (A_i)\subset B_i$ pour tout $i\in \cL$.
\end{definition}

\begin{proposition}\label{injectif} Soient $(\gA,(A_i)_{i\in
    \cL})$ une $C^*$-alg{\`e}bre gradu{\'e}e, $B$ une $C^*$-alg{\`e}bre et
    $\psi:\gA\to B$ un homomorphisme. Pour tout sous-semi-treillis
    fini $\cF$ de $\cL$, notons $\psi_\cF:\gA_\cF\to B$ la
    restriction de $\psi$ {\`a} $\gA_\cF$. Pour $i\in \cL$ on {\'e}crit
    $\psi_i:A_i\to B$ plut{\^o}t que $\psi_{\{i\}}$. 

\begin{enumerate}
\item[$a)$] Pour tout $j,k \in \cL$ et tout $x\in A_j$, $y\in
  A_k$ on a $\psi _{j\wedge k}(xy)=\psi_j(x)\psi_k (y)$.

\item[$b)$] Les conditions suivantes sont {\'e}quivalentes:
\begin{enumerate}
\renewcommand{\theenumii}{\roman{enumii}}
\renewcommand{\labelenumii}{\rm (\theenumii)}
\item L'application $\psi$ est surjective.

\item La famille $\psi_i(A_i)$ est totale.
\end{enumerate}

\item[$c)$] Les conditions suivantes sont {\'e}quivalentes:
\begin{enumerate}
\renewcommand{\theenumii}{\roman{enumii}}
\renewcommand{\labelenumii}{\rm (\theenumii)}
\item L'application $\psi$ est injective.

\item Pour tout sous-semi-treillis fini $\cF$ de $\cL$,  $\psi_\cF$ est injective.

\item Pour tout $i\in \cL$ l'application $\psi_i $ est injective et la famille $\psi_i (A_i)$ de sous-espaces de $B$ est en somme directe.

\end{enumerate}
\end{enumerate}
\end{proposition}

\begin{proof}\begin{enumerate}
\item[$a)$] est clair.

\item[$b)$] Puisque l'image d'un homomorphisme de $C^*$-alg{\`e}bres est ferm{\'e} et que la
  famille ${(A_i)}_{i\in \cL}$ est totale dans $\gA$, l'image de
  $\psi$ est le sous-espace ferm{\'e} de $B$ engendr{\'e} par les
  $\psi_i(A_i)$.

\item[$c)$] Remarquons que la condition $(iii)$ {\'e}quivaut {\`a} : $(iii)'$
  la restriction $\psi_\cL$ de $\psi$ {\`a} $\gA_\cL$ est
  injective. Il est alors clair que $(i)\Rightarrow (iii)\Rightarrow
  (ii)$.\\
Supposons que $(ii)$ soit satisfait. Soit $x\in \gA_\cL$. Il existe un
  sous-semi-treillis fini $\cF\subset \cL$ tel que $x\in
  \gA_\cF\subset \gA_\cL$. Le morphisme $\psi_\cF$ de
  $C^*$-alg{\`e}bres {\'e}tant injectif, il est isom{\'e}trique. On a donc
  $||\psi(x)||=||x||$ pour tout $x\in \gA_\cL$ et, par densit{\'e},
  pour tout $x\in \gA$.  
\end{enumerate}
\end{proof}

\begin{proposition}\label{rest}Soient $(\gA,(A_i)_{i\in
    \cL})$ une $C^*$-alg{\`e}bre gradu{\'e}e et $B$ une
    $C^*$-alg{\`e}bre. Donnons-nous pour tout $i\in \cL$ un homomorphisme
    $\psi_i:A_i\to B$. On suppose que pour tout $j,k \in \cL$ et tout $x\in A_j$, $y\in
  A_k$ on a $\psi _{j\wedge k}(xy)=\psi_j(x)\psi_k (y)$. Alors il
    existe un unique homomorphisme $\psi :\gA\to B$ dont
  la restriction {\`a} $A_i$ soit $\psi_i$ pour tout $i\in
  \cL$.

\end{proposition}

\begin{proof}Puisque les $A_i$ sont ind{\'e}pendants, il existe une
    unique application lin{\'e}aire $\psi_\cL:\gA_\cL\to B$ dont la
    restriction aux $A_i$ soit $\psi_i$. L'application
    $\psi_\cL$ est un homomorphisme d'alg{\`e}bres involutives. Soit
    $x\in \gA_\cL$. Il existe un sous-semi-treillis fini $\cF\subset
    \cL$ tel que $x\in \gA_\cF\subset \gA_\cL$. La restriction de
    $\psi_\cL$ {\`a} $\gA_\cF$ est un homomorphisme de $C^*$-alg{\`e}bres
    $\psi_\cF:\gA_\cF\to B$. On a donc
    $||\psi_\cL(x)||=||\psi_\cF(x)||\leq ||x||$. On en d{\'e}duit
    que $\psi_\cL$ admet un unique prolongement {\`a} l'adh{\'e}rence $\gA$
    de $\gA_\cL$.

\end{proof}

On remplace maintenant la $C^*$-alg{\`e}bre $B$ dans les propositions
\ref{injectif}, \ref{rest} par une $C^*$- alg{\`e}bre $\cL$-gradu{\'e}e
$({\frak B},{(B_i)}_{i\in \cal L})$, donc on a

\begin{corollary}\label{morgrad}Soit ${\cal L}$ un semi-treillis et soient $({\frak A},{(A_i)}_{i\in
  \cal L})$ et $({\frak B},{(B_i)}_{i\in \cal L})$ deux
  $C^*$-alg{\`e}bres $\cal L$-gradu{\'e}es. Soit ${\rho}_i:A_i\to B_i$ un morphisme
  tel que ${\rho}_i(x) {\rho}_j(y)={\rho}_{i\wedge j}(xy)$ pour tout
  $i,j\in \cL$ et pour tout $x\in A_i,y\in A_j$. Alors il existe un
  unique morphisme $\rho:{\frak A}\to {\frak B}$ qui soit un morphisme
  de $C^*$-alg{\`e}bres gradu{\'e}es. De plus si pour tout $i\in {\cal L}$ le
  morphisme $\rho_i$ est injectif
$($resp. surjectif$)$ alors $\rho$ est injectif $($resp. surjectif$)$.

  \end{corollary}

\begin{proof}Cela r{\'e}sulte de \ref{injectif} et \ref{rest}.

 \end{proof}

\section{Sous-semi-treillis finissants et suites exactes scind{\'e}es}

\begin{proposition}\label{suitexacte} Soit $\cL$ un semi-treillis et
  $\cM$ un sous-semi-treillis finissant de $\cL$; notons $\cM'$ son
  compl{\'e}mentaire dans $\cL$. Il existe un unique homomorphisme
  $p:\gA\to \overline {\gA_{\cM}}$ dont la restriction {\`a} $A_i$ soit
  {\'e}gale {\`a} $\mbox{Id}_{A_i}$ si $i\in \cM$ et $0$ si $i\in \cM'$. On a
  $\ker p=\overline {\gA_{\cM'}}$. En d'autres termes, on a
  $\gA=\overline{\gA_{\cM'}} \bm\oplus \overline{\gA_\cM}$ et une suite
  exacte scind{\'e}e $$0\to \overline {\gA_{\cM'}}\xrightarrow{i} \gA \
  \mathop{\rightleftarrows}\limits^p_\sigma\  \overline {\gA_{\cM}}
  \to 0,$$
o{\`u} $i$ et $\sigma$ sont les inclusions naturelles.
\end {proposition}

\begin{proof}  Soit $B_i=\left\{\begin{aligned}&A_i\;\mbox{si}\;i\in
      \cM\\&0\;\;\mbox{si}\;i\in \cM'\end{aligned}\right.$ alors
      l'alg{\`e}bre $(\overline{\gA_\cM},(B_i)_{i\in \cL})$ est une $C^*$-alg{\`e}bre
    $\cL$ gradu{\'e}e. 

Soit $p_i:A_i\to B_i$ le morphisme d{\'e}fini par $p_i=\left\{\begin{aligned}&\mbox{Id}_{A_i}\;\mbox{si}\;i\in
      \cM\\&0\;\;\mbox{si}\;i\in \cM'.\end{aligned}\right.$ 

D'apr{\`e}s le
      corollaire \ref{morgrad} il existe un unique homomorphisme
      $p:\gA\to \overline{\gA_\cM}$ tel que $p_{|_{A_i}}=p_i$
      pour tout $i\in \cL$. Alors $p|_{\overline {\gA_{\cM'}}}=0$
      et $p|_{\overline {\gA_{\cM}}}=\mbox{Id}_{\overline
        {\gA_{\cM}}}$ de sorte que $\overline{\gA_{\cM'}} \bm\cap
      \overline{\gA_\cM}=\{0\}$.

D'autre part, $\overline{\gA_{\cM'}}$ est un id{\'e}al ferm{\'e} de $\gA$ et $\overline{\gA_\cM}$ est une
      sous-$C^*$-alg{\`e}bre de $\gA$ (proposition \ref{ideal}), donc $\overline{\gA_{\cM'}} +
      \overline{\gA_\cM}$ est ferm{\'e} (proposition \ref{idssa}). Comme $\overline{\gA_{\cM'}} + \overline{\gA_\cM}$ contient $\gA_\cL$ qui est dense dans $\gA$
      on a $\gA=\overline{\gA_{\cM'}} + \overline{\gA_\cM}$, autrement
      dit on a $\gA=\overline{\gA_{\cM'}} \bm\oplus \overline{\gA_\cM}$. 

\end{proof}

En particulier, l'ensemble $\cL_a=\{b\in \cL;\ a\le b\}$ {\'e}tant un sous-semi-treillis finissant, on obtient une suite exacte scind{\'e}e $$0\to \overline {\gA_{\cL'_a}}\rightarrow \gA \ \mathop{\rightleftarrows}\limits^{p_a}_{\sigma_a}\  \overline {\gA_{\cL_a}} \to 0$$
o{\`u} l'on a not{\'e} $\cL'_a=\{b\in \cL;\ a\not\le b\}$ le compl{\'e}mentaire de
$\cL_a$.\\

\begin{corollary}\label{surinv}Soit ${\cal L}$ un semi-treillis et soient $({\frak A},{(A_i)}_{i\in
  \cal L})$ et $({\frak B},{(B_i)}_{i\in \cal L})$ des
  $C^*$-alg{\`e}bres $\cal L$-gradu{\'e}es.
 Soit $\rho:{\frak A}\to {\frak B}$ un morphisme de
  $C^*$-alg{\`e}bres gradu{\'e}es. Notons $\rho_i$ sa restriction {\`a} $A_i$ pour tout
  $i\in \cL$. Si $\rho$ est surjectif alors $\rho_i$ est surjectif pour tout $i\in \cL$.
 
 \end{corollary}

\begin{proof} Supposons que $\rho $ soit surjective et soit $i\in
  \cL$. Par la prop. \ref{suitexacte}, on a
  $\gA=\overline{\gA_{\cL_i}}\bm\oplus \overline{\gA_{\cL'_i}}$ de sorte
  que $\gB=\rho(\gA)=\rho \big(\overline{\gA_{\cL_i}}\big)+\rho
  \big(\overline{\gA_{\cL'_i}}\big)$. Remarquons que l'on a $\rho
  \big(\overline{\gA_{\cL_i}}\big)\subset \overline{\gB_{\cL_i}}$ et
  $\rho \big(\overline{\gA_{\cL'_i}}\big)\subset
  \overline{\gB_{\cL'_i}}$  Comme de plus
  $\gB=\overline{\gB_{\cL_i}}\bm \oplus \overline{\gB_{\cL'_i}}$, il vient $\rho \big(\overline{\gA_{\cL_i}}\big)=\overline{\gB_{\cL_i}}$ et $\rho \big(\overline{\gA_{\cL'_i}}\big)=\overline{\gB_{\cL'_i}}$.

Soit alors $b\in B_i\subset \overline{\gB_{\cL_i}}$. Il existe $c\in \overline {\gA_{\cL_i}}$ tel que $\rho(c)=b$. Soit $\varepsilon >0$. Il existe un sous-semi-treillis fini $\cF$ de $\cL_i$ et $x\in \gA_\cF$ tel que $\|x-c\|<\varepsilon$. 

Posons $T=\cF\setminus\{i\}$ et, pour $j\in T$, notons $p_j:\gB\to \overline{\gB_{\cL_j}}$ l'homomorphisme associ{\'e} et $P:\gB\to \displaystyle\prod_{j\in T}\overline{\gB_{\cL_j}}$ l'homomorphisme $y\mapsto (p_j(y))_{j\in T}$. On a $P(b)=0$, de sorte que $\|P(\rho(x))\|=\|P(\rho(x-c))\|<\varepsilon$. 

\begin{lemma} On a $\ker P\cap \gB_\cF\subset B_i$.
\end{lemma}
\begin{proof} On raisonne par r{\'e}currence sur le nombre d'{\'e}l{\'e}ments de $\cF$. Si   $\cF\subset \{i\}$  il n'y a rien {\`a} montrer. Sinon, soit $j$ un {\'e}l{\'e}ment maximal de $T$ et posons $\cF_1=\cF\setminus \{j\}$. Soit $y=\sum_{k\in \cF} y_k\in \ker P\cap \gB_\cF$ . Puisque $y_j=p_j(y)=0$, on a $y= \sum_{k\in \cF_1} y_k$, de sorte que par hypoth{\`e}se de r{\'e}currence $y\in B_i$.
\end{proof}

Comme les $\rho (A_k)$ sont contenus dans $B_k$ qui sont en somme
directe, il vient $\rho(\gA_\cF)\cap  \ker P=\rho(\gA_\cF)\cap
\gB_\cF\cap  \ker P\subset \rho(\gA_\cF)\cap B_i\subset \rho(A_i)$. Par
 la proposition \ref{isom} appliqu{\'e}e {\`a} $P:\rho(\gA_\cF)\to
\displaystyle\prod_{j\in T}\overline{\gB_{\cL_j}}$, il
existe $h\in \rho(\gA_\cF)$ tel que
$P(h)=P(\rho(x))$ et $\|h\|=\|P(\rho(x))\|<\varepsilon$. Posons
$z=\rho(x)-h$. On a $z\in \rho(\gA_\cF)\cap  \ker P\subset \rho(A_i)$
et $\|\rho(x)-z\|<\varepsilon$. Alors $\|z-b\|<2\varepsilon$. Comme $\rho (A_i)$ est ferm{\'e}e, on trouve $b\in \rho(A_i)$.
\end{proof}

Par la proposition \ref{injectif} et le corollaire pr{\'e}cedent on a
alors

\begin{corollary}Soit ${\cal L}$ un semi-treillis et soient $({\frak A},{(A_i)}_{i\in
  \cal L})$ et $({\frak B},{(B_i)}_{i\in \cal L})$ des
  $C^*$-alg{\`e}bres $\cal L$-gradu{\'e}es.
 Soit $\rho:{\frak A}\to {\frak B}$ un morphisme de
  $C^*$-alg{\`e}bres gradu{\'e}es. Notons $\rho_i$ sa restriction {\`a} $A_i$ pour tout
  $i\in \cL$. Alors $\rho$ est surjectif (resp. injectif) si et
  seulement si $\rho_i$ est surjectif (resp. injectif) pour tout $i\in
  \cL$.\newline\indent \hfill$\square$ 
\end{corollary}

On en deduit aussi 

\begin{corollary}\label{sgr}Soient ${\cal L}$ un semi-treillis, $\gA$
  une $C^*$-alg{\`e}bre et $(A_i)_{i\in
  \cal L}$, $(B_i)_{i\in \cal L}$ des familles de sous-$C^*$-alg{\`e}bres
de $\gA$ telles que $(A_i)$ soit totale, $(B_i)$
soit lin{\'e}airement ind{\'e}pendente, $A_iA_j\subset A_{i\wedge j}$ pour tout $i,j\in \cL$,
$B_iB_j\subset B_{i\wedge j}$ pour tout $i,j\in \cL$ et $A_i\subset
B_i$ pour tout $i\in \cL$. Alors $({\frak A},{(A_i)}_{i\in
  \cal L})$ et $({\frak A},{(B_i)}_{i\in \cal L})$ sont $\cal
L$-gradu{\'e}es et $A_i=B_i$ pour tout $i\in \cL$.
\end{corollary}

\begin{proof}C'est clair par la d{\'e}finition d'une $C^*$-alg{\`e}bre $\cL$-gradu{\'e}e
  que $\gA$ est une
  $C^*$-alg{\`e}bre $\cL$-gradu{\'e}e dont les composantes sont les $A_i$ et
  les $B_i$. L'{\'e}galit{\'e} de ces deux composantes r{\'e}sulte du corollaire
 \ref{surinv}.\\
\end{proof}

\begin{remark} Soit $\cM$ un sous-semi-treillis de $\cL$. L'alg{\`e}bre $\overline{\gA_\cM}$ est une
sous-$C^*$-alg{\`e}bre de $\gA$, $\cL$-gradu{\'e}e et on a
$\overline{\gA_\cM}\bm\cap A_i=\left\{\begin{aligned}A_i\;\mbox{si}\;i\in
      \cM\\0\;\mbox{si}\;i\notin \cM.\end{aligned}\right.$

En effet, c'est clair que si $i\in \cM$ on a $\overline{\gA_\cM}\bm\cap
A_i=A_i$.

Supposons que $i\notin \cM$. Soit $\cL_i=\{j\in \cL;\ i\le j\}$ et $\cL'_i=\{j\in \cL;\
i\not\le j\}$. Alors on a une suite exacte scind{\'e}e $$0\to \overline
{\gA_{\cL'_i}}\rightarrow \gA \
\mathop{\rightleftarrows}\limits^{p_i}_{\sigma_i}\  \overline
{\gA_{\cL_i}} \to 0.$$
Prenons $x\in \overline{\gA_\cM}\bm\cap A_i$ non nul. Comme $x\in A_i$ on a
$p_i(x)=x$.

L'ensemble $\cM\bm\cap \cL_i$ est un sous-semi-treillis finissant de
$\cM$ et son compl{\'e}mentaire dans $\cM$ est l'ensemble $\cM\bm\cap
\cL'_i$. En appliquant la proposition pr{\'e}cedente on obtient une suite
exacte scind{\'e}e
 $$0\to {\overline
{\gA}}_{\cM\bm\cap \cL'_i}\rightarrow \overline{\gA_\cM} \
\mathop{\rightleftarrows}\limits^{q_i}_{\sigma_i}\  {\overline
{\gA}}_{\cM\bm\cap \cL_i} \to 0.$$
Or, $q_i(x)=x$ si $x\in A_j$; $j\in \cM\bm\cap \cL_i$ et $q_i(x)=0$ si
$x\in A_j$; $j\in \cM\bm\cap \cL'_i$, donc $q_i$ est la restriction de
$p_i$ {\`a} $\overline{\gA_\cM}$. Puisque $x\in \overline{\gA_\cM}$ et on
a $x=p_i(x)$ alors $x\in p_i(\overline{\gA_\cM})=q_i(\overline{\gA_\cM})={\overline
{\gA}}_{\cM\bm\cap \cL_i}$. 

Comme $\cM\bm\cap \cL_i$ est un
sous-semi-treillis finissant de $(\cM\bm\cap \cL_i)\bm\cup \{i\}$ et $\{i\}$
est son compl{\'e}mentaire dans $(\cM\bm\cap \cL_i)\bm\cup \{i\}$, on a aussi
${\overline{\gA}}_{(\cM\bm\cap \cL_i)\bm\cup \{i\}}={\overline{\gA}}_{\cM\bm\cap
      \cL_i}\bm\oplus A_i$. Puisque $x\in {\overline{\gA}}_{\cM\bm\cap \cL_i}\bm\cap A_i$, on en deduit que $x=0$.\\  
\end{remark}

\section{Morphismes de structure de $C^*$-alg{\`e}bres gradu{\'e}es}

Soit $\cL$ un semi-treillis et $(\frak A,{(A_a)}_{a\in \cal L})$ une
  $C^*$-alg{\`e}bre $\cal L$-gradu{\'e}e. Comme $\{a\}$ est une partie commen{\c c}ante de $\cL_a$, $A_a$ est un id{\'e}al ferm{\'e} dans
  $\overline{\gA_{\cL_a}}$. On obtient donc un morphisme
  $\varphi_a:\overline{\gA_{\cL_a}}\to M(A_a)$ et donc un morphisme
  $\pi_a=\varphi_a\circ p_a:\gA\to M(A_a)$.

Soient $a,b\in \cL$ tels que $a\le b$. Notons $\varphi_{a,b}:A_b\to
M(A_a)$ la restriction {\`a} $A_b$ de $\varphi_a:\overline{\gA_{\cL_a}}\to M(A_a)$. Les morphismes
$\varphi_{i,j}$ s'appellent les \emph{morphismes de structure} de
l'alg{\`e}bre gradu{\'e}e $(\gA,(A_i)_{i\in \cL})$.

\begin{proposition}\label{q}Soit $(\frak A,{(A_i)}_{i\in \cal L})$ une
  $C^*$-alg{\`e}bre $\cal L$-gradu{\'e}e. Alors,
\begin{enumerate}
\item[$(i)$] pour tout $i,j\in \cal L$ tels que $i\leq j$, il existe un unique
morphisme \\${\varphi}_{i,j}:A_j\to M(A_i)$ satisfaisant
${\varphi}_{i,j}(y)x=y x$, $x\in
A_i,y\in A_j$. Les
morphismes ${\varphi}_{i,j}$ v{\'e}rifient$:$
\begin{enumerate}

 \item[$a)$]${\varphi}_{i,i}={\mbox{Id}}_{A_i}$ pour tout $i\in {\cal L}$,

\item[$b)$] pour tout $i,j\in {\cal L}$ posons $k=i\wedge j$ et prenons $m\leq i\wedge j$ alors,
$${\varphi}_{k,i}(x){\varphi}_{k,j}(y)\in A_k$$ et 
 $${\varphi}_{m,k}({\varphi}_{k,i}(x){\varphi}_{k,j}(y))={\varphi}_{m,i}(x){\varphi}_{m,j}(y),$$
pour tout
  $x\in A_i,y\in A_j$.
\end{enumerate}
\item[$(ii)$] pour tout $i,j\in \cal L$ il existe une unique
  application \\$q_{i,j}:A_i\times A_j\to A_{i\wedge j}$ satisfaisant
 $q_{i,j}(x,y)=xy$, $x\in A_i, y\in A_j$. Les applications
$q_{i,j}$ sont bilin{\'e}aires et verifient$:$

\begin{enumerate}
\item[$a')$] $q_{i,i}(x,y)=xy$ pour tout $i\in {\cal L}$,

\item[$b')$] pour tout $i,j\in {\cal L}$ et $x\in A_i,y\in
A_j$,  $q_{i,j}(x,y)={q_{j,i}(y^*,x^*)}^*$ et

\item[$c')$] pour tout $i,j,k\in {\cal L}$ et $x\in A_i,y\in A_j,z\in A_k$
 on a $$q_{i\wedge j,k}(q_{i,j}(x,y),z)=q_{i,j\wedge
  k}(x,q_{j,k}(y,z)).$$
\end{enumerate}
\end{enumerate}
\end{proposition}

\begin{proof} \begin{enumerate}
\item[(i)] Par d{\'e}finition on a $A_i A_j\subset A_{i\wedge j}$ pour
  tout $i,j\in {\cal L}$.\\
Si $i\leq j$ alors $A_i A_j\subset A_i$ donc si $x\in A_i,y\in
A_j\Rightarrow x y\in A_i$ et on peut d{\'e}finir un
morphisme ${\varphi}_{i,j}:A_j\to M(A_i)$ par
${\varphi}_{i,j}(y)x=yx$. On va d{\'e}montrer que ${\varphi}_{i,j}$
v{\'e}rifie les conditions $a)$ et $b)$. En effet,
\begin{enumerate}

\item[$a)$] ${\varphi}_{i,i}$ est le plongement canonique de $A_i\to M(A_i)$.

\item[$b)$] Soit $i,j\in {\cal L}\;\mbox{et}\; x\in A_i,y\in
  A_j$. Posons $k=i\wedge j$ et prenons $z\in
A_k$, alors par hypoth{\'e}se on a $A_i A_j\subset A_k$. Par les
  applications ${\varphi}_{k,i}(x):z\mapsto x z$ et 
${\varphi}_{k,j}(y):z\mapsto y z$ on a,
${\varphi}_{k,i}(x)[{\varphi}_{k,j}(y)(z)]={\varphi}_{k,i}(x)(y z)=x y
  z$, donc ${\varphi}_{k,i}(x){\varphi}_{k,j}(y)={\varphi}_{k,k}(x
  y)=x y\in A_k$.

 Soit maintenant $m\leq k$ et $w\in A_m$, alors\\
${\varphi}_{m,k}({\varphi}_{k,i}(x){\varphi}_{k,j}(y))(w)={\varphi}_{m,k}(x
y)(w)=(x y) w$.\\
D'autre
part ${\varphi}_{m,i}(x)[{\varphi}_{m,j}(y)(w)]={\varphi}_{m,i}(x)[y w]=x (y w)$.\\
Puisque $w\in A_m$ est quelconque on obtient
l'{\'e}galit{\'e}:\\${\varphi}_{m,k}({\varphi}_{k,i}(x){\varphi}_{k,j}(y))={\varphi}_{m,i}(x){\varphi}_{m,i}(y).$
\end{enumerate}

\item[(ii)] Puisque $A_iA_j\subset A_{i\wedge j}$ pour tout $i,j\in {\cal L}$ on peut d{\'e}finir une application $q_{i,j}:A_i\times
A_j\to A_{i\wedge j}$ par $(x,y)\mapsto q_{i,j}(x,y)=x y$ qui est bilin{\'e}aire. On va
d{\'e}montrer qu'elle verifie les conditions $a'),\,b'),\,c')$. En
effet,
\begin{enumerate}

\item[$a')$] c'est {\'e}vident.

\item[$b')$] Soient $i,j\in \cL$ et $x\in A_i,\,y\in A_j$, par
  d{\'e}finition on a
${q_{i,j}(x,y)}^*={(x y)}^*=y^* x^*=q_{j,i}(y^*,x^*)$.

\item[$c')$] Soient $i,j,k\in \cL$ et $x\in A_i,\,y\in A_j,\,z\in A_k$
 alors,

 $\begin{aligned}q_{i\wedge j,k}(q_{i,j}(x,y),z)&=q_{i\wedge j,k}(x y,z)=x y z=q_{i,j\wedge k}(x,y z)\\&=q_{i,j\wedge
 k}(x,q_{j,k}(y,z)).\end{aligned}$
\end{enumerate}
\end{enumerate}
\end{proof}

\begin{remark} Soient $i,j\in \cL$ tels que $i\leq j$. On suppose que le
 morphisme $\varphi_{i,j}$ est non-d{\'e}g{\'e}n{\'e}r{\'e} \ie $\overline
 {\varphi_{i,j}(A_j)A_i}=\overline {A_i\varphi_{i,j}(A_j)}=A_i$. Suite
 {\`a} la proposition \ref{cohen} $\varphi_{i,j}$ admet un unique prolongement en un morphisme
 $\widetilde {\varphi_{i,j}}:M(A_j)\to M(A_i)$.

Si pour tout $a,b\in \cL,\,a\leq b$ les morphismes $\varphi_{a,b}$
 sont non-d{\'e}g{\'e}n{\'e}r{\'e}s, la condition $b)$ de la proposition pr{\'e}cedente
 est {\'e}quivalente {\`a} la condition $$\widetilde {\varphi_{a,b}}\circ
\widetilde {\varphi_{b,c}}=\widetilde {\varphi_{a,c}}\,,\,\mbox{avec}\,a\leq b\leq
c,\,a,b,c\in \cL.$$    
\end{remark}

\begin{proposition} \label{injectivite} Soit $(\gA,(A_\ell)_{\ell\in \cL})$ une $C^*$-alg{\`e}bre gradu{\'e}e.
On suppose que les morphismes de structure
$\varphi_{k,\ell}:A_\ell\to M(A_k),\;k\leq \ell$ satisfont $\varphi_{k,\ell}^{-1}(A_k)=\{0\}$, alors pour tout $i\in \cL$, l'application $\varphi_i:\overline{\gA_{\cL_i}}\to M(A_i)$ est injective.
\end {proposition}

\begin{proof} Soit $\cF$ un sous-semi-treillis fini non vide de $\cL$. Notons $\ell$ son plus petit {\'e}l{\'e}ment. Montrons par r{\'e}currence sur le nombre d'{\'e}l{\'e}ments de $\cF$ que l'application $\varphi_\ell:\gA_{\cF}\to M(A_\ell)$ est injective.

Soient $(x_k)_{k\in \cF}$ une famille d'{\'e}l{\'e}ments de $\gA_\cF$ avec $x_k\in A_k$ et $\sum_{k\in \cF}\varphi_{\ell,k}(x_k)=0$.

Si $\cF$ a deux {\'e}l{\'e}ments, $\ell, k$, avec $\ell\le k$, on a $\varphi_{\ell,k}(x_k)=-\varphi_{\ell,\ell}(x_\ell)=-x_\ell$, donc $\varphi_{\ell,k}(x_k)\in A_\ell$, d'o{\`u} l'on d{\'e}duit $x_k=0$ et $x_\ell=-\varphi_{\ell,k}(x_k)=0$.

Supposons que $\cF$ a au moins trois {\'e}l{\'e}ments et que, pour tout sous-semi-treillis $\cF'$ de $\cF$ distinct de $\cF$ de plus petit {\'e}l{\'e}ment $\ell'$, l'application $\varphi_{\ell'}:\gA_{\cF'}\to M(A_{\ell'})$ est injective.

Soit $k\in \cF$ un {\'e}l{\'e}ment tel que $k\ne \ell$ et que $k$ ne soit pas
le plus grand {\'e}l{\'e}ment de $\cF$. Posons $\cF'=\{j\in \cF;\ j\le
k\}$. Posons aussi $\cF_k=\{j\in \cF;\ k\le j\}$ et $\cF_k'=\{j\in
\cF;\ k\not \le j\}$. Remarquons que $\cF',\,\cF_k,\,\cF_k'$ sont
distincts de $\cF$.

Donnons-nous un {\'e}l{\'e}ment $\sum\limits_{j\in \cF}x_j\in \ker \varphi_\ell$
et soit $x\in A_k$. Le produit $x\sum_{j\in \cF}x_j$ s'{\'e}crit $\sum_{j\in \cF'} y_j$, o{\`u} pour $j\in \cF'$, on a pos{\'e} $y_j=x\sum_{i\in \cF;\ i\wedge k=j}x_i$. Comme $\varphi_\ell\Big(x\sum_{j\in \cF}x_j\Big)=0$, on trouve $\varphi_\ell\Big(\sum_{j\in \cF'} y_j\Big)=0$. Par hypoth{\`e}se de r{\'e}currence, on trouve $y_j=0$ pour tout $j\in {\cF'}$; en particulier, $y_k=0$. Or $y_k=\sum_{j\in \cF_k}xx_j=x\varphi_k(\sum_{j\in \cF_k}x_j)$. Ceci {\'e}tant vrai pour tout $x\in A_k$, le multiplicateur $\varphi_k(\sum_{j\in \cF_k}x_j)$ est nul. Par hypoth{\`e}se de r{\'e}currence l'application $\varphi_{k}:\gA_{\cF_k}\to M(A_{k})$ est injective, donc tous les $x_j$ sont nuls pour $j\in \cF_k$. Appliquant une derni{\`e}re fois l'hypoth{\`e}se de r{\'e}currence {\`a} $\cF'_k$, on trouve que tous les $x_j$ sont nuls.

Soit enfin $i\in \cL$. Pour tout sous-semi-treillis fini $\cF$ de
$\cL_i$, l'ensemble $\cF'=\cF\cup\{i\}$ est un
sous-semi-treillis (fini) de $\cL$ de plus petit {\'e}l{\'e}ment $i$, donc
la restriction de $\varphi_i $ {\`a} $\gA_{\cF'}$ est injective, donc la
restriction de $\varphi_i $ {\`a} $\gA_\cF\subset \gA_{\cF'}$ est
injective. Par la proposition~\ref{injectif},
$\varphi_i:\overline{\gA_{\cL_i}}\to M(A_i)$ est injective.
\end{proof}

Soit $(\gA,(A_\ell)_{\ell\in \cL})$ une $C^*$-alg{\`e}bre gradu{\'e}e. Rappelons qu'on
  obtient des morphismes $\pi_i:\gA\to M(A_i)$ par
  $\pi_i={\varphi}_i\circ p_i$ pour tout $i\in \cL$, o{\`u}
  $\varphi_i:\overline{\gA_{\cL_i}}\to M(A_i)$ et $p_i:\gA\to
  \overline{\gA_{\cL_i}}$. Suite {\`a} la proposition \ref{injectivite} on
  a le corollaire suivant:

\begin{corollary}\label{min}
Supposons que $\cL$ poss{\`e}de un plus petit {\'e}l{\'e}ment, $\ell_0$ et que les
morphismes de structure $\varphi_{k,\ell}$ avec $k\leq \ell$ satisfont
$\varphi_{k,\ell}^{-1}(A_k)=\{0\}$. Alors le morphisme
$\varphi_{\ell_0}=\pi_{\ell_0}$ est injectif.\newline\indent \hfill$\square$
\end{corollary}

\begin{remark}\label{atome} Soit $\cL$ un semi-treillis et
  $(\gA,(A_\ell)_{\ell\in \cL})$ une $C^*$-alg{\`e}bre gradu{\'e}e. Supposons
  que les morphismes de structure $\varphi_{k,\ell},k\leq \ell$
  satisfont \\$\varphi_{k,\ell}^{-1}(A_k)=\{0\}$. 

 La proposition \ref{injectivite} montre que $\ker\pi_i=\overline{\gA_{\cL'_i}}$ o{\`u}
  $\cL'_i=\{j\in \cL;i\not \leq j\}.$

Soit $I$ un ensemble (infini) et $(B_i)_{i\in I}$ une famille de $C^*$-alg{\`e}bres. Rappelons que le produit
$\ell^{\infty}$ des $B_i$ que nous noterons ici $\prod_{i\in
  I}^{\ell^{\infty}}B_i$ est $$\prod_{i\in
  I}^{\ell^{\infty}}B_i=\{(x_i)_{i\in I};\,x_i\in B_i\,\mbox{et}\,\sup_{i\in
  I}||x_i||\,\mbox{soit fini}\}$$ et la somme directe $c_o$ des $B_i$
que nous noterons ici $\bigoplus_{i\in
  I}^{c_o}B_i$ est $$\bigoplus_{i\in
  I}^{c_o}B_i=\{(x_i)_{i\in I};\,x_i\in B_i\,\mbox{et}\,||x_i||\to
0\;\mbox{quand}\;i\to \infty\}.$$ 

Soit $\cM$ une partie de $\cL$. Posons $B_\cM=M(\bigoplus\limits_{m\in \cM}^{c_o}A_m)=\prod_{m\in \cM}^{\ell^{\infty}}M(A_m)$. Notons $\pi={(\pi_m)}_{m\in \cM}:\gA\to B_\cM$ le morphisme donn{\'e} par
$a\mapsto \pi(a)={(\pi_m(a))}_{m\in \cM}$. On en deduit que\\
 $\begin{aligned}\ker \pi=\ker{(\pi_m)}_{m\in \cM}=\bigcap\limits_{m\in
  \cM}\ker\pi_m=\overline{\gA_J}\;\mbox{o{\`u}}\;J=\{j\in \cL:\forall m\in
  M,\,m\nleq j\}\end{aligned}.$\\
Donc le morphisme $\gA/\overline{\gA_J}\to B_\cM$ est injectif.\\

Supposons de plus que $\cL$ a un plus
  petit {\'e}l{\'e}ment, $\ell_0$. On rappelle qu'un \emph{atome} de $\cL$ est
  un {\'e}l{\'e}ment $a\neq \ell_0$ tel que si $b\leq a\Rightarrow
  b=\ell_0$ ou $b=a$. On dit que $\cL$ est \emph{atomique} si tout
  {\'e}l{\'e}ment $i\in \cL$ diff{\'e}rent de $\ell_0$
  est minor{\'e} par un atome. Notons encore $\cM$ l'ensemble des atomes de
  $\cL$. On remarque que l'ensemble $J$ poss{\`e}de un seul {\'e}l{\'e}ment
  $\ell_0$ et que le morphisme $\gA/A_{\ell_0}\to B_\cM$ est injectif. Notons $\widetilde {\gA}=\prod\limits_{m\in
  \cM}^{\ell^{\infty}}\overline{\gA_{\cL_m}}$. On en
  deduit que le morphisme
  $\gA/A_{\ell_0}\to \widetilde {\gA}$ est injectif.
\end{remark}
 
\newpage
\quad \thispagestyle{empty}

\newpage

\chapter{Reconstruction des $C^*$-alg{\`e}bres gradu{\'e}es}

\section{Reconstruction}
Nous d{\'e}montrons ici que les morphismes de structure nous permettent de
construire les
$C^*$-alg{\`e}bres gradu{\'e}es.

Soit $\cL$ un semi-treillis et ${(A_i)}_{i\in {\cal L}}$ une
 famille de $C^*$-alg{\`e}bres indix{\'e}es par $\cL$.

On veut construire une $C^*$-alg{\`e}bre gradu{\'e}e $\gA$ dont la composante
index{\'e}e par $i$ soit $A_i$. Si une telle $C^*$-alg{\`e}bre gradu{\'e}e est
donn{\'e}e, on dispose par la proposition \ref{q}
\begin{enumerate}
\item[$(i)$] des morphismes ${\varphi}_{i,j}:A_j\to M(A_i)$, pour tout
  $i,j\in \cal L$ tels que $i\leq j$, v{\'e}rifiant les conditions
\begin{enumerate}
\item[$a)$]${\varphi}_{i,i}={\mbox{Id}}_{A_i}$ pour tout $i\in {\cal L}$,

\item[$b)$] pour tout $i,j\in {\cal L}, m\leq i\wedge j$ et tout
  $x\in A_i,y\in A_j$ posons $k=i\wedge j$ alors,
$${\varphi}_{k,i}(x){\varphi}_{k,j}(y)\in A_k$$ et 
 $${\varphi}_{m,k}({\varphi}_{k,i}(x){\varphi}_{k,j}(y))={\varphi}_{m,i}(x){\varphi}_{m,j}(y).$$

\end{enumerate}
\item[$(ii)$] des applications bilin{\'e}aires $q_{i,j}:A_i\times A_j\to
  A_{i\wedge j}$, pour tout $i,j\in \cal L$, v{\'e}rifiant les conditions
  
\begin{enumerate}
\item[$a')$] $q_{i,i}(x,y)=xy$ pour tout $i\in {\cal L}$,

\item[$b')$] pour tout $i,j\in {\cal L}$ et $x\in A_i,y\in
A_j$,  $q_{i,j}(x,y)={q_{j,i}(y^*,x^*)}^*$ et

\item[$c')$] pour tout $i,j,k\in {\cal L}$ et $x\in A_i,y\in A_j,z\in A_k$
 on a $$q_{i\wedge j,k}(q_{i,j}(x,y),z)=q_{i,j\wedge
  k}(x,q_{j,k}(y,z)).$$
\end{enumerate}

  \end{enumerate}

\begin{proposition}\label{equiv} Soit $\cL$ un semi-treillis et
  ${(A_i)}_{i\in {\cal L}}$ une famille de $C^*$-alg{\`e}bres. 

\begin{enumerate}
\item[$(1)$] On se donne pour tout $i,j\in \cal L$ tels que $i\leq j$,
  des morphismes ${\varphi}_{i,j}:A_j\to M(A_i)$ qui v{\'e}rifient les
  conditions $a)$ et $b)$. Pour $i,j\in
  \cL$ notons $q_{i,j}:A_i\times A_j\to A_{i\wedge j}$ l'application
  $(x,y)\mapsto \varphi_{i\wedge j,i}(x)\varphi_{i\wedge j,j}(y)$. La
  famille d'applications $q_{i,j}$ v{\'e}rifie les conditions $a')$,
  $b')$, $c')$ ci-dessus.

\item[$(2)$] Inversement, on se donne pour tout $i,j\in \cal L$ des
  applications $q_{i,j}:A_i\times A_j\to A_{i\wedge j}$ qui v{\'e}rifient
  les conditions $a')$, $b')$ et $c')$ ci-dessus. Pour
 $i,j\in \cL$ tels que $i\leq j$ il existe un morphisme $\varphi_{i,j}:A_j\to M(A_i)$
  qui est donn{\'e} par $\varphi_{i,j}(y)x=q_{j,i}(y,x)$  et $x\varphi_{i,j}(y)=q_{i,j}(x,y)$ pour $x\in A_i,\,y\in
  A_j$. La famille de morphismes $\varphi_{i,j}$ v{\'e}rifie les
  conditions $a),\,b)$ ci-dessus.   
  \end{enumerate}
\end{proposition}

\begin{proof}
\begin{enumerate}
\item[$(1)$] Soient $i,j\in \cL$ et $x\in A_i,\,y\in A_j$.\\ La
  condition $a')$ est satisfaite car par $a)$ on a  $q_{i,i}(x,y)=\varphi_{i,i}(x)\varphi_{i,i}(y)=xy$.

On a $b')$ puisque

 $\begin{aligned}q_{j,i}(y^*,x^*)&=\varphi_{i\wedge j,j}(y^*)\varphi_{i\wedge
  j,i}(x^*)={\varphi_{i\wedge j,j}(y)}^*{\varphi_{i\wedge
  j,i}(x)}^*\\&={(\varphi_{i\wedge j,i}(x)\varphi_{i\wedge
  j,j}(y))}^*={q_{i,j}(x,y)}^*.\end{aligned}$

Soit aussi $z\in A_k$, en utilisant $b)$ on a \\
 $\begin{aligned}q_{i,j\wedge k}(x,q_{j,k}(y,z))&=\varphi_{i\wedge
  j\wedge k,i}(x)\varphi_{i\wedge j\wedge k,j\wedge k}(q_{j,k}(y,z))\\&=\varphi_{i\wedge
  j\wedge k,i}(x)\varphi_{i\wedge j\wedge k,j\wedge
  k}(\varphi_{j\wedge k,j}(y)\varphi_{j\wedge
  k,k}(z))\\&=\varphi_{i\wedge j\wedge k,i}(x)\varphi_{i\wedge j\wedge
  k,j}(y)\varphi_{i\wedge j\wedge k,k}(z).\end{aligned}$

De m{\^e}me on a $q_{i\wedge j,k}(q_{i,j}(x,y),z)=\varphi_{i\wedge j\wedge k,i}(x)\varphi_{i\wedge j\wedge
  k,j}(y)\varphi_{i\wedge j\wedge k,k}(z)$ donc la famille $q_{i,j}$
v{\'e}rifie aussi $c')$.

\item[$(2)$]
  Soient $i,j\in \cL$ tels que $i\leq j$ et $x\in A_i,\,y\in
  A_j$. Le morphisme $\varphi_{i,j}(y)$ est un multiplicateur puisque si $x'\in A_i$ alors
on a

 $\begin{aligned}x'(\varphi_{i,j}(y)x)&=x'q_{j,i}(y,x)=q_{i,i}(x',q_{j,i}(y,x))=q_{i,i}(q_{i,j}(x',y),x)\\&=q_{i,j}(x',y)x=(x'\varphi_{i,j}(y))x.\end{aligned}$

 En utilisant les propri{\'e}t{\'e}s
  $b')$ et $c')$ des $q_{i,j}$ c'est facile de montrer que
  $\varphi_{i,j}$ est
  un morphisme de $C^*$-alg{\`e}bres.\\

C'est {\'e}vident qu'on a $a)$.

 Soient $i,j\in \cL$ et posons $k=i\wedge j$ alors $j\wedge
  k=k$. Soient $x\in A_i,\,y\in A_j,\,z\in A_k$, en utilisant $c')$ on
  va d{\'e}montrer $b)$.\\
$\begin{aligned}\varphi_{i\wedge j,i}(x)\varphi_{i\wedge
  j,j}(y)(z)&=\varphi_{i\wedge j,i}(x)(q_{j,i\wedge
  j}(y,z))=q_{i,i\wedge j}(x,q_{j,i\wedge j}(y,z))\\&=q_{i,j\wedge
  k}(x,q_{j,k}(y,z))=q_{i\wedge
  j,k}(q_{i,j}(x,y),z)\\&=q_{k,k}(q_{i,j}(x,y),z)=q_{i,j}(x,y)z.\end{aligned}$

 Donc $\varphi_{i\wedge j,i}(x)\varphi_{i\wedge
  j,j}(y)=q_{i,j}(x,y)\in A_k$.

Soit $m\leq i\wedge j$ et $z\in A_m$, alors\\
$\begin{aligned}\varphi_{m,i\wedge j}(\varphi_{i\wedge
  j,i}(x)\varphi_{i\wedge j,j}(y))z&=q_{i\wedge
  j,m}(q_{i,j}(x,y),z)=q_{i,j\wedge
  m}(x,q_{j,m}(y,z))\\&=q_{i,m}(x,q_{j,m}(y,z))=\varphi_{m,i}(x)q_{j,m}(y,z)\\&=\varphi_{m,i}(x)\varphi_{m,j}(y)z.\end{aligned}$\\
Donc $\varphi_{m,i\wedge j}(\varphi_{i\wedge j,i}(x)\varphi_{i\wedge
  j,j}(y))=\varphi_{m,i}(x)\varphi_{m,j}(y)$.  
  \end{enumerate}
\end{proof}

\begin{theorem} \label{carctgrad}
 Soit $\cal L$ un semi-treillis et ${(A_i)}_{i\in {\cal L}}$ une
 famille de $C^*$-alg{\`e}bres. On se donne des morphismes
  $${\varphi}_{i,j}:A_j\to
 M(A_i)\;\mbox{pour}\;i\leq j,i,j\in {\cal L},$$ qui verifient les 
propri{\'e}tes $a)$ et $b)$ de la proposition \ref{q}. Alors,
il existe {\`a} isomorphisme pr{\`e}s une unique $C^*$-alg{\`e}bre $\cal L$-gradu{\'e}e
  $(\frak A,{(A_i)}_{i\in \cal L})$ dont les morphismes de structure
  soient les $\varphi_{i,j}$.
\end{theorem}

\begin{proof} On va faire la d{\'e}monstration par {\'e}tapes:

$I)$ $\underline {\mbox{existence}}$

 On d{\'e}finit ${\frak A}_{\cal L}=\bigoplus\limits_{i\in {\cal L}}A_i$ et on
id{\'e}ntifie $A_i$ avec son image canonique dans ${\frak A}_{\cal L}$. Un
{\'e}l{\'e}ment $x\in{\frak A}_{\cal L}$ est {\'e}crit sous la forme $x=\sum\limits_{i\in {\cal L}}x_i$ avec $x_i\neq 0$ pour seulement un
nombre fini de $i\in {\cal L}$. On va d{\'e}montrer que le compl{\'e}t{\'e} de
$\frak A_{\cal L}$  pour une norme qu'on d{\'e}finira est une
$C^*$-alg{\`e}bre $\cal L$-gradu{\'e}e. On notera $\frak A$ ce compl{\'e}t{\'e}.

 Soient $x,y\in {\frak A}_{\cal L}$ alors $x=\sum\limits_{i\in {\cal L}}x_i={(x_i)}_{i\in {\cal L}}\;\mbox{et}\;y=\sum\limits_{j\in
 {\cal L}}y_j={(y_j)}_{j\in {\cal L}}$ avec $x_i\in A_i,y_j\in
 A_j$ pour tout $i,j\in {\cal L}$.
\begin{itemize}

\item Pour tout $i,j\in {\cal L}$ on pose $k=i\wedge j$ et on d{\'e}finit le
 produit $:$\\
$$x_i y_j={\varphi}_{k,i}(x_i){\varphi}_{k,j}(y_j)\in A_k.$$ 
Soit $x$ et $y\in \gA_\cL$. On d{\'e}finit le produit dans ${\frak A}_{\cal L}$ par $:$ $$x y=\sum\limits_{i,j\in {\cal L
}}x_i y_j.$$
Il r{\'e}sulte imm{\'e}diatement de la proposition \ref{equiv} que ce produit
 est associatif.

D{\'e}finissons une application $*$ de $\gA_\cL$ dans lui-m{\^e}me par $x\mapsto x^*=\sum\limits_{i\in {\cal L
}}x_i^*$ pour tout $x\in \gA_\cL$. C'est clair qu'on a $x^{**}=x$ et
$(xy)^*=y^*x^*$ pour tout $x,y\in \gA_\cL$, donc cette application est
une involution. 

L'espace $\gA_\cL$ muni de ce
produit et de cette involution est une alg{\`e}bre
 involutive.\\

 D{\'e}finissons la norme de $\gA_\cL$.

\item
Pour tout $i\in \cL$, on d{\'e}finit une application lin{\'e}aire $${\pi}_i:{\frak A}_{\cal L}\to M(A_i)$$ par
 $${\pi}_i(x)=\pi_i(\sum\limits_{j\in {\cal
 L}}x_j)=\sum\limits_{j\in \cL;i\leq j}\varphi_{i,j}(x_j)\,\mbox{pour}\,x=\sum\limits_{j\in {\cal L}}x_j.$$ 
Montrons que ${\pi}_i$ est un morphisme d'alg{\`e}bres involutives.\\

Par bilin{\'e}arit{\'e} il suffit de montrer que pour tout $i,j,k\in {\cal L}$
 et pour tout $x_k\in A_k,y_j\in A_j$ on a ${\pi}_i(x_k y_j)={\pi}_i(x_k){\pi}_i(y_j).$\\

Si $i\in \cL$ est tel que $i\not \leq k\wedge j$ alors $i\not \leq j$ ou
 $i\not \leq k$, donc 

${\pi}_i(x_k){\pi}_i(y_j)=0=\pi_i(x_k y_j)$.\\

Si $i\leq i\wedge j$ alors 

$\begin{aligned}\pi_i(x_k y_j)&=\varphi_{i,k\wedge
  j}(x_k y_j)=\varphi_{i,k\wedge j}(\varphi_{k\wedge
  j,k}(x_k)\varphi_{k\wedge
  j,j}(y_j))=\varphi_{i,k}(x_k)\varphi_{i,j}(y_j)\\&=\pi_i(x_k)\pi_i(y_j).\end{aligned}$\\

De plus on a pour tout $x\in {\frak A}_{\cal L}$ et tout $i\in \cal L$, $${\pi}_i(x^*)=\sum_{\substack{j\in {\cal L}:\\{i\leq
 j}}}{\varphi}_{i,j}(x^*_j)=\sum_{\substack{j\in {\cal L}:\\{i\leq
 j}}}{{\varphi}_{i,j}(x_j)}^*={{\pi}_i(x)}^*.$$
 
\item Soit $\cal S\subset \cL$ une partie quelconque, on d{\'e}finit
 ${||x||}_{\cal S}=\sup\limits_{i\in {\cal
 S}}||{\pi}_i(x)||.$\\
Remarquons tout d'abord que pour tout $\cS$, ${||\cdot||}_{\cal S}$ est finie:\\
$$\begin{aligned}{||x||}_{\cal S}=&\sup\limits_{i\in {\cal
 S}}||{\pi}_i(x)||=\sup\limits_{i\in {\cal S}}||\sum\limits_{j:i\leq
 j}{\varphi}_{i,j}(x_j)||\\&\leq  {\sup\limits_{i\in {\cal S}}\sum\limits_{j:i\leq
 j}||{\varphi}_{i,j}(x_j)||}\leq {\sup\limits_{i\in {\cal S}}\sum\limits_{j:i\leq
 j}||x_j||}\leq {\sum\limits_{j\in {\cal
 S}}||x_j||}.\end{aligned}$$

De plus, ${||\cdot||}_{\cal S}$ est une $C^*$-semi-norme puisque c'est
 le supremum d'une famille de $C^*$-semi-normes.
\end{itemize}

\begin{lemma} Si $m$ est maximal dans le support de $x\in {\frak A}_{\cal L}$ alors
 $x_m={\pi}_m(x)$.\end{lemma}
 
\begin{proof} Le r{\'e}sultat est immediat par d{\'e}finition \ie 

${\pi}_m(x)=\sum\limits_{j\in {\cal
 L}}{\pi}_m(x_j)=\sum\limits_{j\in {\cal L}:m\leq
 j}{\varphi}_{m,j}(x_j)={\varphi}_{m,m}(x_m)=x_m.$ \end{proof}

\begin{proposition} $||\cdot||_\cL$ est une $C^*$-norme sur $\gA_\cL$.
\end{proposition}

\begin{proof}
 Soit $x\in {\frak A}_{\cal L}$ non nul, alors
 $x=\sum\limits_{i\in \cal L} x_i$ avec $x_i\neq 0$ pour seulement un
 nombre fini de $i\in \cal L$. Notons $\Omega=\{j\in {\cal L}:x_j\neq 0\}$
 le support de $x$. Soit $m\in \Omega$ un {\'e}l{\'e}ment maximal alors $x_m\neq 0$
 par d{\'e}finition de $\Omega$. Par le lemme pr{\'e}cedent on a ${\pi}_m(x)\neq 0\Rightarrow ||{\pi}_m(x)||\neq
 0\Rightarrow \sup\limits_{k\in \Omega}||{\pi}_k(x)||\neq 0\Rightarrow
 {||x||}_\cL\neq 0.$
\end{proof}

Le complet{\'e} de  ${\frak A}_{\cal L}$ par rapport {\`a} la norme
${||\cdot||}_{\cal L}$ est une $C^*$-alg{\`e}bre.\\

 Montrons que $\frak
  A$  est $\cal L$-gradu{\'e}e. 
Les $A_i$ sont lin{\'e}airement independents car ils le sont
dans $\gA_\cL$ (par d{\'e}finition) et ${||\cdot||}_\cL$ {\'e}tant une norme sur
$\gA_\cL$, l'application $\gA_\cL\to \gA$ est injective. 

De plus par la d{\'e}finition du produit on a $A_iA_j\subset A_{i\wedge
  j}$ pour tout $i,j\in \cal L$ et finalement par construction on a que $\bigcup\limits_{{\cal F}\in \F_\cL}{\frak A}_{\cal
  F}\equiv \bigoplus\limits_{i\in \cal L}A_i={\frak A}_{\cal
  L}$ est dense dans son compl{\'e}t{\'e} $\frak A$ pour la norme $||\cdot||_\cL.$

$II)$ $\underline {\mbox{unicit{\'e}}}$

 Soit ${\frak B}=\overline {\frak A}$ pour une autre norme ${||\cdot||}_*$. 
 Pour $i\in \cL$ notons ${\rho}_i:A_i\to B_i=A_i$
  l'id{\'e}ntit{\'e} ${\mbox{Id}}_{A_i}$. C'est {\'e}vident qu'on a ${\rho}_i(x){\rho}_j(y)={\rho}_{i\wedge
  j}(xy)$ pour tout $x\in A_i,y\in A_j$ et $i,j\in \cal L$. Donc par le
  corollaire \ref{morgrad} il existe un unique morphisme $\rho:{\frak A}\to
  {\frak B}$ qui est injectif donc isom{\'e}trique, d'o{\`u} l'unicit{\'e}.

\end{proof}

\begin{corollary}\label{kerim}Soient $(\frak A,{(A_i)}_{i\in \cal L})$ et $(\frak B,{(B_i)}_{i\in
  \cal L})$ des $C^*$-alg{\`e}bres \\$\cL$-gradu{\'e}es. On se donne un morphisme $\rho:{\frak
    A}\to {\frak B}$ de $C^*$-alg{\`e}bres gradu{\'e}es. Notons $\rho_i$ la
    restriction de $\rho$ {\`a} $A_i$. Alors
\begin{enumerate} 
\item $\ker \rho=\overline{\bigoplus\limits_{i\in \cal L}\ker {\rho}_i}$
et 
\item $\mbox{Im}\rho=\overline{\bigoplus\limits_{i\in \cal
      L}\mbox{Im}{\rho}_i}$. 
\end{enumerate}
\end{corollary}

\begin{proof}
\begin{enumerate}
\item 
 Pour tout $i\in \cL$ on a $\ker \rho_i\subset \ker \rho$, donc
$\bigoplus\limits_{i\in \cL}\ker \rho_i\subset \ker \rho$.

 Soit $\cal F$ un sous-semi-treillis fini de $\cal L$. Montrons
 $\bigoplus\limits_{i\in \cal F}\ker {\rho}_i=\ker {\rho}_{\cal F}.$ 

Soit $x=\sum\limits_{i\in \cF}x_i\in \ker
{\rho}_{\cal F}$ alors ${\rho}_{\cal F}(x)=0$ donc $\sum\limits_{i\in \cal
  F}{\rho}_i(x_i)=0$. Puisque ${\rho}_i(x_i)\in B_i$ pour tout $i\in \cal F$
et ${(B_i)}_{i\in \cal F}$ est une famille lin{\'e}airement ind{\'e}pendente
on a ${\rho}_i(x_i)=0$ pour tout $i\in {\cal F}$. Donc $x_i\in \ker
{\rho}_i$ pour tout $i\in {\cal F}$ et $x\in \bigoplus\limits_{i\in \cal
  F}\ker {\rho}_i.$

 Si $x\in \ker {\rho}$ alors pour tout $\varepsilon>0$ il existe $\cal F$ un sous-semi-treillis fini de $\cal L$ et
il existe $y\in {\frak A}_{\cal F}$ tel que $||x-y||<{{\varepsilon}\over
  2}$. On a $||\rho(y)||=||\rho(y-x)||\leq {\varepsilon\over 2}.$
Puisque $\rho(y)\in \rho({\frak A}_{\cal F})$ par la proposition \ref{isom}
il existe $z\in {\frak A}_{\cal F}$ tel que $\rho(z)=\rho(y)$ et
$||z||=||\rho(y)||<{\varepsilon\over 2}$. \\Mais $\rho(y-z)=0\Rightarrow y-z\in
\ker \rho\cap{\frak A}_{\cal F}\Rightarrow y-z\in \ker {\rho}_{\cal
  F}=\bigoplus\limits_{i\in \cal F}\ker {\rho}_i$. Alors
$||x-(y-z)||\leq ||x-y||+||z||<\varepsilon$ \ie on a trouv{\'e} un
{\'e}l{\'e}ment $y-z$ dans $\ker {\rho}_{\cal F}$ qui est proche de $x\in \ker
\rho$, donc $\ker \rho=\overline{\bigoplus\limits_{i\in \cal L}\ker
  {\rho}_i}$.

\item 
 Comme ${(A_i)}_{i\in \cL}$ est totale dans $\gA$,
${(\rho(A_i))}_{i\in \cL}$ est totale dans $\rho(\gA)$. 
\end{enumerate} 
\end{proof} 

\section{Id{\'e}aux et quotients de $C^*$-alg{\`e}bres gradu{\'e}es}
De la proposition \ref{injectif} et du corollaire \ref{kerim} on
deduit immediatement:
\begin{corollary}\label{exactgr} Soient $(\mathfrak J,(J_\ell)_{\ell\in \cL})$,
  $(\gA,(A_\ell)_{\ell\in \cL})$ et $(\gB,{(B_\ell)}_{\ell\in
  \cL})$ des $C^*$-alg{\`e}bres $\cL$-gradu{\'e}es. On se donne $i:\mathfrak J\to
  \gA$ et $p:\gA\to \gB$ des morphismes de $C^*$-alg{\`e}bres gradu{\'e}es. Si
  pour tout $\ell\in \cL$ la suite $0\to J_\ell\to A_\ell\to B_\ell\to
  0$ est exacte, alors la suite $0\to \mathfrak J\to \gA\to \gB\to 0$
  est exacte.\newline\indent \hfill$\square$
\end{corollary}

\begin{proposition} Soit $(\gA,(A_i)_{i\in \cL})$ une $C^*$-alg{\`e}bre
  gradu{\'e}e et prenons la surjection canonique $\pi:\gA\to \gA/I$ o{\`u} $I$
  est un id{\'e}al ferm{\'e} de $\gA$. On pose $\pi(A_i)=B_i$. Alors
  on a une {\'e}quivalence entre:
\begin{enumerate}
\item $I$ est une sous-$C^*$-alg{\`e}bre $\cL$-gradu{\'e}e de $\gA$ (et dans ce cas
  on a \\$\overline {\bigoplus\limits_{i\in \cL}I\cap A_i}=I$ d'apr{\`e}s la
  d{\'e}finition \ref{sg}). 

\item La famille $(B_i)$ est lin{\'e}airement ind{\'e}pendente, autrement dit\\
  $(\gA/I,\;{(B_i)}_{i\in \cL})$ est une $C^*$-alg{\`e}bre $\cL$-gradu{\'e}e.

\end{enumerate}

\end{proposition}

\begin{proof}$b)\Rightarrow a)$ On applique le corollaire \ref{kerim} pour $\ker\pi=I$.\\

$a)\Rightarrow b)$. Comme $\gA$ est une
$C^*$-alg{\`e}bre $\cL$ gradu{\'e}e, par la proposition \ref{q} il existe une
unique application bilin{\'e}aire $q_{i,j}:A_i\times A_j\to A_{i\wedge j}$
d{\'e}finie par $q_{i,j}(x,y)=xy$ pour tout $x\in A_i$ et $y\in
A_j$ qui v{\'e}rifie les conditions $a'),\,b'),\,c')$ de cette proposition. On a le diagramme suivant:
$$\xymatrix{
           A_i\times A_j \ar[r]^{q_{i,j}} \ar[d]_{\pi_{i,j}} &
            A_{i\wedge j} \ar[d]^{\pi_{i\wedge j}} \\
           B_i\times B_j\ar[r]_{q'_{i,j}} & B_{i\wedge
           j}
         }$$ o{\`u} $q'_{i,j}$ est une application bien d{\'e}finie. En effet,
soit $x_i=\pi(x'_i)\in B_i$ et $y_j=\pi(y'_j)\in B_j$ (avec $x'_i\in
           A_i$ et $y'_j\in A_j$) alors $x_i y_j=\pi(x'_i y'_j)\in
           B_{i\wedge j}$ et $q'_{i,j}(x_i,y_j)=x_i y_j$ est bien
           d{\'e}finie.

 Comme
           les applications $\pi_{i,j}$ et $\pi_{i\wedge j}$ sont
           surjectives $q'_{i,j}$ v{\'e}rifie aussi les conditions
           $a'),\,b'),\,c')$ de la proposition \ref{q}. 

Alors d'apr{\`e}s le th{\'e}or{\`e}me \ref{carctgrad} il
           existe {\`a} isomorphisme pr{\`e}s une unique $C^*$-alg{\`e}bre
           $\cL$-gradu{\'e}e $(\gB,\,{(B_i)}_{i\in \cL})$ dont les
           applications de structure soient les $q'_{i,j}$.

Pour tout $i\in \cL$, notons $\pi_i$ la restriction de $\pi$ {\`a}
$A_i$. D'apr{\`e}s le corollaire \ref{morgrad} il existe un unique
morphisme $\pi':\gA\to \gB$ dont la restriction {\`a} $A_i$ soit $\pi_i$
pour tout $i\in \cL$. Comme $\ker \pi_i=I\cap A_i$ pour tout $i\in
\cL$ et par hypoth{\`e}se, on a $\ker
\pi'=\overline{\bigoplus\limits_{i\in \cL}\ker
  \pi_i}=\overline{\bigoplus\limits_{i\in \cL}I\cap A_i}=I=\ker
\pi$. Les applications $\pi$ et $\pi'$ sont surjectives et ont m{\^e}me
noyau. Il existe donc un unique isomorphisme $\rho:\gB\to \gA/I$ tel
que $\rho\circ  \pi'=\pi$. 

Comme  $(\gB,\,{(B_i)}_{i\in \cL})$ est une $C^*$-alg{\`e}bre gradu{\'e}e
$(\gA/I,(\rho(B_i))_{i\in \cL})$ est gradu{\'e}e. Or, en consid{\'e}rant la
restriction de $\rho\circ \pi'$ et $\pi$ {\`a} $A_i$, on trouve
$\rho_{|_{B_i}}=\mbox{Id}_{B_i}$ de sorte que $\rho(B_i)=B_i$, d'o{\`u} le
r{\'e}sultat. 
\end{proof}

\newpage
\quad \thispagestyle{empty}

\newpage
\chapter{Produits tensoriels et produits crois{\'e}s de $C^*$-alg{\`e}bres gradu{\'e}es}

\section{Produits tensoriels}
Les produits tensoriels de deux (un nombre fini) $C^*$-alg{\`e}bres
gradu{\'e}es par des semi-treillis finis ont {\'e}t{\'e} {\'e}tudi{\'e}s dans
\cite{gi4}. Nous g{\'e}n{\'e}ralisons ici leur r{\'e}sultat dans le cas de
semi-treillis quelconques.

Notons $\alpha$ la $C^*$-norme minimale \emph{min} ou maximale \emph{max}.
\begin{lemma}\label{ten}Soit $(\gA,(A_\ell)_{\ell\in \cL})$ une
    $C^*$-alg{\`e}bre gradu{\'e}e et $C$ une $C^*$-alg{\`e}bre. Alors $(\gA\otimes_{\alpha}
    C,(A_\ell\otimes_{\alpha} C)_{\ell\in \cL})$ est aussi une $C^*$-alg{\`e}bre
    $\cL$-gradu{\'e}e. 
\end{lemma}

\begin{proof}Tout d'abord, montrons que pour tout $\ell\in \cL$
    l'inclusion naturelle $\varphi_\ell:A_\ell\hookrightarrow \gA$ induit un morphisme
    injectif $\varphi_\ell\otimes_{\alpha} \mbox{Id}_C:A_\ell\otimes_{\alpha} C\to
    \gA\otimes_{\alpha} C$. En d'autres termes, montrons que $(A_\ell\otimes_{\alpha} C)_{\ell\in \cL}$ est une
    famille de sous-$C^*$-alg{\`e}bres de $\gA\otimes_{\alpha} C$.

 C'est imm{\'e}diat pour le produit tensoriel
    minimal.

 Montrons-le pour le produit tensoriel maximal. Soit $\ell\in \cL$,
 puisque $\gA$ est une $C^*$-alg{\`e}bre gradu{\'e}e on a par la
proposition \ref{suitexacte} une suite exacte scind{\'e}e
$$0\to
 \overline
 {\gA_{{\cL}'_\ell}}\xrightarrow{} \gA\
 \mathop{\rightleftarrows}\limits^{p_\ell}_{{\sigma}_\ell}\ \overline
 {\gA_{\cL_\ell}}\to 0.$$
Par la proposition \ref{mormax} on a des morphismes 

$$\xymatrix{
\gA\otimes_{\max} C \;\;\; \ar@<2pt>[r]^{p_\ell\otimes_{\max}
 \mbox{Id}_C} & \;\;\; \overline
 {\gA_{\cL_\ell}}\otimes_{\max} C  \ar@<2pt>[l]^{{\sigma}_\ell\otimes_{\max}
   \mbox{Id}_C}
}$$

et $(p_\ell\otimes_{\max} \mbox{Id}_C)\circ
 ({\sigma}_\ell\otimes_{\max} \mbox{Id}_C)=(p_\ell\circ
 \sigma_\ell)\otimes \mbox{Id}_C=\mbox{Id}_{\overline
 {\gA_{\cL_\ell}}\otimes_{\max} C}$. Donc le morphisme
 ${\sigma}_\ell\otimes_{\max} \mbox{Id}_C$ est injectif. Puisque
 $A_\ell\otimes_{\max}C$ est un id{\'e}al de $\overline{\gA_{\cL_\ell}}\otimes_{\max}C$ on a le
 diagramme suivant: 

  $$\xymatrix{
    A_\ell\otimes_{\max}C \ar@{^{(}->}[r]^{}  \ar[rd]_{\varphi_\ell\otimes_{\max}\mbox{Id}_C} & \overline {\gA_{\cL_\ell}}\otimes_{\max}C  \ar[d]^{{\sigma}_{\ell}\otimes_{\max}\mbox{Id}_C} \\
    {} & \gA\otimes_{\max}C  
  }$$ 
et par cons{\'e}quent $\varphi_\ell\otimes_{\max}\mbox{Id}_C$ est injectif.\\

Montrons ensuite que $(A_\ell\otimes_{\alpha} C)_{\ell\in \cL}$ est une
famille lin{\'e}airement ind{\'e}pendente. Soit $(x_i)_{i\in \cL}$ une famille
telle que $x_i\in A_i\otimes_{\alpha} C$ pour tout $i\in \cL$ et notons $I=\{j\in {\cal L}:x_j\neq 0\}$
 le support de $x=\sum\limits_{i\in \cL}x_i$. On va montrer par
 recurrence sur le nombre d'{\'e}l{\'e}ments de $I$ que si $\sum\limits_{i\in
   I}x_i=0$ alors $x_i=0$ pour tout $i\in I$. 

C'est evident si $I$ a un seul {\'e}l{\'e}ment. Supposons que c'est vrai pour
tout sous ensemble de $I$ distinct de $I$. Soit $m\in I$ un {\'e}l{\'e}ment
maximal et $x=\sum\limits_{i\in I}x_i=0$. Puisque $\gA$ est gradu{\'e}e on a 
$$0\to
 \overline
 {\gA_{{\cL}'_m}}\xrightarrow{} \gA\
 \mathop{\rightleftarrows}\limits^{p_m}_{{\sigma}_m}\ \overline
 {\gA_{\cL_m}}\to 0.$$
Par les propositions \ref{exmax} et \ref{exspa} on a aussi:

$$\xymatrix{
0 \ar[r]^{} &  \overline
 {\gA_{{\cL}'_m}}\otimes_{\alpha} C \ar[r]^{} & \gA\otimes_{\alpha} C\;\;
 \ar@<2pt>[r]^{p_m\otimes_{\alpha}
 \mbox{Id}_C} & \;\;\; \overline
 {\gA_{\cL_m}}\otimes_{\alpha} C \ar@<2pt>[l] ^{{\sigma}_m\otimes_{\alpha}
   \mbox{Id}_C} \ar[r] & 0.
}$$
Donc $(p_m\otimes_{\alpha} \mbox{Id}_C)(x)=x_m=0$ et $\sum\limits_{i\in
 I\backslash \{m\}}x_i=0$. Comme $I\backslash \{m\}$ est un
 sous-ensemble de $I$ par hypoth{\`e}se de recurrence on a $x_i=0$ pour
 tout $i\in I\backslash \{m\}$, d'o{\`u} le r{\'e}sultat.

La famille $(A_\ell\otimes_{\alpha}
C)_{\ell\in \cL}$ est clairement totale.

Soient $\ell$, $\ell'\in \cL$, on v{\'e}rifie facilement avec les tenseurs
{\'e}l{\'e}mentaires qu'on a
$(A_\ell\otimes_{\alpha} C)(A_\ell'\otimes_\alpha C)\subset A_{\ell\wedge
  \ell'}\otimes_\alpha C$.

\end{proof}

On va g{\'e}n{\'e}raliser ce r{\'e}sultat par la proposition suivante

\begin{proposition}\label{tens}Soient
  $\cL$ et $\cM$ deux semi-treillis. Soit $(\gA,(A_\ell)_{\ell\in \cL})$
  une $C^*$-alg{\`e}bre $\cL$-gradu{\'e}e et
  $(\gB,(B_m)_{m\in \cM})$ une $C^*$-alg{\`e}bre $\cM$-gradu{\'e}e, alors
  $(\gA\otimes_\alpha \gB,(A_\ell\otimes_\alpha B_m)_{(\ell,m)\in
    \cL\times \cM})$ est une $C^*$-alg{\`e}bre $\cL\times \cM$-gradu{\'e}e.
\end{proposition}

\begin{proof} Soit $\ell\in \cL,m\in \cM$, puisque $\gA$ et $\gB$ sont des $C^*$-alg{\`e}bres gradu{\'e}es on a par la
proposition \ref{suitexacte} des suites exactes scind{\'e}es
$$0\to
 \overline
 {\gA_{{\cL}'_\ell}}\xrightarrow{} \gA\
 \mathop{\rightleftarrows}\limits^{p_\ell}_{{\sigma}_\ell}\ \overline
 {\gA_{\cL_\ell}}\to 0,$$

$$0\to
 \overline
 {\gB_{{\cM}'_m}}\xrightarrow{} \gB\
 \mathop{\rightleftarrows}\limits^{p'_m}_{{\sigma}'_m}\ \overline
 {\gB_{\cM_m}}\to 0.$$

On remarque que si $\varphi_\ell:A_\ell\hookrightarrow
\gA$ et $\psi_m:B_m\hookrightarrow \gB$ sont les inclusions naturelles
pour tout $\ell\in \cL$ et $m\in \cM$, alors par les propositions
\ref{morspatial} et \ref{mormax} on
peut d{\'e}finir un morphisme
$\varphi_\ell\otimes_\alpha \psi_m:A_\ell\otimes_\alpha B_m\to
\gA\otimes_\alpha \gB$ pour tout $\ell\in \cL$ et $m\in \cM$. Montrons que
$\varphi_\ell\otimes_\alpha \psi_m$ est injectif.\\ 

Par la
proposition \ref{morspatial} le morphisme $\varphi_\ell\otimes_{\min}\psi_m$ est
injectif donc $A_\ell\otimes_{\min}B_m$ est une sous-$C^*$-alg{\`e}bre
de $\gA\otimes_{\min}\gB$ pour tout $\ell\in \cL$ et $m\in \cM$.

Montrons-le pour le produit crois{\'e} maximal. Soit $\ell\in \cL$ et
$m\in \cM$. Puisque $\gB$ est une $C^*$-alg{\`e}bre gradu{\'e}e en appliquant
le lemme \ref{ten} pour $C=A_\ell$ on a un morphisme injectif
$$\mbox{Id}_{A_\ell}\otimes_{\max} \psi_m:A_\ell\otimes_{\max} B_m\to
A_\ell\otimes_{\max} \gB.$$
 Maintenant on utilise la graduation de
$\gA$ et on applique une deuxi{\`e}me fois le lemme \ref{ten} avec
$C=\gB$. On obtient donc un morphisme injectif
$$\varphi_\ell\otimes_{\max} \mbox{Id}_\gB:A_\ell\otimes_{\max} \gB\to
\gA\otimes_{\max} \gB.$$ La composition de ces deux morphismes
injectifs $(\varphi_\ell\otimes_{\max} \mbox{Id}_\gB)\circ
(\mbox{Id}_{A_\ell}\otimes_{\max}
\psi_m)=\varphi_\ell\otimes_{\max}\psi_m$ est aussi un morphisme injectif.   
Autrement dit $A_\ell\otimes_{\max}B_m$ est une sous-$C^*$-alg{\`e}bre
de $\gA\otimes_{\max}\gB$ pour tout $\ell\in \cL$ et $m\in \cM$.\\

Montrons ensuite que la famille $(A_\ell\otimes_\alpha B_m)_{(\ell,m)\in \cL\times
  \cM}$ est lin{\'e}airement ind{\'e}pendente. Il est clair (sur les tenseurs
el{\'e}mentaires) que 
  $A_\ell\otimes_\alpha B_m\subset (A_\ell\otimes_\alpha \gB)\cap (\gA\otimes_\alpha B_m)$
  pour tout $\ell\in \cL$ et $m\in \cM$. Mais par le lemme \ref{ten} les familles
  $(A_\ell\otimes_\alpha \gB)_{\ell\in \cL}$ et $(\gA\otimes_\alpha B_m)_{m\in \cM}$
  sont lin{\'e}airement ind{\'e}pendentes, d'o{\`u} le r{\'e}sultat.

La famille $(A_\ell\otimes_\alpha B_m)_{(\ell,m)\in \cL\times
  \cM}$ est clairement totale.

Finalement, soient $(\ell,m),\,(\ell',m')\in \cL\times \cM$ et $\alpha_\ell\otimes
\beta_m$, $\alpha_{\ell'}\otimes \beta_{m'}$ des tenseurs el{\'e}mentaires de
$A_\ell\otimes_\alpha B_m$, $A_{\ell'}\otimes_\alpha B_{m'}$. Alors $(\alpha_\ell\otimes
\beta_m) (\alpha_{\ell'}\otimes \beta_{m'})=\alpha_\ell \alpha_{\ell'}\otimes
\beta_m \beta_{m'}\in A_{\ell\wedge \ell'}\otimes_\alpha B_{m\wedge m'}$. On en
deduit que $(A_\ell\otimes_\alpha B_m)(A_{\ell'}\otimes_\alpha B_{m'})\subset
A_{\ell\wedge \ell'}\otimes_\alpha B_{m\wedge m'}$.

\end{proof}

\begin{remark}Par la proposition pr{\'e}cedente on a une
  $C^*$-alg{\`e}bre $\cL\times \cM$-gradu{\'e}e $(\gA\otimes_\alpha \gB,(A_\ell\otimes_\alpha B_m)_{(\ell,m)\in
    \cL\times \cM})$ et on a montr{\'e} que $A_\ell\otimes_\alpha B_m\subset (A_\ell\otimes_\alpha \gB)\cap (\gA\otimes_\alpha
  B_m)$ pour tout $\ell\in \cL$ et $m\in \cM$. Montrons l'{\'e}galit{\'e}.

 Par le
  lemme \ref{ten} la $C^*$-alg{\`e}bre $(\gA\otimes_\alpha
  \gB,(A_\ell\otimes_\alpha \gB)_{\ell\in \cL})$ est $\cL$-gradu{\'e}e et $(\gA\otimes_\alpha
  \gB,(\gA\otimes_\alpha B_m)_{m\in \cM})$ est $\cM$-gradu{\'e}e. Puisque
  $\sum\limits_{(\ell,m)\in \cL\times \cM}(A_\ell\otimes_\alpha B_m)$ est
  dense dans $\gA\otimes_\alpha \gB$, en appliquant
  la proposition \ref{inter}, on a que $(\gA\otimes_\alpha \gB,((A_\ell\otimes_\alpha
  \gB)\cap (\gA\otimes_\alpha B_m))_{(\ell,m)\in \cL\times \cM})$ est aussi une
  $C^*$-alg{\`e}bre gradu{\'e}e. L'{\'e}galit{\'e} $A_\ell\otimes_\alpha B_m=(A_\ell\otimes_\alpha \gB)\cap (\gA\otimes_\alpha
  B_m)$ r{\'e}sulte maintenant du corollaire \ref{sgr}.\\  
\end{remark}

\begin{remark}On reprend les hypoth{\`e}ses de la proposition
  \ref{tens}. Suite a la proposition \ref{q}
  on remarque que les morphismes de structure de la $C^*$-alg{\`e}bre
  $\cL\times \cM$-gradu{\'e}e $(\gA\otimes_\alpha \gB,(A_\ell\otimes_\alpha B_m)_{(\ell,m)\in
    \cL\times \cM})$ sont donn{\'e}s par
  $\varphi_{(\ell,m),(\ell',m')}=\varphi_{\ell,\ell'}\otimes_\alpha
  \varphi_{m,m'}$ pour $(\ell,m),(\ell',m')\in \cL\times \cM$ tels que
  $\ell\leq \ell'$ et $m\leq m'$ o{\`u} $\varphi_{\ell,\ell'}$, $\varphi_{m,m'}$
   sont les morphismes de structure de la $C^*$-alg{\`e}bre gradu{\'e}e
   $(\gA,(A_\ell)_{\ell\in \cL})$ et $(\gB,(B_m)_{m\in \cM})$
   respectivement. Plus precisement, si $\alpha_\ell\otimes
\beta_m$, $\alpha_{\ell'}\otimes \beta_{m'}$ sont des tenseurs el{\'e}mentaires de
$A_\ell\otimes_\alpha B_m$, $A_{\ell'}\otimes_\alpha B_{m'}$ respectivement avec $\ell\leq \ell'$
et $m\leq m'$ alors le morphisme 

$$\varphi_{\ell,\ell'}\otimes_\alpha
  \varphi_{m,m'}:A_{\ell'}\otimes_\alpha B_{m'}\to M(A_\ell\otimes_\alpha B_m)$$

est donn{\'e} par
 $$(\varphi_{\ell,\ell'}\otimes_\alpha
  \varphi_{m,m'})(\alpha_{\ell'}\otimes \beta_{m'})=
 \varphi_{\ell,\ell'}(\alpha_{\ell'})\otimes_\alpha
 \varphi_{m,m'}(\beta_{m'}):\alpha_{\ell}\otimes \beta_m\mapsto
 \alpha_{\ell'}\alpha_{\ell}\otimes \beta_{m'}\beta_m.$$

On peut voir facilement que les morphismes $\varphi_{(\ell,m),(\ell',m')}$ v{\'e}rifient
les conditions $a)$ et $b)$ de la proposition \ref{q}.
\end{remark} 

\begin{remark}Soient $\cL$ et $\cM$ des semi-treillis qui chacun poss{\`e}de
  un plus petit {\'e}l{\'e}ment not{\'e} $\ell_0$ et $m_0$ respectivement. Alors le
  semi-treillis produit $\cL\times \cM$ poss{\`e}de aussi un plus petit
  {\'e}l{\'e}ment not{\'e} $(\ell_0,m_0)$.

 Soit $(\gA\otimes_\alpha \gB,(A_\ell\otimes_\alpha B_m)_{(\ell,m)\in
    \cL\times \cM})$ une $C^*$-alg{\`e}bre $\cL\times \cM$-gradu{\'e}e. On
  suppose que les morphismes de structure
  $\varphi_{\ell,\ell'}$ avec $\ell\leq \ell'$ et $\varphi_{m,m'}$
  avec $m\leq m'$ de l'alg{\`e}bre
  $(\gA,(A_\ell)_{\ell\in \cL})$ et $(\gB,(B_m)_{m\in \cM})$ respectivement v{\'e}rifient
  $\varphi_{\ell,\ell'}^{-1}(A_\ell)=\{0\}$ et
  $\varphi_{m,m'}^{-1}(B_m)=\{0\}$.  Alors, par le
  corollaire \ref{min} les morphismes $\varphi_{\ell_0}:\gA\to
  M(A_{\ell_0})$ et $\varphi_{m_0}:\gB\to
  M(B_{m_0})$ sont injectifs. 

D'apr{\`e}s le corollaire \ref{mormult}
  le morphisme
  $$\varphi_{\ell_0}\otimes_{\min}\varphi_{m_0}:\gA\otimes_{\min}\gB\to
  M(A_{\ell_0}\otimes_{\min} B_{m_0})$$
   est injectif.\\
\end{remark}

On peut g{\'e}n{\'e}raliser ces r{\'e}sultats par recurrence pour une famille
finie de semi-treillis $(\cL^k)_{k=1}^n,$ $n\in \N$ et une famille
finie de $C^*$-alg{\`e}bres $(\gA^k)_{k=1}^n$ gradu{\'e}es par
$\cL^k$, $k\in \N$ respectivement. 

\section{Produits crois{\'e}s}

\begin{proposition}\label{proexa}
Soit $\cal L$ un semi-treillis et soit $(\frak A,{(A_i)}_{i\in\cal L})$ une $C^*$-alg{\`e}bre $\cal
L$-gradu{\'e}e, munie d'une action continue d'un groupe localement
compact $G$. Supposons que pour tout $i\in {\cal L}$ la
$C^*$-alg{\`e}bre $A_i$ est $G$-invariante \ie stable par
les automorphismes de $\frak A$, ${\alpha}_g\;,g\in G$.

Alors le produit
crois{\'e} maximal $({\frak
  A}{\rtimes}_{\alpha}G,{(A_{i}{\rtimes}_{\alpha}G)}_{i\in {\cal
    L}})$ et le produit crois{\'e} reduit $({\frak A}{\rtimes}_{r,{\alpha}}G,{(A_{i}{\rtimes}_{r,{\alpha}}G)}_{i\in {\cal L}})$ sont des $C^*$-alg{\`e}bres
$\cal L$-gradu{\'e}es.
\end{proposition}

\begin{proof} On va faire la d{\'e}monstration pour le produit crois{\'e}
  maximal puisque le m{\^e}me raisonment d{\'e}montre que le produit crois{\'e} reduit
 est aussi une $C^*$-alg{\`e}bre $\cal L$-gradu{\'e}e.\\

Soit $k\in \cal L$ et ${\cal L}_k=\{j\in
 {\cal L}|k\leq j\}$, ${\cal L}'_k=\{j\in {\cal L}|k\nleq j\}$. Puisque $(\frak A,{(A_i)}_{i\in\cal L})$ une
$C^*$-alg{\`e}bre $\cal L$-gradu{\'e}e par la proposition \ref{suitexacte} on obtient
que ${\frak A}= \overline {\gA_{{\cL}'_k}}\bm\oplus \overline
{\gA_{\cL_k}}$. Par le corollaire \ref{exactmax} on a donc une
suite exacte scind{\'e}e des produits crois{\'e}s maximaux,$$0\to
 \overline
 {\gA_{{\cL}'_k}}{\rtimes}_{\alpha}G\xrightarrow{i_{*}}{\frak
   A}{\rtimes}_{\alpha}G\
 \mathop{\rightleftarrows}\limits^{p_{*}}_{{\sigma}_{*}}\ \overline {\gA_{\cL_k}}{\rtimes}_{\alpha}G\to 0.$$
 Autrement dit on a 
$${\frak A}{\rtimes}_{\alpha}G\simeq
( \overline {\gA_{{\cL}'_k}}{\rtimes}_{\alpha}G)\bm\oplus (\overline {\gA_{\cL_k}}{\rtimes}_{\alpha}G)$$
On va d{\'e}montrer que :
\begin{enumerate}
\item[$i)$] $h_k:A_k{\rtimes}_{\alpha}G\to \gA{\rtimes}_{\alpha}G$ est
  injectif pour tout $k\in \cL$.
\item[$ii)$] La famille de $C^*$-alg{\`e}bres
${(A_{i}{\rtimes}_{\alpha}G)}_{i\in {\cal L}}$ est lin{\'e}airement ind{\'e}pendante, 
\item[$iii)$] $(A_{i}{\rtimes}_{\alpha}G) (A_{j}{\rtimes}_{\alpha}G)\subset A_{i\wedge
  j}{\rtimes}_{\alpha}G\;\;,\forall i,j\in \cal L,$

\item[$iv)$] $\bigcup\limits_{{\cal F}\in \F_\cL}{\frak A}_{\cal F}{\rtimes}_{\alpha}G\equiv \sum\limits_{i\in
  {\cal L}}A_{i}{\rtimes}_{\alpha}G$ est dense dans
${\frak A}{\rtimes}_{\alpha}G$.
\end{enumerate}
\begin{enumerate}
 \item[$i)$] Comme $A_k$ est un id{\'e}al de $\overline {\gA_{\cL_k}}$
    alors $A_k{\rtimes}_{\alpha}G$ est un id{\'e}al de $\overline
    {\gA_{\cL_k}}{\rtimes}_{\alpha}G$ et on a
  $$\xymatrix{
    A_k{\rtimes}_{\alpha}G \ar@{^{(}->}[r]^{}  \ar[rd]_{h_k} & \overline {\gA_{\cL_k}}{\rtimes}_{\alpha}G  \ar[d]^{\sigma_*} \\
    {} & \gA{\rtimes}_{\alpha}G  .
  }$$

\item[$ii)$] Soit $f\in \gA{\rtimes}_{\alpha}G$, $f=\sum\limits_{i\in \cal L} f_i=0$ avec $f_i\neq 0$ pour seulement un
 nombre fini de $i\in \cal L$. Notons $\Omega=\{j\in {\cal L}:f_j\neq 0\}$
 le support de $f$. Supposons que $\Omega\neq \emptyset$ et soit $m\in \Omega$ un {\'e}l{\'e}ment maximal alors
 $f_m\neq 0\Rightarrow ||f_m||\neq 0$.

On a $\cL_m\cap \Omega=\{m\}$ et puisque $p_*$ est un morphisme $||f_m||\leq
||\sum\limits_{i\in \cL}f_i||$. Mais $\sum\limits_{i\in \cal L}
f_i=0\Rightarrow f_m=0$, ce qui est absurde. 

\item[$iii)$] Soient $i,j\in \cal L$ et soit $f_i\in C_c(G,A_i)\,,\,f_j\in
C_c(G,A_j)$ alors par la d{\'e}finition du produit dans l'alg{\`e}bre
involutive $C_c(G,{\frak A})$, pour tout $s\in G$ on a $(f_i*f_j)(s)\in
A_{i\wedge j}$. Donc $C_c(G,A_i) C_c(G,A_j)\subset
C_c(G,A_{i\wedge j})$ et par continuit{\'e} du produit on a
$(A_i{\rtimes}_{\alpha}G) (A_j{\rtimes}_{\alpha}G)\subset (A_{i\wedge
  j}{\rtimes}_{\alpha}G)$.

\item[$iv)$] Par hypothese on a que l'alg{\`e}bre $\bigcup\limits_{{\cal F}\in
  \F_\cL}{\frak A}_{\cal F}\equiv \sum\limits_{i\in \cal L}A_i$  est dense dans la
  $C^{*}$-alg{\`e}bre $\frak A$. Puisque le produit crois{\'e} des limites
  inductives est la limite inductive des produits crois{\'e}s on a $iv)$.
\end{enumerate}

\end{proof}

\begin{remark}On reprend les hypoth{\`e}ses de la proposition
  \ref{proexa}. Les morphismes de structure de $\gA$,
  $\varphi_{i,j}:A_j\to M(A_i)$ avec $i\leq j$, sont
  {\'e}quivariants. Suite {\`a} la proposition \ref{mu}, on
  remarque que les morphismes de structure de la $C^*$ alg{\`e}bre gradu{\'e}e $({\frak
  A}{\rtimes}_{\alpha}G,{(A_{i}{\rtimes}_{\alpha}G)}_{i\in {\cal
    L}})$ sont les morphismes $$\varphi_{i,j}^{\alpha}:A_{j}{\rtimes}_{\alpha}G\to
M(A_{i}{\rtimes}_{\alpha}G)$$
 pour tout $i,j\in \cL$ tels que $i\leq
j$ satisfaisant $\varphi_{i,j}^{\alpha}(f)g=f*g$ pour tout $f\in
A_{j}{\rtimes}_{\alpha}G$, $g\in A_{i}{\rtimes}_{\alpha}G$.  

Egalement, pour le
produit crois{\'e} reduit $({\frak
  A}{\rtimes}_{r,{\alpha}}G,{(A_{i}{\rtimes}_{r,{\alpha}}G)}_{i\in
  {\cal L}})$ on remarque que ses morphismes de structure sont les
morphismes
 $$\varphi_{i,j}^{r,\alpha}:A_{j}{\rtimes}_{r,\alpha}G\to
M(A_{i}{\rtimes}_{r,\alpha}G)$$ satisfaisant $\varphi_{i,j}^{r,\alpha}(f)g=f*g$ pour tout $f\in
A_{j}{\rtimes}_{r,\alpha}G$, $g\in A_{i}{\rtimes}_{r,\alpha}G$.\\ 
\end{remark}

\begin{remark}Soit $(\gA,(A_i)_{i\in \cL})$ une $C^*$-alg{\`e}bre
  gradu{\'e}e. On suppose que $\cL$ a un plus petit {\'e}l{\'e}ment qu'on note
  $\ell_0$ et que les morphismes de
  structure $\varphi_{k,\ell}:A_\ell\to M(A_k)$ satisfont
  $\varphi_{k,\ell}^{-1}(A_k)=\{0\}$. Alors, par le corollaire \ref{min} on obtient
  que $\varphi_{\ell_0}:\gA\to M(A_{\ell_0})$ est injectif.

Soit $G$ un groupe
  localement compact agissant sur $\gA$ et pr{\'e}servent les $A_i$,
  notons $\alpha$ cette action. Par la proposition \ref{proexa} on a que ${\frak A}{\rtimes}_{r,{\alpha}}G$ est une
  $C^*$-alg{\`e}bre gradu{\'e}e. Le morphisme {\'e}quivariant $\varphi_{\ell_0}$
  est injectif, alors par la proposition \ref{mu}, le morphisme $\varphi_{\ell_0}^{r,\alpha}:{\frak A}{\rtimes}_{r,{\alpha}}G\to
  M(A_{\ell_0}{\rtimes}_{r,{\alpha}}G)$ est injectif.
\end{remark}

\chapter{Propri{\'e}t{\'e}s des $C^*$-alg{\`e}bres gradu{\'e}es}
On suppose dans la suite qu'on se donne un semi-treillis $\cal L$ et
une $C^*$-alg{\`e}bre $\cal L$-gradu{\'e}e $(\frak A,{(A_i)}_{i\in\cal L})$.

\section{Commutativit{\'e}}

\begin{proposition} La $C^*$-alg{\`e}bre $\gA$ est commutative si et seulement si les composantes $A_i$ sont
    commutatifs pour tout $i\in \cL$.
\end{proposition}
 
\begin{proof}C'est {\'e}vident que si $\gA$ est une $C^*$-alg{\`e}bre
    commutative alors toute sous-$C^*$-alg{\`e}bre de $\gA$ est
    commutative. Donc $A_i$ est commutative pour tout $i\in \cL$.

Supposons que $A_i$ est commutatif pour tout $i\in \cL$. Par densit{\'e}, il suffit de montrer
que pour tout $i,j\in \cL$ et $x_i\in A_i$, $y_j\in A_j$ on a
$x_i y_j=y_j x_i$. Mais $x_i y_j=\varphi_{i\wedge
  j,i}(x_i)\varphi_{i\wedge j,j}(y_j)$ et comme $A_{i\wedge j}$ est
une $C^*$-alg{\`e}bre commutative, l'alg{\`e}bre de multiplicateurs
$M(A_{i\wedge j})$ l'est aussi, donc 
$x_i y_j=\varphi_{i\wedge
  j,i}(x_i)\varphi_{i\wedge j,j}(y_j)\\=\varphi_{i\wedge
  j,j}(y_j)\varphi_{i\wedge j,i}(x_i)=y_j x_i$.

\end{proof}

Etudions le spectre de $C^*$-alg{\`e}bres
gradu{\'e}es commutatives. 
\begin{remark}\label{sp}Soit $(\gA,(A_\ell)_{\ell\in \cL})$ une $C^*$-alg{\`e}bre
  gradu{\'e}e et commutative. Soit $i\in \cL$ et $\chi_i$ un caract{\`e}re de
  $A_i$. Notons $\widetilde{\chi_i}$ son extension dans l'alg{\`e}bre des
  multiplicateurs  $M(A_i)$. Alors, on d{\'e}finit de fa{\c c}on 
  unique un caract{\`e}re de $\gA$ par $\chi=\widetilde{\chi_i}\circ
  \pi_i$ o{\`u} $\pi_i=\varphi_i\circ p_i$ est le morphisme $\gA\to
  M(A_i)$. On obtient ainsi une application continue injective 
  $\psi_i:\mbox{Sp}(A_i)\to \mbox{Sp}(\gA)$. 

Soit $\chi$ un caract{\`e}re
  de $\gA$ et supposons qu'il existe $i,j\in \cL$ tels que
  $\chi=\widetilde{\chi_i}\circ \pi_i=\widetilde{\chi_j}\circ
  \pi_j$. On peut supposer que $i\not \leq j$. Alors
  ${\pi_i}_{|_{A_j}}=0$ donc ${\widetilde{\chi_i}\circ
    \pi_i}_{|_{A_j}}=0$ or ${\widetilde{\chi_j}\circ
    \pi_j}_{|_{A_j}}=\chi_j\neq 0$ d'o{\`u} la contradiction. Autrement
  dit les images des
  $\psi_i$ sont disjoints.  

\end{remark}

\begin{proposition}\label{spectre}Soit $(\gA,(A_\ell)_{\ell\in \cL})$ une $C^*$-alg{\`e}bre
  gradu{\'e}e et commutative. Si $\cL$ est un bon semi-treillis alors $$\mbox{Sp}(\gA)=\bigcup\limits_{i\in \cL}\psi_i(\mbox{Sp}(A_i)).$$ 
\end{proposition}

\begin{proof}
Soit $\chi\in \mbox{Sp}(\gA)$. On consid{\`e}re l'ensemble $I=\{i\in
\cL;\chi_{|_{A_i}}\neq 0\}$. On remarque que $I$ est un
sous-semi-treillis de $\cL$ car si $i,j\in I$ prenons $a_i\in A_i$, $a_j\in
A_j$ tels que $\chi(a_i)\neq 0$ et $\chi(a_j)\neq 0$, alors puisque $\chi$ est un
morphisme de $C^*$-alg{\`e}bres on a $\chi(a_ia_j)=\chi(a_i)\chi(a_j)\neq
0$ donc $\chi_{|_{A_{i\wedge j}}}\neq 0$. Par hypoth{\`e}se et par la remarque
\ref{bon} l'ensemble $I$ admet un plus petit
{\'e}l{\'e}ment qu'on note $i_0$. On pose $\chi_{i_0}=\chi_{|_{A_{i_0}}}$,
alors $\chi_{i_0}\neq 0$. On va montrer que
$\chi=\psi_{i_0}(\chi_{i_0})$. Pour cela il suffit de montrer que pour tout $j\in
\cL$ et $a_j\in A_j$, $\chi(a_j)=\psi_{i_0}(\chi_{i_0})(a_j)$ car
alors par lin{\'e}arit{\'e} l'{\'e}galit{\'e} sera vraie pour tout $x\in \gA_\cL$ et
par densit{\'e} pour tout $x\in \gA$. 

Soit $j\in \cL$. Si $i_0\nleq j$, puisque $i_0$ est le plus petit
{\'e}l{\'e}ment de $I$, alors $j\notin I$ et par consequent
$\chi_{|_{A_j}}=0$. D'autre part si $a_j\in A_j$, alors
$\psi_{i_0}(\chi_{i_0})(a_j)=\widetilde{\chi_{i_0}}\circ
\pi_{i_0}(a_j)=0$. 

Si $i_0\leq j$, puisque $\chi_{i_0}\neq 0$ il
existe $a_{i_0}\in A_{i_0}$ tel que
$\chi_{i_0}(a_{i_0})=\chi(a_{i_0})=1$\\ et comme
$\pi_{i_0}(a_{i_0})=a_{i_0}$ alors
$\psi_{i_0}(\chi_{i_0})(a_{i_0})=1$. Soit $a_j\in A_j$, alors

$\begin{aligned}\psi_{i_0}(\chi_{i_0})(a_j)&=\psi_{i_0}(\chi_{i_0})(a_j)\psi_{i_0}(\chi_{i_0})(a_{i_0})=\psi_{i_0}(\chi_{i_0})(a_ja_{i_0})=\chi_{i_0}(a_ja_{i_0})\\&=\chi(a_ja_{i_0})=\chi(a_j)\chi(a_{i_0})=\chi(a_j)\end{aligned}$\\
d'o{\`u} le r{\'e}sultat.
\end{proof}
  
\begin{remark}\label{morsp}Soit $\cL$ un bon semi-treillis et
 $(\gA,(A_\ell)_{\ell\in \cL})$ une $C^*$-alg{\`e}bre gradu{\'e}e,
  commutative. On suppose que les morphismes de structure
  $\varphi_{i,j}$ sont non-d{\'e}g{\'e}n{\'e}r{\'e}s pour tout $i,j\in \cL$ tels que
  $i\leq j$.

 Soit $\cM$ un sous-semi-treillis de $\cL$. On a
  $\overline{\gA_\cM}\subset \gA$. Supposons que $\overline{\gA_\cM}\gA=\gA$,
  autrement dit que tout {\'e}l{\'e}ment de $\cL$ est major{\'e} par un {\'e}l{\'e}ment de
  $\cM$. A l'inclusion
  $\iota:\overline{\gA_\cM}\to \gA$ correspond alors l'application continue
  $\Psi:\mbox{Sp}(\gA)\to \mbox{Sp}(\overline{\gA_{\cM}})$ donn{\'e}e par
  $\chi\mapsto \chi\circ \iota$. Etudions cette application {\`a} la
  lumi{\`e}re des d{\'e}compositions $\mbox{Sp}(\gA)=\bigcup\limits_{i\in
  \cL}\psi_i(\mbox{Sp}(A_i))$ et comme un sous-semi-treillis d'un bon
semi-treillis est un bon semi-treillis {\'e}crivons de m{\^e}me 
  $\mbox{Sp}(\overline{\gA_{\cM}})=\bigcup\limits_{m\in
  \cM}\rho_m(\mbox{Sp}(A_m))$. 

Soit $\chi\in \mbox{Sp}(\gA)$. Il existe un unique $i\in \cL$ tel que
$\chi=\psi_i(\chi_i)\in \psi_i(\mbox{Sp}(A_i))$. Soit $\cL_i=\{j\in
\cL;i\leq j\}$. Comme $\chi_{|_{A_i}}\neq 0$ et les morphismes $\varphi_{i,j}$ sont
non-d{\'e}g{\'e}n{\'e}r{\'e}s pour tout $j\in \cL$ tel que $i\leq j$, $\chi$ est non
nul sur $A_jA_i$ donc $\chi_{|_{A_j}}\neq 0$. En d'autres termes on a  $\cL_i=\{j\in
\cL;i\leq j\}=\{j\in \cL;\chi_{|_{A_j}}\neq 0\}$. Posons
$\cM_i=\cM\cap \cL_i=\{m\in
\cM;i\leq m\}=\{j\in \cM;\chi_{|_{A_j}}\neq 0\}$, alors $\cM_i$ est un
sous-semi-treillis finissant de $\cM$
. Puisque $\cL$ est un bon semi-treillis, $\cM$ admet un plus
petit {\'e}l{\'e}ment (remarque \ref{bon}). Soit $m_0\equiv m(i)=\inf\cM_i$. Alors
$\chi_{|_{A_{m_0}}}=\chi_{m_0}\neq 0$. 

 Le morphisme
$\varphi_{i,m_0}:A_{m_0}\to M(A_i)$ d{\'e}termine une application
continue $\Psi_i:\mbox{Sp}(A_i)\to \mbox{Sp}(A_{m_0})$ qui est donn{\'e}e
par $\Psi_i(\chi_i)=\widetilde{\chi_i}\circ \varphi_{i,m_0}$. Notons
$\widetilde{\varphi_{i,m_0}}$ l'extension de $\varphi_{i,m_0}$ {\`a}
l'alg{\`e}bre des multiplicateurs $M(A_{m_0})$. Le diagramme
 
$$\xymatrix{
     \overline{\gA_\cM}\ar[r]^{}  \ar[d]_{} &  \gA \ar[d]^{} \\
     M(A_{m_0}) \ar[r]^{} & M(A_i)
  }$$
est commutatif (on le v{\'e}rifie sur les composantes de $\overline{\gA_\cM}$). Par cons{\'e}quent le diagramme  

$$\xymatrix{
     \mbox{Sp}(A_i)\ar[r]^{\Psi_i}  \ar[d]_{\psi_i} &  \mbox{Sp}(A_{m_0}) \ar[d]^{\rho_{m_0}} \\
     \mbox{Sp}(\gA) \ar[r]^{\Psi} & \mbox{Sp}(\overline{\gA_{\cM}})
  }$$ est commutatif. En d'autres termes l'application
  $\Psi$ est entierement d{\'e}crite par l'application $\Psi_i$.
\end{remark}

\section{Nucl{\'e}arit{\'e}}

\begin{proposition}La $C^*$-alg{\`e}bre $\gA$ est nucl{\'e}aire si et
  seulement si les composantes $A_i$ sont nucl{\'e}aires pour tout $i\in \cL$. 
\end{proposition}

\begin{proof}Supposons que pour tout $i\in \cL$, $A_i$ est
  nucl{\'e}aire. Soit $B$ une $C^*$-alg{\`e}bre. Par le lemme
  \ref{ten}, $(\gA\otimes_{\max}B,(A_i\otimes_{\max}B)_{i\in \cL})$ et\\
  $(\gA\otimes_{\min}B,(A_i\otimes_{\min}B)_{i\in \cL})$ sont des
  $C^*$-alg{\`e}bres $\cL$-gradu{\'e}es. 

Notons $\rho:\gA\otimes_{\max}B\to
  \gA\otimes_{\min}B$ l'homomorphisme naturel. Pour tout $i\in \cL$,
  $\rho(A_i\otimes_{\max}B)\subset A_i\otimes_{\min}B$. Donc $\rho$
  est un homomorphisme de $C^*$-alg{\`e}bres gradu{\'e}es. Notons $\rho_i$ la
  restriction de $\rho$ sur $A_i\otimes_{\max}B$. Par hypoth{\`e}se
  $\rho_i:A_i\otimes_{\max}B\to A_i\otimes_{\min}B$ est un
  isomorphisme de $C^*$-alg{\`e}bres donc en appliquant la proposition
  \ref{injectif} on a que $\rho$ est un isomorphisme.

Supposons que $\gA$ est une $C^*$-alg{\`e}bre nucl{\'e}aire. Soit $i\in \cL$,
alors $\overline{\gA_{\cL_i}}$ qui est un quotient de $\gA$ est
nucl{\'e}aire. Donc $A_i$ qui est un id{\'e}al ferm{\'e} de
$\overline{\gA_{\cL_i}}$ est nucl{\'e}aire.
\end{proof}

\section{Exactitude}

\begin{proposition}
La $C^*$-alg{\`e}bre $\gA$ est exacte si et seulement si $A_i$ est exacte pour
tout $i\in \cL$.
\end{proposition}

\begin{proof}Si $\gA$ est exacte, comme pour tout $i\in \cL$, $A_i$ est une
  sous-$C^*$-alg{\`e}bre de $\gA$, alors $A_i$ est exacte.

Supposons que $A_\ell$ est exacte pour tout $\ell\in \cL$. Soit $0\to
  A\xrightarrow{i}B\xrightarrow{p}C\to 0$ une suite exacte
de $C^*$-alg{\`e}bres. Par le lemme \ref{ten} les
$C^*$-alg{\`e}bres $(A\otimes_{\min}\gA,(A\otimes_{\min}A_\ell)_{\ell\in \cL})$,
$(B\otimes_{\min}\gA,(B\otimes_{\min}A_\ell)_{\ell\in \cL})$ et 
$(C\otimes_{\min}\gA,(C\otimes_{\min}A_\ell)_{\ell\in \cL})$ sont $\cL$-gradu{\'e}es. On v{\'e}rifie facilement que les
morphismes induits $A\otimes_{\min}\gA\xrightarrow{i\otimes_{\min}
  Id_{\gA}}B\otimes_{\min} \gA$ et $B\otimes_{\min}
\gA\xrightarrow{p\otimes_{\min} Id_{\gA}}C\otimes_{\min}\gA$ sont des
morphismes de $C^*$-alg{\`e}bres gradu{\'e}es. Comme $A_\ell$ est exacte la suite $$0\to A\otimes_{\min}A_\ell\xrightarrow{i\otimes_{\min} Id_{A_\ell}}B\otimes_{\min}
A_\ell\xrightarrow{p\otimes_{\min} Id_{A_\ell}}C\otimes_{\min}A_\ell\to 0$$ est
exacte donc par le corollaire \ref{exactgr} la suite 
$$0\to A\otimes_{\min}\gA\xrightarrow{i\otimes_{\min} Id_{\gA}}B\otimes_{\min}
\gA\xrightarrow{p\otimes_{\min} Id_{\gA}}C\otimes_{\min}\gA\to 0$$ est
exacte. Cela prouve que $\gA$ est exacte.

\end{proof}

\section{K-th{\'e}orie}
  L'inclusion naturelle de $A_i\hookrightarrow \gA$ pour
  tout $i\in \cL$ d{\'e}finit un morphisme de groupes $K$, $K_k(A_i)\to
  K_k(\gA),\,k=0,1$. Donc on en deduit un morphisme $\Phi$ des groupes $K_0$ et $K_1$,
$$\Phi:\bigoplus\limits_{i\in \cL}K_k(A_i)\to K_k(\gA),\,k=0,1.$$

\begin{proposition} $\Phi$ est un isomorphisme.

\end{proposition}

\begin{proof}
Soit $\cF\in \F_\cL$. On va montrer par recurrence sur le nombre
d'{\'e}l{\'e}ments de $\cF$ que $\Phi_\cF:\bigoplus\limits_{i\in \cF}K_k(A_i)\to
K_k(\gA_\cF),\,k=0,1$ est un isomorphisme de groupes. 

C'est clair si $\cF$
poss{\`e}de un seul {\'e}l{\'e}ment. Supposons que $\cF$ a au moins deux {\'e}l{\'e}ments
et que pour tout sous-semi-treillis distinct de $\cF$, $\cF'$ le
r{\'e}sultat est
vrai.

Soit $\ell\neq \min\cF$, puisque $\gA_\cF$ est une $C^*$-alg{\`e}bre
$\cF$-gradu{\'e}e, par la proposition \ref{suitexacte} on a une suite
exacte scind{\'e}e de $C^*$-alg{\`e}bres $$0\to
 \gA_{{\cF}'_{\ell}}\xrightarrow{i_{\ell}}{\gA_\cF}
 \mathop{\rightleftarrows}\limits^{p_{\ell}}_{{\sigma}_{\ell}}\
 \gA_{\cF_{\ell}}\to 0.$$
   Donc, on en deduit une suite exacte scind{\'e}e de groupes
$$0\to
 K_k(\gA_{{\cF}'_{\ell}})\xrightarrow{i}K_k({\gA_\cF})
 \mathop{\rightleftarrows}\limits^{p}_{\sigma}\
 K_k(\gA_{\cF_{\ell}})\to 0,\,k=0,1,$$
et on a le diagramme suivant:
 $$\xymatrix{
           0 \ar[r]^{} & K_k(\gA_{{\cF}'_{\ell}}) \ar[r]^{i}
            & K_k({\gA_\cF}) \ar@<2pt>[r]^{p}
            & K_k(\gA_{\cF_{\ell}})
           \ar@<2pt>[l]^{\sigma} \ar[r]^{} 
            &0 \\
           0 \ar[r]^{} & \bigoplus\limits_{j\in \cF'_{\ell}}K_k(A_j)
           \ar[r]^{i'} \ar[u]_{\Phi_{\cF'_{\ell}}} &
           \bigoplus\limits_{j\in \cF}K_k(A_j) \ar@<2pt>[r]^{p'} \ar[u]_{\Phi_\cF}&
           \bigoplus\limits_{j\in \cF_{\ell}}K_k(A_j)
           \ar@<2pt>[l]^{\sigma'} \ar[r]^{} \ar[u]^{\Phi_{\cF_\ell}} & 0. 
         }$$
Par hypoth{\`e}se de recurrence, on a que $\Phi_{\cF'_{\ell}}$ et $\Phi_{\cF_\ell}$
           sont des isomorphismes, donc $\Phi_\cF$ l'est aussi.

Puisque la $K$-th{\'e}orie preserve la limite inductive, on a le m{\^e}me r{\'e}sultat
pour un semi-treillis infini $\cL$.
\end{proof}

\newpage
\quad \thispagestyle{empty}

\newpage

\chapter{Exemples de $C^*$-alg{\`e}bres gradu{\'e}es}

\section{Semi-treillis de sous-groupes}

Soit $G$ un groupe localement compact. L'ensemble $\cG$ des sous-groupes ferm{\'e}s de $G$ est un treillis (complet) pour l'inclusion: si $(H_i)_{i\in I}$ est une famille \emph{quelconque} de sous-groupes ferm{\'e}s de $G$, leur intersection $\bigcap_{i\in I}H_i$ est bien le plus grand {\'e}l{\'e}ment de $\{H\in \cG;\ \forall i\in I,\ H\subset H_i\}$; le plus petit {\'e}l{\'e}ment de $\{H\in \cG;\ \forall i\in I,\ H_i\subset H\}$ est l'adh{\'e}rence du sous-groupe engendr{\'e} par les $H_i$.

Pour $H\in \cG$, posons $A_H=C_0(G/H)$. Si $H,K\in \cG$ sont tels que $K\subset H$, toute classe {\`a} gauche $x\in G/K$ est contenue dans une classe {\`a} gauche $p_{K,H}(x)\in G/H$, d'o{\`u} une  application continue $p_{K,H}:G/K\to G/H$. On en d{\'e}duit un homomorphisme $\varphi_{K,H}:A_H=C_0(G/H)\to C_b(G/K)=M(A_K)$.

\begin{lemma}\label{conequiv} Soient $H,K\in \cG$. Les conditions suivantes  sont {\'e}quivalentes:\begin{enumerate}
\renewcommand{\theenumi}{\roman{enumi}}
\renewcommand{\labelenumi}{\rm (\theenumi)}
\item On a $\varphi_{H\cap K,H}(A_H)\varphi_{H\cap K,K}(A_K)\subset A_{H\cap K}$.
\item L'application $q:x\mapsto (p_{H\cap K,H}(x),p_{H\cap K,K}(x))$ de $G/(H\cap K)$ dans $G/H\times G/K$ est ferm{\'e}e.
\item Le sous-ensemble $HK=\{xy;\ x\in H,\ y\in K\}$ est ferm{\'e} dans
  $G$ et l'application produit $H\times K\to HK$, $(x,y)\mapsto xy$, est ouverte.
\item $p_H(K)$ est ferm{\'e} (dans $G/H$) et l'application $K\to p_H(K)$
  obtenue par restriction de $p_H$ est ouverte.
\item $p_K(H)$ est ferm{\'e} (dans $G/K$) et l'application $H\to p_K(H)$
  obtenue par restriction de $p_K$ est ouverte.
\end{enumerate}
\end{lemma}

\begin{proof}
Tout d'abord on remarque que $q$ est une application injective.

Montrons $(ii)\Rightarrow (i)$. Puisque $q$ est continue, injective
et ferm{\'e}e, c'est une application propre donc pour tout $C$ sous-espace
compact de $G/H\times G/K$, $q^{-1}(C)$ est un sous-espace compact de
$G/(H\cap K)$. Soit $f\in C_0(G/H),\,g\in C_0(G/K)$ et supposons que
${||f||}_{\infty}\leq 1,\;{||g||}_{\infty}\leq 1$, on va montrer
que $fg=\varphi_{H\cap K,H}(f)\varphi_{H\cap K,K}(g)\in C_0(G/(H\cap
K))$. On doit montrer que pour tout
$\varepsilon>0$ le sous ensemble $M:=\{x\in G/(H\cap K):|(fg)(x)|\geq \varepsilon\}$ de $G/(H\cap K)$ est compact. 

Soient $C_1=\{x\in G/H:|f(x)|\geq \varepsilon\}$, $C_2=\{x\in
G/K:|g(x)|\geq \varepsilon\}$. Par hypoth{\`e}se, ce sont des
sous-ensembles compacts de $G/H,\,G/K$. En utilisant le fait que
$\varphi_{H\cap K,H}$ et $\varphi_{H\cap K,K}$ sont des morphismes de
$C^*$-alg{\`e}bres on a $M\subset
q^{-1}(C_1\times C_2)$. Puisque $C_1\times C_2$ est
compact et $q$ est propre, alors $q^{-1}(C_1\times C_2)$ est compact
et par suite $M$ est compact.\\

Montrons $(i)\Rightarrow (ii)$. Montrons en fait que $q$ est
propre. Soit $C$ un compact de $G/H\times G/K$ et montrons que
$q^{-1}(C)$ est compact. Il existe $C_1,\,
C_2$ compacts de $G/H$ et $G/K$ tels que $C\subset C_1\times
C_2$. Comme $q^{-1}(C)$ est un
ferm{\'e} il suffit de montrer que $q^{-1}(C_1\times C_2)$ est un
compact. 

Soit $f\in C_0(G/H)$ telle que $f=1$ sur $C_1$ et soit $g\in
C_0(G/K)$ telle que $g=1$ sur $C_2$. Par hypoth{\`e}se on a $f g=\varphi_{H\cap K,H}(f)\varphi_{H\cap K,K}(g)\in
C_0(G/(H\cap K))$. L'ensemble $S=\{x\in G/(H\cap K):\,(f g)(x)=1\}$ est
un compact de $G/(H\cap K)$, or $q^{-1}(C_1\times C_2)\subset S$ donc
$q^{-1}(C_1\times C_2)$ est compact.\\

Montrons que $(ii)\Leftrightarrow (iii)$. Soit $r:G\to G/H\times G/K$
l'application $x\mapsto (p_H(x),p_K(x))$. Puisque $q$ est injective,
$p:G\to G/(H\cap K)$ est surjective et ouverte et $r=q\circ p$, on a les {\'e}quivalences suivantes:
\begin{enumerate}
\item $q:G/(H\cap K)\to G/H\times G/K$ est ferm{\'e}e.
\item $q(G/(H\cap K))$ est ferm{\'e} et $G/(H\cap K)\to q(G/(H\cap K))$
  est un hom{\'e}omorphisme.
\item $r(G)$ est ferm{\'e} et $r:G\to r(G)$ est ouverte.
\end{enumerate}
\bigskip

Rappelons d'abord qu'un \emph{carr{\'e} cart{\'e}sien} (d'espaces topologiques) est un carr{\'e}
commutatif d'applications continues
$$\xymatrix{
          X \ar[r]^{} \ar[d]_{} & Y \ar[d]^{g} \\
           Z \ar[r]_{f} & T
         }$$
tel que l'application $X\to Y{\times}_{T}Z$ d{\'e}duite soit un
          hom{\'e}omorphisme de $X$ sur le produit fibr{\'e} $\{(y,z)\in
          Y\times Z;g(y)=f(z)\}$.

On consid{\'e}re le diagramme suivant  $$\xymatrix{
           G \ar[r]^{} \ar[d]_{p} \ar[rd]_{r} & G\times G \ar[d]^{t} \\
           G/(H\cap K) \ar[r]_{q} & G/H\times G/K
         }$$

Puisque les applications $p_H:G\to G/H$ et $p_K:G\to G/K$
           sont ouvertes, l'application $t:G\times G\to G/H\times G/K$ est
           ouverte.  En d'autres termes l'application naturelle $(G\times
           G)/(H\times K)\to G/H\times G/K$ est un hom{\'e}omorphisme.

$$\begin{aligned}\mbox{On a}\;r(G)&=q(G/(H\cap K))\\&=\{(x,y)\in
           G/H\times G/K;\,\exists g\in G:x=gH,\;y=gK\},\,\mbox{alors}\end{aligned}$$ 
           $$\begin{aligned}t^{-1}(r(G))&=\{(g_1,g_2)\in G\times G;\,\exists g\in
           G,\,\exists h\in H,\,k\in
           K:g_1=gh,\,g_2=gk\}\\&=\{(g_1,g_2)\in G\times G;\,\exists h\in H,\,k\in
           K:g_1^{-1}g_2=h^{-1}k\}\\&=\{(g_1,g_2)\in G\times
           G:g_1^{-1}g_2\in HK\}.\end{aligned}$$

On en d{\'e}duit que $HK$ est ferm{\'e} si et seulement si $t^{-1}(r(G))$ est
           ferm{\'e} dans $G\times G$, autrement dit si et seulement si
           $r(G)$ est ferm{\'e}.

Remarquons que le produit fibr{\'e}
de $G$ par $t^{-1}(r(G))$ au-dessus de $r(G)$ s'{\'e}crit 
$$\begin{aligned}G{\times}_{r(G)}t^{-1}(r(G))&=G{\times}_{G/H\times
    G/K}(G\times G)\\&=\{(x,y,z)\in G\times G\times
  G;r(x)=t(y,z)\}\\&=\{(x,y,z)\in G\times G\times G;x^{-1}y\in
  H,\,x^{-1}z\in K\}\end{aligned}$$

Pour $(g,h,k)\in G\times H\times K$ posons $s(g,h,k)=g$ et
$\varphi(g,h,k)=(gh,gk)\in t^{-1}(r(G))$.
L'application $\psi:G\times H\times K\to G{\times}_{r(G)}t^{-1}(r(G))$
donn{\'e} par $(g,h,k)\mapsto (g,gh,gk)$ est un hom{\'e}omorphisme. Donc le
diagramme suivant  est un carr{\'e} cart{\'e}sien

$$\xymatrix{
           G\times H\times K \ar[r]^{\phi} \ar[d]_{s} \ar[rd]_{} & t^{-1}(r(G)) \ar[d]^{t_{r(G)}} \\
           G \ar[r]_{r} & r(G).
         }$$

 Les applications du
           diagramme ci-dessus sont surjectives. De plus $t_{r(G)}$ est
           ouverte car $t$ l'est. Puisque $G\times
           H\times K$ est le produit fibr{\'e}
           $G{\times}_{r(G)}t^{-1}(r(G))$ alors $r$ est
           ouverte si et seulement si $\phi$ l'est. L'application
           $t':t^{-1}(r(G))\to G\times HK$ donn{\'e}e par $(g_1,g_2)\mapsto
           (g_1,g_1^{-1}g_2)$ est un hom{\'e}omorphisme. Donc
           l'application $\phi'=t'\circ \phi:G\times
           H\times K\to G\times HK$ d{\'e}finie par $\phi'(g,h,k)=(g,hk)$
           est ouverte si et seulement si l'application $H\times K\to
           HK$ est ouverte d'o{\`u} l'{\'e}quivalence $(ii)\Leftrightarrow (iii)$.\\

On va montrer que $(iii)\Leftrightarrow (iv)$. Tout d'abord on remarque que 
$$\begin{aligned}p_H^{-1}p_H(K)&=\{g\in
           G:p_H(g)\in p_H(K)\}=\{g\in G;\,\exists k\in
           K,\,gH=kH\}\\&=\{g\in G;\,\exists k\in K,\,\exists h\in
           H:g=kh\}=KH.\end{aligned}$$

 Donc $p_H(K)$ est ferm{\'e} dans $G/H$ si et seulement si $KH$ est ferm{\'e}
 et par passage aux {\'e}l{\'e}ments inverses de $KH$ on a que $p_H(K)$ ferm{\'e}
 si et seulement si $HK$ est ferm{\'e}.

D'autre part, l'application $K\times H\to K{\times}_{p_H(K)}KH$ donn{\'e}e
           par $(k,h)\mapsto (k,kh)$ est un hom{\'e}omorphisme. Donc on a
           un carr{\'e} cart{\'e}sien:
$$\xymatrix{
            K\times H \ar[r]^{} \ar[d]_{} & KH \ar[d]^{} \\
           K \ar[r]_{} & p_H(K).
         }$$  
Donc l'application $K\to p_H(K)$ est
           ouverte si et seulement si l'application $K\times H\to KH$
            est ouverte, ce qui a lieu  si et seulement si l'application $H\times K\to
            HK$ est ouverte.

Par symmetrie on obtient aussi que $(iii)\Leftrightarrow (v)$.
 
\end{proof}

Donnons-nous un sous-semi-treillis $\cL$ de $\cG$ et supposons que
tout $H,K\in \cL$ satisfont les conditions {\'e}quivalentes du lemme \ref{conequiv}. On d{\'e}duit du th{\'e}or{\`e}me~\ref{carctgrad}:

\begin{proposition}\label{gr} Il existe une alg{\`e}bre gradu{\'e}e $(\gA,(A_H)_{H\in \cL})$ dont les morphismes de structure soient les $\varphi_{K,H}$ d{\'e}finis ci-dessus.\newline\indent \hfill$\square$
\end{proposition}

L'alg{\`e}bre des multiplicateurs de $C_0(G)$ est l'alg{\`e}bre $C_b(G)$ des fonctions continues born{\'e}es sur $G$.

Pour $H\in \cL$. Notons $p_H:G\to G/H$ l'application quotient. L'application $f\mapsto f\circ p_H$ identifie $C_0(G/H)$ {\`a} une sous-$C^*$-alg{\`e}bre de $C_b(G)$.

On obtient des morphismes $\varphi_H$ qui sont non-d{\'e}n{\'e}g{\'e}r{\'e}s pour tout
$H\in \cL$ donc ils se prolongent {\`a} des morphismes $\widetilde{\varphi_H}$
d{\'e}finis sur l'alg{\`e}bre de multiplicateurs $C_b(G/H)$. Il est clair que
$\varphi_H=\widetilde{\varphi_K}\circ \varphi_{K,H}$ pour tout $H,K\in \cL$ tels
que $K\subset H$. Ceci implique que
$\varphi_H(f)\varphi_K(g)=\varphi_{H\cap K}(fg)$ pour tout $K,H\in
\cL$ et $f\in C_0(G/H),g\in C_0(G/K)$. On en deduit un homomorphisme $\varphi :\gA\to C_b(G)$.

\begin{proposition}  \label{injmult} L'homomorphisme $\varphi :\gA\to C_b(G)$ est injectif si et seulement si, pour tout $H,K\in \cL$ tels que $K\subset H$ et $H\ne K$ l'espace quotient $H/K$ n'est pas compact.

\end {proposition} 

Nous aurons besoin d'un lemme:

\begin{lemma} \label{K/H} Soient $K,H$ deux sous-groupes de $G$ avec $K\subset H$. 
\begin{enumerate} 
\item Si $H/K$ est compact on a $\varphi_{K,H}(C_0(G/H))\subset
  C_0(G/K)$.
\item Si $H/K$ n'est pas compact on a $\varphi_{K,H}(C_0(G/H))\cap C_0(G/K)=\{0\}$.
\end{enumerate}
\end{lemma}

Pour cela on d{\'e}montre le lemme suivant:

\begin{lemma}\label{compact} Soit $G$ un groupe localement compact et soient $H$ et $K$ des sous-groupes ferm{\'e}s de $G$ tels que $K\subset H$.
\begin {enumerate}
\item Si $H/K$ est compact, l'application $p_{K,H}$ est propre.
\item Si $H/K$ n'est pas compact, pour tout $a\in G/H$ la partie $p_{K,H}^{-1}(\{a\})$ n'est pas relativement compacte dans $G/K$.
\end{enumerate}
\end{lemma}
\begin{proof} \begin {enumerate}
\item Soit $C$ une partie compacte de $G/H$, puisque $G$ est
  localement compact et l'application $p_H$ est ouverte et surjective, il existe une partie compacte $A$ de $G$ telle que $p_H(A)=C$. On a $p_{K,H}^{-1}(C)=\{ax; \ a\in A,\ x\in H/K\}$, c'est l'image du compact $A\times H/K$, donc une partie compacte de $G/K$.
\item Soit   $g\in G$ dont la classe est $a$, l'application $x\mapsto gx$ est un hom{\'e}omor\-phisme de $G/K$ qui envoie $H/K$ dans $p_{K,H}^{-1}(\{a\})$. Or, comme $H/K$ est ferm{\'e} dans $G/K$ et n'est pas compact, il n'est pas relativement compact. 
\end{enumerate}
\end{proof}

\medskip {\it D{\'e}monstration du lemme \ref{K/H}.}\/ 

\begin {enumerate} 
\item Soit $f\in C_0(G/H)$. Pour tout $\varepsilon>0$, l'ensemble $\{x\in G/K;\ |f\circ p_{K,H}(x)|\ge \varepsilon\}=p_{K,H}^{-1}(\{x\in G/H;\ |f (x)|\ge \varepsilon\})$ est compact, donc $f\circ p_{K,H}\in C_0(G/K)$.
\item Soit $f\in C_0(G/K)\cap \varphi_{K,H}(C_0(G/H))$. Pour tout $a\in G/H$, l'application $f$, constante sur $p_{K,H}^{-1}(\{a\})$, y est donc nulle.\newline\indent \hfill$\square$
\end{enumerate}

\medskip {\it D{\'e}monstration de la proposition~\ref{injmult}.}\/ S'il existe $H,K$ avec $K\subset H$, $H\ne K$ et $H/K$ compact, alors pour $f\in C_0(G/H)$ on a $$\varphi (f-\varphi_{K,H}(f))=\varphi_H(f)-\varphi_K(\varphi_{K,H}(f))=0,$$ donc $\varphi $ n'est pas injective.

\medskip Supposons inversement que pour tout $H,K\in \cL$ tels que $K\subset H$ et $H\ne K$ l'espace quotient $H/K$ n'est pas compact.
Nous distinguons trois cas:
\begin{enumerate}
\item Si $\{e\}\in \cL$. Ce cas r{\'e}sulte imm{\'e}diatement du lemme~\ref{K/H} et de la prop.~\ref{injectivite}.

\item Si aucun $H\in \cL$ n'est compact, on pose $\cL'=\cL\cup  \{\{e\}\}$ et on revient au premier cas.

\item S'il existe $H\in \cL$ avec $H$ compact, alors pour tout $K\in \cL$, l'espace $H/H\cap K$ est compact. Cela n'est possible que si $H\cap K=H$. Donc $H$ est le plus petit {\'e}l{\'e}ment de $\cL$. Par a prop.~\ref{injectivite}, l'application $\pi_H:\gA\to C_b(G/H)\subset C_b(G)$ est injective.
\newline\indent \hfill$\square$
\end{enumerate}

\section{Exemple non-commutatif}
On en d{\'e}duit que si on a un semi-treillis de sous-groupes ferm{\'e}s v{\'e}rifiant les
conditions du lemme \ref{conequiv} et de la proposition \ref{injmult},
les sous-$C^*$-alg{\`e}bres $C_0(G/H)\rtimes G\subset \B(L^2(G))$ forment une
sous-$C^*$-alg{\`e}bre de  $\B(L^2(G))$ gradu{\'e}e ({\`a} l'aide de la
proposition \ref{proexa}).

\begin{example}\label{espvec}Soit $X$ un espace vectoriel de dimension
            finie et soit $\cG$ le treillis pour l'inclusion de tous les sous-espaces vectoriels
            de $X$. En dimension finie, un sous-espace vectoriel est ferm{\'e} et toute
            application lin{\'e}aire et surjective est ouverte. Donc la condition $(iii)$ du lemme
            \ref{conequiv} est satisfaite et on a pour tout
            $Y,Z\in \cG$ que $C_0(X/Y) C_0(X/Z)\subset
            C_0(X/(Y\cap Z))$. Par la proposition \ref{gr}, on a pour
            tout sous-semi-treillis $\cL$ de $\cG$ une
            $C^*$-alg{\`e}bre gradu{\'e}e 
$$\gA=\overline{\bigoplus\limits_{Y\in \cL}C_0(X/Y)}.$$
 De plus, un espace vectoriel n'est pas
            compact, donc la condition $b)$ du lemme \ref{compact} est
            aussi satisfaite et on peut considerer $\gA$ comme une
 sous-$C^*$-alg{\`e}bre de $C_b(X)=M(C_0(X))\subset \B(L^2(X))$ o{\`u} la
 deuxi{\`e}me inclusion est donn{\'e}e par l'application qui {\`a} une fonction $u\in
 C_b(X)$ associe un op{\'e}rateur de multiplication par $u$ d{\'e}fini sur
 $L^2(X)$.

 On remarque que $\gA$ est stable par
 translations. On consid{\`e}re le produit crois{\'e} $\gA\rtimes X$ pour
 l'action continue de $X$ donn{\'e}e par les translations. Par la proposition \ref{proexa} ce produit crois{\'e} est aussi une
 $C^*$-alg{\`e}bre $\cG$-gradu{\'e}e \ie 
$$\gA\rtimes X=\overline{\bigoplus\limits_{Y\in
    \cL}(C_0(X/Y)\rtimes X)}.$$ En appliquant une nouvelle fois
    la proposition \ref{injmult} et par l'exemple \ref{comt} on a\\
 $$\begin{aligned}\gA\rtimes X&=\overline{\bigoplus\limits_{Y\in
    \cL}(C_0(X/Y)\rtimes X)}&\subset M(C_0(X)\rtimes
    X)&=M(\K(L^2(X)))\\&=\B(L^2(X)).\end{aligned}$$

En d'autres termes $\gA\rtimes X$ est aussi une
   $C^*$-alg{\`e}bre d'op{\'e}rateurs d{\'e}finis sur l'espace de Hilbert $L^2(X)$.\\

\end{example}

L'exemple pr{\'e}cedent est une autre approche du produit crois{\'e} consid{\'e}r{\'e}
et {\'e}tudi{\'e} par M. Damak et
V. Georgescu dans \cite{dg1} et par V. Georgescu et A. Iftimovici dans
\cite{gi1} (voir en particulier th{\'e}or{\`e}me $3.12$).

\section{Exemples commutatifs}

\begin{example}On se donne un plan $\cP$ et $n$ droites
  $\delta_i,\,i=1,...,n$ de $\cP$. On d{\'e}finit $n$ formes lin{\'e}aires
  $f_i,\,i=1,...,n$ de $\cP$ telles que $\ker f_i=\delta_i,\,i=1,...,n$. On a un treillis d'espaces
  vectoriels $\cF'=\{\{0\},\delta_1,...,\delta_n,\cP\}$. Par l'exemple
  pr{\'e}cedent on a une $C^*$-alg{\`e}bre $\cF'$-gradu{\'e}e:
 $$\gB=\bigoplus\limits_{Y\in \cF'}C_0(\cP/Y)\subset C_b(\cP).$$ On se propose
  de d{\'e}terminer le spectre de la $C^*$-alg{\`e}bre commutative et unif{\`e}re
  $\gB$.

\bigskip

On commence par traiter le cas $n=1$. Alors le treillis $\cF'=\{\{0\},\delta_1,\cP\}$ et l'alg{\`e}bre $\gB=C_0(\cP)\bm\oplus
C_0(\cP/\delta_1)\bm\oplus \C$.

On consid{\`e}re le sous-treillis de $\cF'$,
$\cF_1=\{\{0\},\delta_1\}.$ Alors, on a un id{\'e}al $\gB_1=C_0(\cP)\bm\oplus
C_0(\cP/\delta_1)\subset C_b(\cP).$ Par le th{\'e}or{\`e}me de Gelfand-Naimark
on a $\gB_1 \simeq C_0(\mbox{Sp}\,\gB_1)$. 

Id{\'e}ntifions $\cP$ {\`a} $\R^2$ et $\delta_1$ {\`a} l' ``axe des x'' (\ie
$\{(x,y);\,y=0\}$). Alors un {\'e}l{\'e}ment de $C_0(\cP/\delta_1)$ est de la
forme $(x,y)\mapsto g(y)$ o{\`u} $g\in C_0(\R)$. Posons $\T=P^1(\R)=\R\cup \{\infty\}$ qui est hom{\'e}omorphe au
cercle. On a alors les inclusions $\psi_0:C_0(\cP)\hookrightarrow C_0(\T\times
\R)$ et $\psi_{\delta_1}:C_0(\cP/\delta_1)\hookrightarrow C_0(\T\times
\R)$ donn{\'e}es par:

$$\psi_0(f):(x,y)\mapsto
\left\{\begin{aligned}f(x,y)\;\mbox{si}\;x\in
    \R\\0\;\mbox{si}\;x=\infty\end{aligned}\right.\,\mbox{et}\,\psi_{\delta_1}(g):(x,y)\mapsto g(y).$$ 

On v{\'e}rifie facilement que les morphismes $\psi_0$ et $\psi_{\delta_1}$ satisfont les
hypoth{\`e}ses de la proposition \ref{rest}, \ie $$\psi_{0\wedge
  \delta_1}(f g)=\psi_0(fg)=\psi_0(f)\psi_{\delta_1}(g),\forall f\in C_0(\cP),g\in
C_0(\cP/\delta_1),$$
 o{\`u} le produit $fg=\varphi_{0,0}(f)\varphi_{0,\delta_1}(g)$ est donn{\'e} par $(fg)(x,y)=f(x,y)g(y)$ pour tout $(x,y)\in \T\times \R$. Donc il existe un unique
homomorphisme $\psi:\gB_1\to C_0(\T\times \R)$ dont la restriction {\`a}
$C_0(\cP)$ et $C_0(\cP/\delta_1)$ soit $\psi_0$ et $\psi_{\delta_1}$
respectivement. De plus, comme le diagramme
$$\xymatrix{
            \gB_1 \ar[r]^{\psi} \ar@{^{(}->}[d]_{\varphi} &
            C_0(\T\times \R) \ar[dl]^{}   \\
           C_b(\cP).    &  {}
         }$$
 est commutatif, le morphisme $\psi$ est injectif. 

L'alg{\`e}bre $\gB$ est obtenue en adjoignant une unit{\'e} {\`a} $\gB_1$, \ie
            $\gB=\widetilde {\gB_1}$ et le morphisme $\psi$ s'{\'e}tend
            de mani{\`e}re unique {\`a} un morphisme $\widetilde \psi:\gB\to
            C((\T\times \R)\cup \{\infty\})$ o{\`u} $(\T\times \R)\cup
            \{\infty\}$ est le compactifi{\'e} d'Alexandroff de $\T\times
            \R$. C'est {\'e}vident que le morphisme $\widetilde \psi$ est
            injectif et de plus il est surjectif d'apr{\`e}s le
            th{\'e}or{\`e}me de Stones-Weiestrass car l'alg{\`e}bre $\widetilde {\psi}(\gB)$
            separe les points de $(\T\times \R)\cup \{\infty\}$.
 Donc $\gB \simeq C((\T\times \R)\cup \{\infty\}),$ autrement dit le
            spectre de $\gB$ est hom{\'e}omorphe {\`a} $(\T\times \R)\cup
            \{\infty\}$ qui topologiquement est hom{\'e}omorphe {\`a} une
            sph{\`e}re pinc{\'e}e.\\

Pour $n\geq 2$, on commence par consid{\'e}rer un autre exemple d'alg{\`e}bres
gradu{\'e}es : consid{\'e}rons une base de l'espace vectoriel $E=\R^n$,
  $(e_1,...,e_n)$. Soit $I$ un sous-ensemble de $\{1,...,n\}$. Notons $P_I$ le sous-espace vectoriel engendr{\'e} par ${(e_i)}_{i\in I}$. Soit $\cF$ l'ensemble
  $\{P_I;I\subseteq \{1,...,n\}\}$, alors $\cF$ est un
  treillis dont le plus petit {\'e}l{\'e}ment est $\{0\}$. On
  remarque que l'application $I\mapsto P_I$ est un isomorphisme de treillis.

 On va plonger $\cP$ dans $E$ par l'application
  lin{\'e}aire $(f_1,...,f_n):\cP\to \R^n$. On a $\cP\cap P_I=\{0\}$ pour tout $I$ tel que $\dim P_I\leq n-2$. Notons
 $H_i=P_{\{1,...,n\}\backslash \{i\}}$, $i=1,...,n$ les hyperplans de $E$, alors 
  les droites $\delta_i=\cP\cap
  \{x_i=0\}=\cP\cap H_i$ avec $i=1,...,n$.

Puisque les conditions {\'e}quivalentes du lemme \ref{conequiv}
  sont satisfaites pour tout $H,\,K\in \cF$ on a par la proposition
  \ref{gr} une $C^*$-alg{\`e}bre $\cF$-gradu{\'e}e \\
$(\gA,(C_0(E/H))_{H\in
  \cF})$ dont les morphismes de structure soient les $\varphi_{K,H}$
  d{\'e}finis ci-dessus. Elle est donn{\'e}e par:
  $$\gA=\bigoplus_{H\in \cF}C_0(E/H).$$

On remarque que $\gA$ est commutative et unif{\`e}re donc par le
th{\'e}or{\`e}me de Gelfand-Naimark on a $$\gA\simeq C(\mbox{Sp}(\gA)).$$

Puisque les hypoth{\`e}ses de la proposition \ref{injmult} sont
satisfaites on peut considerer $\gA$ et par consequence
$C(\mbox{Sp}(\gA))$ comme une sous-$C^*$-alg{\`e}bre de l'alg{\`e}bre
$C_b(E)$. On appelle encore $\varphi$ cette inclusion. En d'autres
termes, $\mbox{Sp}(\gA)$ est un compactifi{\'e} de $E$.

 On va montrer que l'image de $\gA$ dans $C_b(E)$ est l'espace de fonctions continues sur le tore de
dimension $n$, ${\mathbb T}^n$, autrement dit que Sp $\gA=\mathbb
  T^n$.

Puisque $\T^n$ est un compactifi{\'e} de $E$, la $C^*$-alg{\`e}bre $C_0(E)$ est un
id{\'e}al essentiel de $C(\T^n)$. On obtient donc un morphisme injectif $p:C(\T^n)\to C_b(E)=M(C_0(E))$.

Notons $J$ l'ensemble $\{1,...,n\}$ et soit $I\subset J$. Notons
$I'=J\backslash I$ et soit $p_{P_I}$ l'application quotient 
$$p_{P_I}:E\to E/P_{I}\simeq P_{I'}.$$
Soient $\rho_E$ et $\rho_{P_{I'}}$ les inclusions naturelles de $E$
et $P_{I'}$ dans leurs compactifi{\'e}s $\T^n$ et $\T^{I'}$, donc par
continuit{\'e} on peut prolonger $p_{P_I}$ en une application
$\widetilde{p_{P_I}}:\T^n\to \T^{I'}$. Autrement dit, on obtient le
diagramme commutatif

$$\xymatrix{
            E \ar[r]^{p_{P_I}} \ar[d]_{\rho_E} & P_{I'} \ar[d]^{\rho_{P_{I'}}} \\
           \T^n  \ar[r]^{\widetilde {p_{P_I}}} & \T^{I'}. 
         }$$

Soit $H\in \cF$ alors $H=P_I$, $I\subseteq J$. Par le diagramme
pr{\'e}cedent on en deduit un morphisme $\psi_H:C_0(E/H)\to C(\T^n)$ et on
a le diagramme commutatif:

$$\xymatrix{
            C_0(E/H) \ar[r]^{\psi_H} \ar[d]_{\varphi_{0,H}} & C(\T^n) \ar[dl]^{p}   \\
           C_b(E).    &  {}
         }$$     
 Comme $p$ est injective, le morphisme $\psi_H=p^{-1}\circ \varphi_{0,H}$ satisfait les hypoth{\`e}ses de la proposition
            \ref{rest}, alors il existe un unique homomorphisme
            $\psi:\gA\to C(\T^n)$ dont la restriction sur $C_0(E/H)$
            soit $\psi_H$. Donc on a le diagramme

$$\xymatrix{
            \gA \ar[r]^{\psi} \ar@{^{(}->}[d]_{\varphi} & C(\T^n) \ar@{^{(}->}[dl]^{p}   \\
           C_b(E)    &  {}
         }$$
qui est commutatif et par consequent le morphisme $\psi$ est injectif.

Pour montrer qu'il est surjectif, par le th{\'e}or{\`e}me de Stone-Weiestrass,
il suffit de montrer que $\psi(\gA)$ s{\'e}pare les points de
$\T^n$. Soient $x,y\in \T^n$ tels que $x\neq y$, alors il existe $j\in
\{1,...,n\}$ tel que $x_j\neq y_j$ avec $x_j,y_j\in \T$. Donc il existe
$f\in C_0(\R)$ telle que $f(x_j)\neq f(y_j)$. Or l'application
$\psi_j(f):(z_1,...,z_n)\mapsto f(z_j)$ est un {\'e}l{\'e}ment de $\psi(\gA)$
et $\psi_j(f)(x)\neq \psi_j(f)(y)$ d'o{\`u} la surjectivit{\'e} voulue.\\

On revient dans notre cas de l'alg{\`e}bre $\gB$. On va montrer que $\gB$ est l'alg{\`e}bre des fonctions continues
d{\'e}finies sur le tore ${\mathcal T}_g$ {\`a} $g$ trous o{\`u} $g=E({n\over 2})$ qui, dans
 le cas o{\`u} $n$ est impair, est pinc{\'e} \ie deux de ses points sont identifi{\'e}s.\\

 A l'inclusion (propre) $\cP\to E$ correspond un homomorphisme
 (surjectif) $C_b(E)\to C_b(\cP)$. Puisqu'on a $\gA\subset C_b(E)$,
on obtient donc un morphisme $\Phi:\gA\to C_b(\cP)$. 

Pour tout $P_I;\,I\subseteq \{1,...,n\}$, soit $\Phi_{P_I}$ la restriction de $\Phi$ sur les composants de
    $\gA$, $C_0(E/P_I)$.

Si $\dim P_I\leq n-2$ on a
    $\cP\hookrightarrow E/P_I$ et par suite l'image par $\Phi_{P_I}$
    de $C_0(E/P_I)$ dans $C_b(\cP)$ est $C_0(\cP)$.

Si $P_I=H_i,\,i=1,...,n$ alors pour tout $i$, $E/H_i$ est hom{\'e}omorphe {\`a}
    $\cP/{\delta}_i$ et on a un isomorphisme donn{\'e} par $\Phi_{H_i}$:
    $C_0(E/H_i)\simeq C_0(\cP/\delta_i).$

Donc on a montr{\'e} que $\Phi(\gA)=\gB$.

\bigskip

D'autre part $\gA=C(\T^n)$ et en consid{\`e}rant $C(\overline{\cP})$ comme
une sous-alg{\`e}bre de $C_b(\cP)$ alors $\Phi$ est la r{\'e}striction de
$C(\T^n)$ {\`a} ${\cP}$ \ie $C(\overline{\cP})=\Phi(\gA)=\gB$. En d'autres
terms, $\mbox{Sp}(\gB)$ est l'adh{\'e}rence de $\cP$ dans $\T^n$.

\bigskip

D{\'e}crivons un peu plus cette adh{\'e}rence. Avec les hypoth{\`e}ses de d{\'e}part le plan $\cP$ sera de la forme suivante:
$$\cP=\{(x_1,...,x_n)\in \R^n;\,x_k=A_{ijk}x_i+B_{ijk}x_j\},$$
 avec $i,k,j\in
\{1,...,n\},i,j\,\mbox{fix{\'e}s},\,k\notin \{i,j\}$
et les coefficients r{\'e}els $A_{ijk},\,B_{ijk}$ sont non nuls et ils
 v{\'e}rifient $B_{ijk}A_{ij\ell}\neq A_{ijk}B_{ij\ell}$ pour tout
 $k,\ell\notin \{i,j\}$ distincts.

Notons par $Q$ l'adh{\'e}rence de $\cP$ dans $(\overline{\R})^n$. Alors
l'adherence de $\cP$ dans $\T^n$, $\overline{\cP}=\rho(Q)$ o{\`u} $\rho$ est la surjection canonique
$(\overline{\R})^n\to \T^n$. 

Soit $(a_1,...,a_n)\in Q$ alors il existe une suite de
vecteurs $(x_{1,k},...,x_{n,k})_k\in \cP$ telle que
$(x_{1,k},...,x_{n,k})\to (a_1,...,a_n)$ quand $k\to
\infty$.

Supposons qu'au moins deux co{\'e}fficients de $(a_1,...,a_n)$
soient finis, $a_i,\,a_j$, et soit $m\in \{1,...,n\}\backslash \{i,j\}$. Alors si $x_{i,k}\to a_i$ et $x_{j,k}\to
a_j$ quand $k\to \infty$ on obtient que
$x_{m,k}=A_{ijm}x_{i,k}+B_{ijm}x_{j,k}\to A_{ijm}a_i+B_{ijm}a_j$, donc $a_m$ est
fini. On en deduit que la limite de $(x_{1,k},...,x_{n,k})$ quand
$k\to \infty$ est dans $\cP$.

Supposons qu'un seul co{\'e}fficient $a_i$ est fini. Alors si $j\in
\{1,...,n\}\backslash \{i\}$ et $x_{i,k}\to a_i$ quand $k\to \infty$
et si $x_{j,k}\to a_j=\infty$, o{\`u} $\infty$ signifie $+\infty$
ou$-\infty$, pour $m\in \{1,...,n\}\backslash \{i,j\}$ on a $x_{m,k}\to
A_{ijm}a_i+B_{ijm}a_j=A_{ijm}a_i+B_{ijm}(\infty)=\infty$, o{\`u} le signe de l'infini
depend du signe de l'infini limite de $x_{j,k}$ et du signe de
$B_{ijm}$, \ie $a_m=B_{ijm}(a_j)$. Donc dans $Q$ on a des vecteurs dont un coefficient est fini
et les autres sont infinis. On les appelle les droites {\`a} l'infini et
puisque $i$ parcourt l'ensemble $\{1,...,n\}$ il y en a $2n$. On va
montrer que le nombre des points d'intersection de ces droites, qu'on
appelera points {\`a} l'infini, est {\'e}galement $2n$. Par consequent $Q$
aura la forme d'un polygone {\`a} $2n$ c{\^o}t{\'e}s.

En effet, choisissons une orientation de $\cP$ et un point $A\in \cP\backslash
\cup_{i=1}^n\delta_i$. On r{\'e}indexe les droites de la fa{\c c}on suivante:
on choisit d'appeller la premi{\`e}re droite qu'on rencontre dans le sens
 direct par $\delta_1$, la deuxi{\`e}me $\delta_2$ jusqu'{\`a} la ni{\`e}me
 $\delta_n$, \ie on obtient

\begin{figure}[http]
\begin{center}
\includegraphics[width=8.025cm,height=5.715cm]{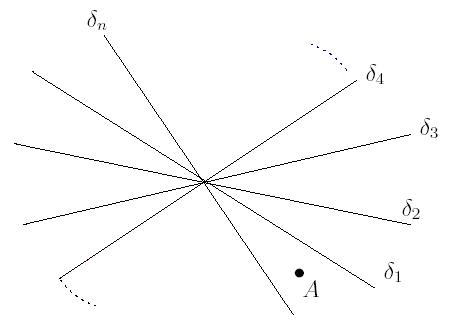} 
\end{center}
\end{figure}

On choisit les $n$ formes lin{\'e}aires $f_i,i=1,...,n$ telles que
$f_i(A)=1$ pour tout $i$. Donc, dans le secteur form{\'e} par les droites
$\delta_1$ et $\delta_n$ qui contient le point $A$, les signes du
point {\`a} l'infini seront tous positifs. Si on traverse la droite
$\delta_1$ et on passe dans le secteur form{\'e} par $\delta_1$ et
$\delta_2$ au-dessus du point $A$, les valeurs de $f_1$ deviennent
negatives, donc on a un changement du premier signe en negatif \ie le
point {\`a} l'infini sera $(-\infty,+\infty,...,+\infty)$. Avec la m{\^e}me 
procedure on ne trouve que $2n$ possibilit{\'e}s de signes qui sont de la forme
$(+,...,+,-,...,-)$ et $(-,...,-,+,...,+)$.

\begin{figure}[http]
\begin{center}
\includegraphics[width=14.11cm,height=7.02cm]{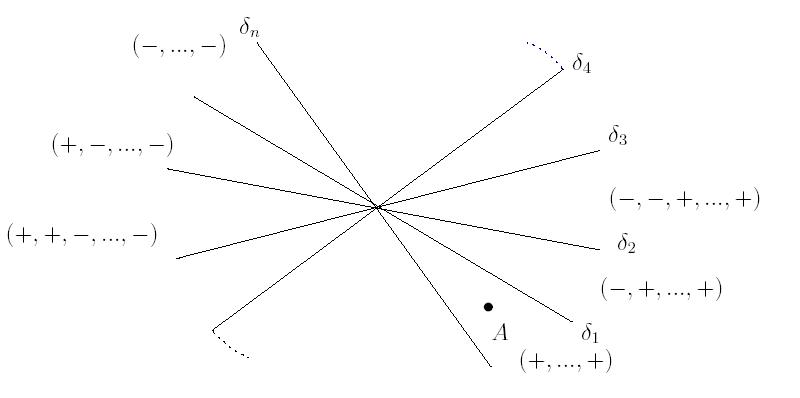} 
\end{center}
\end{figure}

Donc l'adh{\'e}rence de $\cP$ dans $(\overline{\R})^n$ est un polygone {\`a}
$2n$ c{\^o}t{\'e}s. En passant de $(\overline{\R})^n$ {\`a} $\T^n$, on id{\'e}ntifie
$+\infty$ et $-\infty$, ce qui revient {\`a} 
\begin{enumerate}
\item identifier chaque c{\^o}t{\'e} au c{\^o}t{\'e} oppos{\'e}

\item identifier tous les sommets.

\end{enumerate}

\begin{figure}[http]
\begin{center}
\includegraphics[width=10.568cm,height=6.88cm]{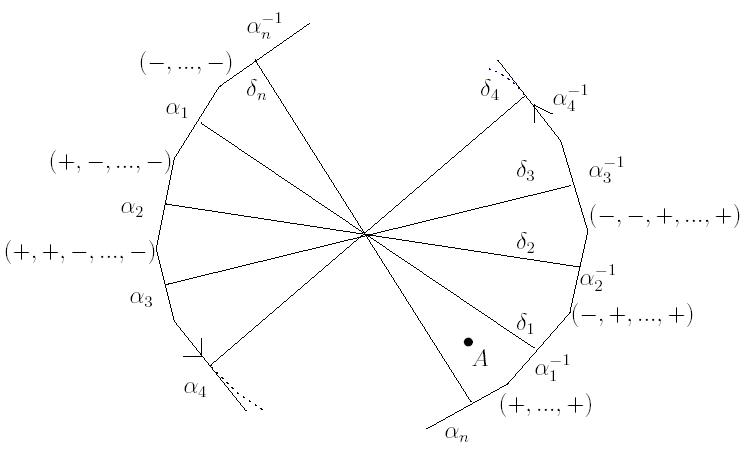} 
\end{center}
\end{figure}

La premi{\`e}re identification, donne lieu {\`a} la surface de Riemann $Y_g$ dont la
forme normale est
$\alpha_1\alpha_2...\alpha_n\alpha_1^{-1}\alpha_2^{-1}...\alpha_n^{-1}$
(voir par exemple dans \cite{ma}, \cite{spr}) qui est {\'e}quivalente {\`a}
$$\alpha_1\alpha_2\alpha_1^{-1}\alpha_2^{-1}...\alpha_{n-1}\alpha_n\alpha_{n-1}^{-1}\alpha_n^{-1}$$
si $n$ est pair et {\`a}
$$\alpha_1\alpha_2\alpha_1^{-1}\alpha_2^{-1}...\alpha_{n-2}\alpha_{n-1}\alpha_{n-2}^{-1}\alpha_{n-1}^{-1}\alpha_n\alpha_n^{-1}$$
si $n$ est impair. C'est donc un tore {\`a} $g$ trous o{\`u} $g=E({n\over
  2})$.

Indexons les sommets par $\Z/2n\Z$. En faisant l'identification des
c{\^o}t{\'e}s on identifie aussi le $j$-i{\`e}me sommet avec le sommet $j+n-1$
pour tout $j\in \Z/2n\Z$. En d'autres termes deux sommets sont
identifi{\'e}s si et seulement si ils sont dans la m{\^e}me orbite pour
l'addition de $(n-1)$ dans $\Z/2n\Z$. Or, si $n$ est pair
PGCD$(n-1,2n)=1$ et si $n$ est impair PGCD$(n-1,2n)=2$.

Donc l'application $Y_g\to \mbox{Sp}(\gB)$ est un hom{\'e}omorphisme
lorsque $n$ est pair.

Lorsque $n$ est impair, on passe de $Y_g$ {\`a} $\mbox{Sp}(\gB)$ en
id{\'e}ntifiant ces deux points, $\mbox{Sp}(\gB)$ est donc une surface de
genre $g$ pinc{\'e}e.\\

\end{example}

Appliquons la remarque \ref{morsp} {\`a} cet exemple.
\begin{remark} Notons $\cF$ le treillis donn{\'e} par
  $\cF=\{\{0\},\delta_1,...,\delta_n,\cP\}$ et $\cF'$ le treillis qui
  contient une droite de moins, par exemple choisissons
  $\cF'=\{\{0\},\delta_1,...,\delta_{n-1},\cP\}=\cF\backslash
  \{\delta_n\}.$ 

Soit
  $\gA_{\cF}=\bigoplus\limits_{H\in \cF}C_0(\cP/H)$ et 
  $\gA_{\cF'}=\bigoplus\limits_{K\in \cF'}C_0(\cP/K)$ les
  $C^*$-alg{\`e}bres gradu{\'e}es correspondantes. Elles sont commutatives et unif{\`e}res. On a
  $\gA_{\cF'}\subset \gA_{\cF}$ et par cons{\'e}quent on a une
  application $\Psi:\mbox{Sp}(\gA_{\cF})\to \mbox{Sp}(\gA_{\cF'})$. Notons
  $A_H=C_0(\cP/H)$ pour tout $H\in \cF$. 

Puisque $\cF$ et $\cF'$ sont finis, ce sont de bons treillis et on a $\mbox{Sp}(\gA_\cF)=\bigcup\limits_{H\in
  \cF}\mbox{Sp}(A_H)$ et $\mbox{Sp}(\gA_{\cF'})=\bigcup\limits_{K\in
  \cF'}\mbox{Sp}(A_K)$ (proposition \ref{spectre}). 

Soit $\chi\in \mbox{Sp}(\gA_\cF)$. Il existe
  un unique $H\in \cF$ tel que $\chi=\chi_H\in \mbox{Sp}(A_H)$. En appliquant la remarque
  \ref{morsp} on obtient que: Si $H\in \cF'$, l'image de $\chi_H$ par
  l'application $\Psi$ est lui-m{\^e}me. Si $H\in \cF\backslash \cF'$ \ie
  $H=\delta_n$, alors $\inf\{K\in \cF';\,\delta_n\subseteq K\}=\cP$ et
  l'image de $\chi_{\delta_n}$ par l'application $\Psi$ est dans
  $\mbox{Sp}(A_\cP)$ \ie c'est le point {\`a} l'infini.

En d'autres termes, on passe du spectre de $\gA_\cF$ {\`a} $\gA_{\cF'}$ en
contractant la droite $\cP/\delta_n$ en un seul point: le point {\`a} l'infini.\\    
\end{remark}

\begin{example}Soit $\cL$ un semi-treillis et $(A_i)_{i\in \cL}$ la famille
  de $C^*$-alg{\`e}bres d{\'e}finie par $A_i=\C$ pour tout $i\in \cL$. Pour $i,j\in \cL$
  tels que $i\leq j$ posons $\varphi_{i,j}={\mbox{Id}}_{\C}$. Ces morphismes v{\'e}rifient les propri{\'e}t{\'e}s $a)$ et $b)$ de
  la proposition \ref{q}, donc par le th{\'e}or{\`e}me \ref{carctgrad} il
  existe {\`a} isomorphisme pr{\`e}s une unique $C^*$-alg{\`e}bre $\cL$-gradu{\'e}e
  $(\gA,(A_i)_{i\in \cL})$ dont les morphismes de structure soient les
  $\varphi_{i,j}$. On veut ici calculer son spectre
  $\mbox{Sp}(\gA)$.

Montrons que l'ensemble des sous-semi-treillis finissants non vides de
$\cL$ est en bijection avec l'ensemble des caract{\`e}res de
$\gA$.

Notons $e_i$, $i\in \cL$ l'image de $1\in A_i$ dans $\gA$. C'est clair
que $e_i$ est un id{\'e}mpotent pour tout $i\in \cL$ donc si $\chi$ est
un caract{\`e}re de $\gA$ on a $\chi(e_i)\in \{0,1\}$. De plus, par la
d{\'e}finition du produit dans $\gA$ on a $e_ie_j=e_{i\wedge j}$ pour tout
$i,j\in \cL$.

Soit $\chi$ un caract{\`e}re de $\gA$. Posons $\cM_\chi=\{i\in
\cL;\chi(e_i)=1\}$. C'est un sous-semi-treillis de $\cL$ car si
$i,j\in \cM_\chi$ on a $\chi(e_{i\wedge
  j})=\chi(e_ie_j)=\chi(e_i)\chi(e_j)=1$ donc $i\wedge j\in \cM_\chi$. C'est un sous-semi-treillis
finissant car pour $i\in \cM_\chi$ et $j\in \cL$ tel que $i\leq j$ on
a $\chi(e_j)=\chi(e_i)\chi(e_j)=\chi(e_ie_j)=\chi(e_i)=1$ \ie $j\in \cM_\chi$.

Soit $\cM$ un sous-semi-treillis finissant non vide de $\cL$. D'apr{\`e}s
la proposition \ref{rest} il existe un
unique homomorphisme $\chi_\cM:\gA\to \C$ satisfaisant 
$\chi_\cM(\lambda e_i)=\left\{\begin{aligned}&\lambda,\;i\in \cM
\\&0,\;i\notin \cM\end{aligned}\right.$ pour tout $\lambda\in \C$ et
tout $i\in \cL$. En effet, il est facile de voir que, puisque $\cM$
est un sous-semi-treillis finissant, $\chi_\cM(\lambda
e_i)\chi_\cM(\mu e_j)=\chi_\cM(\lambda \mu e_{i\wedge j})$. Comme $\cM\neq
\emptyset$, $\chi_\cM$ est un caract{\`e}re de $\gA$.

Remarquant que tout caract{\`e}re est d{\'e}termin{\'e} par sa valeur sur $e_i$
($i\in \cL$) on a montr{\'e} que les applications $\chi\mapsto \cM_\chi$
et $\cM\mapsto \chi_\cM$ sont des bijections inverses l'une de
l'autre.

Remarquons enfin que, puisque $A_i=\C$, $\mbox{Sp}(A_i)$ a un seul
point $\mbox{Id}_\C$. On a
$\psi_i(\mbox{Id}_\C)(e_j)=\left\{\begin{aligned}&1,\;i\leq j
\\&0,\;i\nleq j\end{aligned}\right.$ \ie
$\psi_i(\mbox{Id}_\C)=\chi_{\cL_i}$. Il en r{\'e}sulte qu'un caract{\`e}re
$\chi$ de $\gA$ est dans $\bigcup\limits_{i\in
  \cL}\psi_i(\mbox{Sp}(A_i))$ si et seulement si $\cM_\chi$ est de la
forme $\cL_i$ \ie si et seulement s'il poss{\`e}de un plus petit
{\'e}l{\'e}ment. Finalement $\mbox{Sp}(\gA)=\bigcup\limits_{i\in
  \cL}\psi_i(\mbox{Sp}(A_i))$ si et seulement si $\cL$ est un bon
semi-treillis.     
\end{example}
${}$\\

Etudions un cas particulier de cet exemple. 

\begin{example}On consid{\`e}re l'ensemble
$\Q$ avec son ordre. Remarquons que $\Q$ {\'e}tant totalement ordonn{\'e},
c'est un semi-treillis mais ce n'est pas un \emph{bon} semi-treillis. On prend alors
$\cL=\Q$ et on va {\'e}tudier le spectre de la $C^*$-alg{\`e}bre gradu{\'e}e
correspondante $(\gA,(A_i)_{i\in \Q})$.

Notons $\cB_b(\R)$ la $C^*$-alg{\`e}bre des fonctions complexes
boreliennes born{\'e}es d{\'e}finies sur $\R$.

Pour tout $i\in \Q$ on d{\'e}finit une application $\tau_i:A_i\to \cB_b(\R)$
 par $\tau_i(\lambda e_i)=\lambda {\mathbf{1}}_{]-\infty,i]}$ o{\`u} $\fc$
 est la fonction caract{\'e}ristique. On v{\'e}rifie facilement que $\tau_i$
 est un morphisme injectif de $C^*$-alg{\`e}bres satisfaisant l'{\'e}galit{\'e} $\tau_{i\wedge
 j}(\lambda e_i\mu e_j)=\tau_i(\lambda e_i)\tau_j(\mu e_j)$ pour tout $i,j\in \Q$ et
 $\lambda,\mu\in \C$. Par cons{\'e}quent, suite {\`a} la proposition
 \ref{rest} il existe un unique morphisme $\tau:\gA\to \cB_b(\R)$ dont
 la restriction {\`a} $A_i$ soit $\tau_i$ pour tout $i\in \Q$. 

Montrons que $\tau$ est injectif. D'apr{\`e}s la proposition
\ref{injectif} il suffit de l'{\'e}tudier pour tout sous-semi-treillis
fini $\cF$ de $\Q$. Soit $\cF\in \F_\Q$. Supposons que $\cF$ poss{\`e}de $n$
 {\'e}l{\'e}ments. Puisque l'ordre est total on peut aussi supposer que 
 $\cF=\{i_1,...,i_n\}$ avec $i_1<i_2<...<i_n$. Soit
 $\tau_\cF$ la restriction de $\tau$ {\`a} $\gA_\cF$ et $x\in
 \ker\tau_\cF$. Alors $x=\sum\limits_{i\in \cF}\lambda_i e_i$ et
 $\sum\limits_{i\in \cF}\lambda_i\fc_{]-\infty,i]}(s)=0$ pour tout
 $s\in \R$. Pour $s$ tel que $i_{n-1}<s\leq i_n$ on trouve
 $\lambda_n=0$. Si on prend $s$ tel que
 $i_{n-2}<s\leq i_{n-1}$ on a $\lambda_{n-1}+\lambda_n=0$ et par
 cons{\'e}quent $\lambda_{n-1}=0$. On continue de la m{\^e}me mani{\`e}re et on
 d{\'e}montre que $\lambda_i=0$ pour tout $i\in \cF$. Autrement dit $x=0$
 et $\tau_\cF$ est injectif.\\

On peut id{\'e}ntifier donc $\gA$ avec son image par $\tau$. D{\'e}crivons cette
image. On va montrer que $\tau(\gA)=\cA$ o{\`u} $\cA$ est l'espace des
fonctions complexes boreliennes born{\'e}es d{\'e}finies sur $\R$ qui sont regl{\'e}es,
continues {\`a} gauche en tout point de $\R$, continues sur $\R\backslash
\Q$, nulles en $+\infty$ et qui admettent une limite en $-\infty$.

Montrons d'abord que $\cA$ est une sous-$C^*$-alg{\`e}bre de
$\cB_b(\R)$. Il suffit de montrer qu'elle est ferm{\'e}e. Ceci est clair
car une limite uniforme de fonctions regl{\'e}es (resp. continues {\`a} gauche
en tout point de $\R$, resp. continues sur $\R\backslash
\Q$, resp. nulles en $+\infty$, resp. qui admettent une limite en
$-\infty$) est regl{\'e}e (resp. continue {\`a} gauche
en tout point de $\R$, resp. continue sur $\R\backslash
\Q$, resp. nulle en $+\infty$, resp. admet une limite en
$-\infty$). 

 Pour tout $i\in \cL$, $\tau(e_i)\in \cA$ donc $\tau(\gA_\cL)\subset
\cA$. Puisque $\cA$ est ferm{\'e} et $\gA_\cL$ est dense dans $\gA$ on a
$\tau(\gA)\subset \cA$.

Soit $f\in \cA$. Supposons que la limite de $f$ en $-\infty$ est {\'e}gale
{\`a} $\lambda$, donc soit $\varepsilon>0$, il existe $a\in \Q$ tel que
pour tout $x\leq a$ on ait $|f(x)-\lambda|<\varepsilon$. D'autre part
$f$ est nulle en $+\infty$ donc il existe $A\in \Q$ tel que pour tout
$x\geq A$ on ait $|f(x)|<\varepsilon$. Autrement dit, pour tout $x\in
]-\infty,a]\cup ]A,+\infty[$ on a
$|f(x)-\lambda\fc_{]-\infty,a]}(x)|<\varepsilon$.

 Puisque $f$ est une fonction regl{\'e}e, elle est regl{\'e}e sur $[a,A]$ et
 il existe une fonction en escaliers $g$ sur $[a,A]$ telle que
 $||f-g||_{\infty}=\sup\{|f(t)-g(t)|,t\in [a,A]\}<\varepsilon$. Il
 existe donc une subdivision 
 $i_0,i_1,...,i_s$ de $[a,A]$ telle que $i_0=a$ et $i_s=A$ et $g$ est
 {\'e}gale {\`a} une constante $c_r$ sur $]i_{r-1},i_r[$ pour tout $r\in
 \{1,...,s\}$.

Par hypoth{\`e}se $f$ est continue {\`a} gauche en $i_r$ pour tout $r\in \{0,...,s\}$ donc
$|f(i_r)-c_r|=\lim\limits_{t\to i_r-}|f(t)-g(t)|\leq \varepsilon$.

 Si $i_r\notin \Q$, la fonction $f$ est continue en $i_r$ donc il
existe $j_r\in ]i_r,i_{r+1}[\cap \Q$ tel que pour tout $t\in [i_r,j_r]$ on ait
$|f(t)-f(i_r)|<\varepsilon$. Pour $t\in [i_r,j_r]$ on a
$|f(t)-c_r|\leq |f(t)-f(i_r)|+|f(i_r)-c_r|<2\varepsilon$.

Si $i_r\in \Q$, on pose $j_r=i_r$.

Maintenant on a une nouvelle subdivision $j_0,...,j_s$ de
$[a,A]$. Soit $g$ la fonction d{\'e}finie par

 $$g(t)=\left\{\begin{aligned}&0,\;t>A
\\&\lambda,\;t\leq a
\\&c_r,\;t\in
    ]j_{r-1},j_r]\;\mbox{avec}\;r\in
    \{1,...,s\}.\end{aligned}\right.$$

On a $g=\lambda {\mathbf{1}}_{]-\infty,a]}+\sum\limits_{r=1}^s c_r
{\mathbf{1}}_{]j_{r-1},j_r]}=\lambda \tau(e_a)+\sum\limits_{r=1}^s
c_r\tau(e_{j_r}-e_{j_{r-1}})\in \tau(\gA_\cF)$ o{\`u}
$\cF=\{j_0,...,j_s\}\in \F_\Q$. Comme
$||f-g||_{\infty}=\sup\limits_{t\in \R}|f(t)-g(t)|<2\varepsilon$ et
$\tau(\gA)$ est ferm{\'e} on en
deduit que $f\in \tau(\gA)$.\\

 On a $\mbox{Sp}(\gA)=\mbox{Sp}(\tau(\gA))=\mbox{Sp}(\cA)$. Par
 l'exemple pr{\'e}cedent il y a une bijection entre les caract{\`e}res de
 $\gA$ et les sous-semi-treillis finissants non vides $\cM$ de 
 $\Q$. Remarquons que les seules formes possibles de $\cM$ sont les
 suivantes:
 $\cM=\Q$ ou $\cM=\Q\cap [t,+\infty[\equiv \cL_t$ ou
 $\cM=\Q\cap ]t,+\infty[\equiv \cM_t$ avec $t\in \Q$. Remarquons que si $t\in \R\backslash \Q$
 alors $\cL_t=\cM_t$. On a donc une bijection naturelle entre les
 caract{\`e}res de $\gA$ et $\{\Q\}\cup \{\cL_t;\,t\in \R\}\cup
 \{\cM_t;\,t\in \Q\}$.

Id{\'e}ntifions le caract{\`e}re qui correspond {\`a} chacun de ces
sous-semi-treillis finissants {\`a} travers l'isomorphisme $\tau:\gA\to
\cA$ d{\'e}crit ci-dessus. En v{\'e}rifiant qu'il co{\"\i}ncide sur les g{\'e}n{\'e}rateurs
$e_i$ ($i\in \cL$) on trouve : 
\begin{itemize}

\item $\chi_\Q(a)=\lim\limits_{s\to -\infty}\tau(a)(s)$,

\item $\chi_{\cL_t}(a)=\tau(a)(t)$ pour tout $t\in \R$,

\item $\chi_{\cM_t}(a)=\lim\limits_{s\to t+}\tau(a)(s)$ pour tout
  $t\in \Q$.

\end{itemize}

\end{example}

\newpage

\addcontentsline{toc}{chapter}{Bibliographie}

\bibliographystyle{amsalpha}

\end{document}